\magnification=\magstep1 
\input amstex
\UseAMSsymbols        
\hoffset=0truecm \voffset=0truecm \hsize=15truecm
\NoBlackBoxes
   
\input pictex

\mathsurround=1pt
   \newcount\notenumber
   
   \def\note{\advance\notenumber by 1 
       \plainfootnote{$^{\the\notenumber}$}}

\def\Hom{\operatorname{Hom}}
\def\soc{\operatorname{soc}}
\def\top{\operatorname{top}}
\def\rad{\operatorname{rad}}

\def\GL{\operatorname{GL}}
\def\SL{\operatorname{SL}}
\def\SI{\operatorname{SI}}

\def\add{\operatorname{add}}

\def\Ext{\operatorname{Ext}}

\def\mod{\operatorname{mod}}
\def\projdim{\operatorname{proj-dim}}
\def\injdim{\operatorname{inj-dim}}
\def\gldim{\operatorname{gl.dim.}}
\def\bdim{\operatorname{\bold{dim}}}

\def\End{\operatorname{End}}
\def\arr#1#2{\arrow <1.5mm> [0.25,0.75] from #1 to #2}

\def\semi{\;\beginpicture
\setcoordinatesystem units <.2cm,.2cm>
\put{} at 0 0
\put{} at 1 1
\plot 1 1  0 0  0 1  1 0 /
\endpicture\,}

\def\under{\beginpicture
\setcoordinatesystem units <.2cm,.2cm>
\put{} at 0.4 2
\put{} at 1.6 -1
\plot 0.3 2.5  1 2.5  1 -1.6  1.7 -1.6 /
\endpicture}

\def\underkl{\beginpicture
\setcoordinatesystem units <.2cm,.2cm>
\put{} at 0.4 2
\put{} at 1.6 -1
\plot 0.3 1.5  1 1.5  1 -.5  1.7 -.5 /
\endpicture}

\def\T#1{\qquad\text{#1}\qquad}

\def\vecta#1#2#3#4#5#6
  {\ssize {\sssize#1\, #2}{#3 \atop #4}{\sssize #5\,#6}}
\def\vectb#1#2#3#4#5#6
  {\ssize {#1 \atop #2}{\sssize #3\,#4}{#5 \atop #6}}
\def\vectc#1#2#3#4
  {\sssize #1{#2 \atop #3}#4}
\def\vectd#1#2#3#4
  {\ssize #1#2#3}

\centerline{\bf Some Remarks Concerning Tilting Modules and Tilted Algebras.}
	\medskip
\centerline{\bf  Origin. Relevance. Future.}
	\bigskip
\centerline{(An appendix to the Handbook of Tilting Theory)}
	\bigskip
\centerline{Claus Michael Ringel}
	\bigskip
The project to produce a Handbook of Tilting Theory was discussed
during the Frauen{}insel Conference {\it 20 Years of Tilting Theory},
in November 2002. A need was felt to make available
surveys on the basic properties of tilting modules,
tilting complexes and tilting functors, 
to collect outlines of the
relationship to similar constructions in algebra and geometry,
as well as reports on the growing number of generalizations.
At the time the Handbook was conceived, there was a general consensus 
about the overall frame of tilting theory, with the tilted algebra as the core,
surrounded by a lot of additional considerations and with many
applications in algebra and geometry. One was still looking
forward to further generalizations (say something 
like ``pre-semi-tilting procedures
for near-rings''), but the core of tilting theory seemed to be in a
final shape. The Handbook was supposed to provide a full account of the 
theory as it was known at that time.
The editors of this Handbook have to be highly praised for what they have
achieved. But the omissions which were necessary
in order to bound the size of the volume clearly indicate that there should 
be a second volume.
	
Part I will provide an outline of this core of tilting theory. Part II
will then be devoted to topics where tilting
modules and tilted algebras 
have shown to be relevant. I have to apologize that these parts will repeat
some of the considerations of various chapters of the Handbook,
but such a condensed version may be helpful as a sort of guideline. 
Both parts I and II contain historical annodations and reminiscences.
The final Part III will 
be a short report on some striking recent developments which are
motivated by the cluster theory of Fomin and Zelevinsky. In particular, we will guide the reader to the basic properties of cluster tilted algebras, to the
relationship between tilted algebras and cluster tilted algebras, but also
to the cluster categories which provide a universal setting for all the related tilted
and cluster tilted algebras. In addition, we will focus the attention to the
complex of cluster tilting objects and exhibit a quite elementary description
of this complex. In Part I some problems
concerning tilting modules and tilted algebras are raised and one may jump directly
to Part III in order to see in which way these questions have been answered
by the cluster tilting theory. We stress that it should be possible to 
look at the parts II and III independently. 
	\bigskip
\centerline{\bf I.}
	\medskip
The setting to be exhibited is the following:
We start with a hereditary artin algebra $A$
and a tilting $A$-module $T$. It is the endomorphism ring $B = \End(T)$,
called a {\it tilted} algebra, which attracts the attention. The main interest lies
in the comparison of the categories $\mod A$ and $\mod B$ (for any ring $R$,
let us denote by $\mod R$ the category of all $R$-modules of finite
length). We may assume that $A$ is connected (this means that $0$ and $1$
are the only central idempotents), and we may distinguish whether $A$
is representation-finite, tame, or wild; for hereditary algebras, this
distinction is well understood: the corresponding quiver (or better species)
is a Dynkin diagram, a Euclidean diagram, or a wild diagram, respectively.
There is a parallel class of algebras: if we
start with a canonical algebra $A$ instead of a finite-dimensional
hereditary algebra
(or, equivalently, with a weighted projective line, or a so called ``exceptional curve'' 
in the species case),
there is a corresponding tilting procedure.
Again the representation theory distinguishes three different
cases: $A$ may be domestic, tubular, or wild. Now two of the six cases
coincide: the algebras obtained from the domestic canonical algebras via tilting
are precisely those which can be obtained from a Euclidean algebra via tilting. 
Thus, there are 5 possibilities which are best displayed
as the following ``T'': the upper horizontal line refers to the hereditary
artin algebras, the middle vertical line to the canonical algebras.
$$
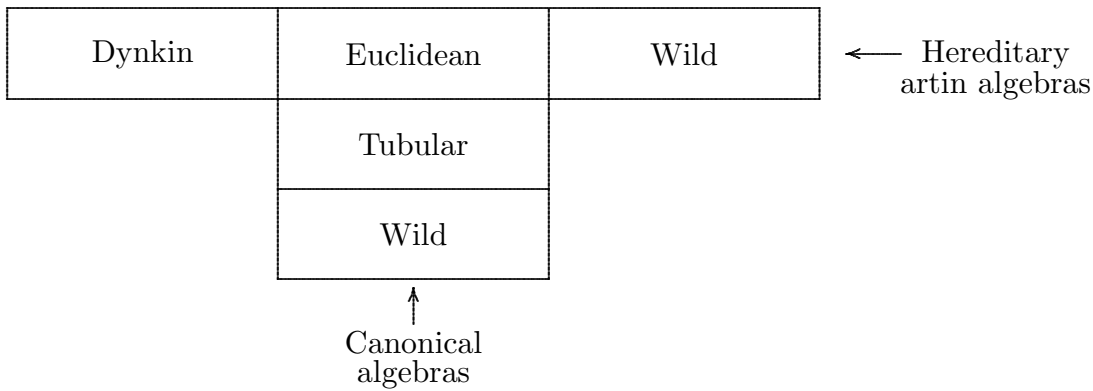
\hbox{\beginpicture
\setcoordinatesystem units <3cm,1cm>
\plot 0 2  0 3  3 3  3 2  0 2 /
\plot 1 3  1 0  2 0  2 3 /
\plot 1 1  2 1 /
\put{} at 0 0
\put{} at 3 3
\put{Dynkin}      at 0.5 2.5
\put{Euclidean}   at 1.5 2.5
\put{Wild} at 2.5 2.5
\put{Tubular}     at 1.5 1.5
\put{Wild}    at 1.5 0.5

\put{Canonical} at 1.5 -0.7
\put{algebras} at 1.5 -1.05
\put{Hereditary} at 3.65 2.5
\put{artin algebras} at 3.65 2.15
\arr{1.5 -0.5}{1.5 -0.1}
\arr{3.3 2.5}{3.1 2.5}

\endpicture}
$$
 
There is a common frame for the five different classes: start with
an artin algebra $A,$ such that the bounded derived category $D^b(\mod A)$
is equivalent to the bounded derived category $D^b(\Cal H)$ of
a hereditary abelian category $\Cal H$. Let $T$ be a tilting
object in $D^b(\mod A)$ and $B$ its endomorphism ring. Then $B$ has
been called a {\it quasi-tilted} algebra by Happel-Reiten-Smal{\o}, and
according to Happel and Happel-Reiten these categories $D^b(\mod A)$
are  just the derived categories of artin algebras which are hereditary
or canonical. In the ``T'' displayed above, the upper horizontal line concerns 
the derived categories with a slice,
the middle vertical line those with a separating tubular family. More
information can be found in Chapter 10 by Lenzing. 

Most of the further considerations will be formulated 
for tilted algebras only. However usually 
there do exist corresponding results for all the quasi-tilted algebras.
To restrict the attention to the tilted algebras has to be seen as
an expression just of laziness, and does not correspond to the high
esteem which I have for the remaining algebras (and the class
of quasi-tilted algebras in general). 

Thus, let us fix again a hereditary artin algebra $A$ and let $D$ be the 
standard duality of $\mod A$ (if $k$ is the center of $A$, then $D = \Hom_k(-,k)$;
note that $k$ is semisimple). Thus $DA$ is an injective cogenerator in $\mod A$. 
We consider a tilting
$A$-module $T$, and let $B = \End(T).$ The first feature which comes to mind
and which was the observation by Brenner and Butler which started the
game\note{Of course, we are aware that examples of tilting modules and tilting
functors had been studied before. Examples to be mentioned are first the
Coxeter functors introduced by Gelfand and Ponomarev in their paper on the
four-subspace-problem (1970), then the BGP-reflection functors (Bernstein, Gelfand
and Ponomarev, 1973), their generalisation by Auslander, Platzeck and Reiten (1979),
now called the APR-tilting functors,
and also a lot of additional ad-hoc constructions used around the globe,
all of which turn out to be special tilting functors. But the proper start of tilting
theory is clearly the Brenner-Butler paper (1980). The axiomatic approach of
Brenner and Butler was considered at that time as quite 
unusual and surprising in a theory which still was in an experimental stage. 
But it soon turned out to be a milestone in the development of representation theory.}, 
is the following: the functor $\Hom_A(T,-)$ yields an equivalence
between the category\note{Subcategories like $\Cal T$ and $\Cal F$ will play a role
everywhere in this appendix. In case we want to stress that they are defined using
the tilting module $T$, we will write $\Cal T(T)$ instead of $\Cal T$, and so on.}
 $\Cal T$ of all $A$-modules generated by $T$
and the category $\Cal Y$ of all $B$-modules cogenerated by the $B$-module $\Hom_A(T,DA)$. 
Now the dimension vectors of the
indecomposable $A$-modules in $\Cal T$ generate the Grothendieck group
$K_0(A)$. If one tries to use $\Hom_A(T,-)$ in order to identify the
Grothendieck groups $K_0(A)$ and $K_0(B)$, one observes that the positivity
cones overlap, but differ: the new axes which define the positive cone for $B$
are ``tilted'' against those for $A$. This was the reason for Brenner and
Butler to call it a tilting procedure. But there is a second ``tilting'' phenomenon
which concerns the corresponding torsion 
pairs\note{In contrast
to the usual convention in dealing with a torsion pair or a ``torsion theory'', 
we name first the
torsion-free class, then the torsion class: this fits to the rule that in a rough
thought, maps go from left to right, and a torsion pair concerns regions with ``no maps
backwards''.}. 
In order to introduce these torsion pairs, we have to look not only at
the functor $\Hom_A(T,-)$, but also at 
$\Ext^1_A(T,-)$. 
The latter functor yields an equivalence between the category $\Cal F$
of all $A$-modules $M$ with $\Hom_A(T,M) = 0$ and the category $\Cal X$
of all $B$-modules $N$ with $T\otimes_B N = 0$. Now the pair $(\Cal F,\Cal T)$
is a torsion pair in the category of $A$-modules, and the pair $(\Cal Y,\Cal X)$
is a torsion pair in the category of $B$-modules: 
$$
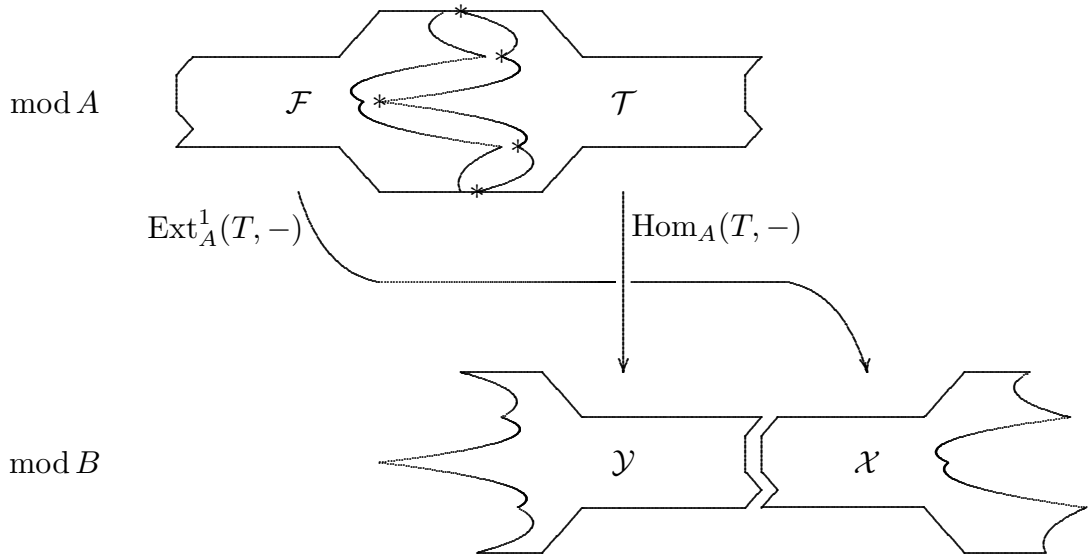
\hbox{\beginpicture
\setcoordinatesystem units <.9cm,1cm>
\put{} at 1 0
\put{} at 13.7 6

\plot 2.5 4.5  2.7 4.7  2.5 4.9  2.5 5.3  2.7 5.5 /
\plot 2.5 4.5  4.5 4.5  5 4  7 4  7.5 4.5  9.5 4.5 / 
\plot 9.5 4.5  9.7 4.7  9.5 4.9  9.5 5.3  9.7 5.5 /
\plot 2.7 5.5  4.5 5.5  5 6  7 6  7.5 5.5  9.7 5.5 /

\plot 9.5 0.5  9.7 0.7  9.5 0.9  9.5 1.3  9.7 1.5 /
\plot 9.7 0.5  9.9 0.7  9.7 0.9  9.7 1.3  9.9 1.5 /
\plot 6 2  7 2  7.5 1.5  9.7 1.5 /
\plot 6.2 0  7 0  7.5 0.5  9.5 0.5 /

\put{$*$} at 6 6
\put{$*$} at 6.5 5.5
\put{$*$} at 5 5
\put{$*$} at 6.7 4.5
\put{$*$} at 6.2 4

\setquadratic
\plot 6 6      6.65  5.75     6.5 5.5  /
\plot 6.5 5.5  6.55  5.25    5 5  /
\plot 5 5      6.55  4.75   6.7 4.5 /
\plot 6.7 4.5  6.85  4.25   6.2 4 /

\plot 5.8 6    5.65 5.75     6.3 5.5  /
\plot 6.3 5.5  4.85 5.25  4.8 5  /
\plot 4.8 5    5.05 4.75    6.5 4.5 /
\plot 6.5 4.5  6.05 4.25   6 4 /

\plot 6 2      6.65  1.75     6.5 1.5  /
\plot 6.5 1.5  6.55  1.25    5 1  /
\plot 5 1    6.55  0.75   6.7 0.5 /
\plot 6.7 0.5  6.85  0.25   6.2 0 /

\plot 13 2    12.85 1.75     13.5 1.5  /
\plot 13.5 1.5  12.05 1.25  12 1  /
\plot 12 1    12.25 0.75    13.7 0.5 /
\plot 13.7 0.5  13.25 0.25   13.2 0 /

\setlinear
\plot 9.9 1.5  11.7 1.5  12.2 2  13   2 / 
\plot 9.7 0.5  11.7 0.5  12.2 0  13.2 0 /
\put{$\mod A$} at 1 5
\put{$\mod B$} at 1 1

\put{$\Cal F$} at 4 5
\put{$\Cal T$} at 8 5
\put{$\Cal Y$} at 8 1
\put{$\Cal X$} at 11 1

\arr{8 4}{8 2}
\arr{10.98 2.1}{11 2}
\put{$\Hom_A(T,-)$} at 9.15  3.6 
\put{$\Ext^1_A(T,-)$} at 3.1 3.6

\setquadratic
\plot 4 4  4.4 3.3  5 3 7.5 3  7.7 3 /

\plot 8.1 3  9 3  10 3  10.6 2.7  11 2 /

\endpicture}
$$
and one encounters the amazing fact that under the pair of functors $\Hom_A(T,-)$
and $\Ext^1_A(T,-)$ the torsion-free class of a torsion pair is 
flipped over
the torsion class in order to form a new torsion pair in reversed 
order\note{The discovery of this phenomenon was based on a detailed
examination of many examples (and contributions by Dieter Vossieck, 
then a student at Bielefeld, should be acknowledged).
At that time only the equivalence of $\Cal T(T)$ and $\Cal Y(T)$ was well
understood. The obvious question was to relate the remaining
indecomposable $B$-modules (those in $\Cal X(T)$) to suitable
$A$-modules. 
As Dieter Happel recalls, the first examples leading to a full undestanding
of
the whole tilting process were tilting modules for the $E_6$-quiver with
subspace orientation.}.
The stars
$\ast$ indicate a possible distribution of the indecomposable direct summands $T_i$ of $T$,
and one should keep in mind that for any $i$, the Auslander-Reiten translate  $\tau T_i$ of
$T_i$ belongs to $\Cal T$ (though it may be zero). We have said that the modules in $\Cal T$ are those generated by $T$, but similarly the modules in $\Cal F$ are those cogenerated by $\tau T$.

According to 
Happel\note{When he propagated this in 1984, it was the first clue
that the use of derived categories may be of interest when dealing with question in the
representation theory of finite dimensional algebras. The derived categories had been
introduced by Grothendieck in order to construct derived functors when dealing with
abelian categories which have neither sufficiently many projective nor sufficiently
many injective objects, and at that time they were considered as useless in case there are
enough projectives and enough injectives, as in the cases $\mod A$ and $\mod B$.}, 
the category $\mod B$ should be seen as being embedded into the
derived category $D^b(\mod A)$
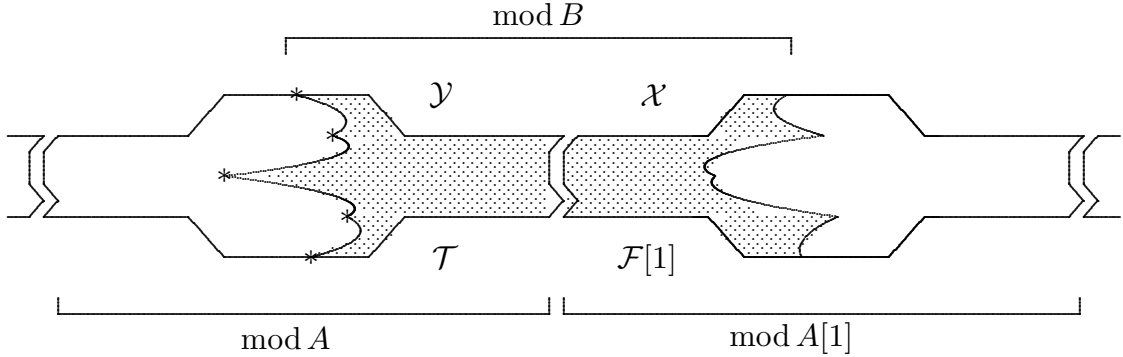
$$
\hbox{\beginpicture
\setcoordinatesystem units <.8cm,.9cm>
\put{\beginpicture 
\plot 2.5 4.5  2.7 4.7  2.5 4.9  2.5 5.3  2.7 5.5 /
\plot 2.5 4.5  4.5 4.5  5 4  7 4  7.5 4.5  9.5 4.5 / 
\plot 9.5 4.5  9.7 4.7  9.5 4.9  9.5 5.3  9.7 5.5 /
\plot 2.7 5.5  4.5 5.5  5 6  7 6  7.5 5.5  9.7 5.5 /

\plot 2 4.5 2.3 4.5  2.5 4.7  2.3 4.9  2.3 5.3  2.5 5.5 2 5.5 /

\put{$*$} at 6 6
\put{$*$} at 6.5 5.5
\put{$*$} at 5 5
\put{$*$} at 6.7 4.5
\put{$*$} at 6.2 4

\setshadegrid span <.5mm>
\hshade 4    6   7   <,,,z>
        4.3  7   7.4 <,,z,z>
        4.5  6.7 7.4 <,,z,z>
        4.55 7   9.5 <,,z,z>
        5    5.5 9.5 <,,z,z>
        5.45 7 9.5 <,,z,z>
        5.5  6.5 7.4 <,,z,z>
        5.7  6.8 7.4 <,,z,>
        6    6   7  /

\setquadratic
\plot 6 6      6.65  5.75   6.5 5.5  /
\plot 6.5 5.5  6.55  5.25   5.0 5  /
\plot 5 5      6.55  4.75   6.7 4.5 /
\plot 6.7 4.5  6.85  4.25   6.2 4 /

\endpicture} at 0 0

\put{\beginpicture 

\plot 2.5 4.5  2.7 4.7  2.5 4.9  2.5 5.3  2.7 5.5 /
\plot 2.5 4.5  4.5 4.5  5 4  7 4  7.5 4.5  9.5 4.5 / 
\plot 9.5 4.5  9.7 4.7  9.5 4.9  9.5 5.3  9.7 5.5 /
\plot 2.7 5.5  4.5 5.5  5 6  7 6  7.5 5.5  9.7 5.5 /
\plot 4.3 4.5  4.5 4.5  5 4  7 4  7.5 4.5  7.7 4.5 / 
\plot 4.3 5.5  4.5 5.5  5 6  7 6  7.5 5.5  7.7 5.5 /

\plot 10.2 4.5  9.7 4.5  9.9 4.7  9.7 4.9  9.7 5.3  9.9 5.5  10.2 5.5 /

\put{} at 6 6
\put{} at 6.5 5.5
\put{} at 5 5
\put{} at 6.7 4.5
\put{} at 6.2 4

\setshadegrid span <.5mm>
\hshade 4   5   5.8  <,,,z>
        4.3 4.6 6  <,,z,z>
        4.5 4.7 6.3  <,,z,z>
        4.55 2.5 5.2  <,,z,z>
        5   2.5 4.5  <,,z,z>
        5.45 2.5 5  <,,z,z>
        5.5 4.7 6  <,,z,z>
        5.7 4.6 5.4  <,,z,>
        6   5   5.6   /

\setquadratic

\plot 5.6 6    5.45 5.75    6.1 5.5  /
\plot 6.1 5.5  4.65 5.25    4.6 5  /
\plot 4.6 5    4.85 4.75    6.3 4.5 /
\plot 6.3 4.5  5.85 4.25    5.8 4 /
\endpicture} at 7.2 0

\plot   -3.15 -1.5   -3.15 -1.7  3.65 -1.7  3.65 -1.5 /
\put{$\mod A$} at 0 -2 

\plot   3.85 -1.5   3.85 -1.7  11 -1.7  11 -1.5 /
\put{$\mod A[1]$} at 7 -2 

\plot   0 1.5   0 1.7  7 1.7  7 1.5 /
\put{$\mod B$} at 3.5 2

\put{$\Cal Y$} at 2.15 1
\put{$\Cal X$} at 5.1 1
\put{$\Cal T$} at 2.2 -1
\put{$\Cal F[1]$} at 5 -1.03
\endpicture}
$$
Under this embedding, $\Cal Y = \Cal T$ is the intersection of $\mod B$
with $\mod A$, whereas $\Cal X = \Cal F[1]$ is the intersection of $\mod B$ with
$\mod A[1]$ (the shift functor in a triangulated category will always be denoted by $[1]$). 
This embedding functor $\mod B \to D^b(\mod A)$ extends to an equivalence of $D^b(\mod B)$ and
$D^b(\mod A)$, and this equivalence is one of the essential features of tilting theory.

Looking at the torsion pair $(\Cal Y,\Cal X)$, there is a sort of asymmetry due to the
fact that $\Cal Y$ is always sincere (this means that every simple module occurs as a 
composition factor of some module in $\Cal Y$), 
whereas $\Cal X$ does not have to be sincere (this happens if $\Cal Y$ contains 
an indecomposable injective module). 
As a remedy, one should divide $\Cal Y$ further
as follows: $\Cal Y$ contains the slice module
$S = \Hom_A(T,DA),$  let $\Cal S = \add S$, and denote by $\Cal Y'$ the class
of all $B$-modules in $\Cal Y$ without an indecomposable direct summand in $\Cal S$.
It is the triple $(\Cal Y',\Cal S,\Cal X)$, which really should be kept in mind:
$$
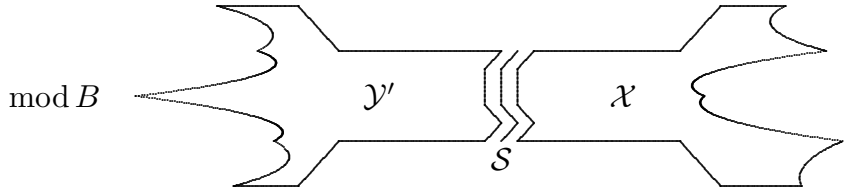
\hbox{\beginpicture
\setcoordinatesystem units <.9cm,1cm>
\put{} at 6 0
\put{} at 13.7 2

\plot 9.3 0.5  9.5 0.7  9.3 0.9  9.3 1.3  9.5 1.5 /
\plot 9.5 0.5  9.7 0.7  9.5 0.9  9.5 1.3  9.7 1.5 /
\plot 9.7 0.5  9.9 0.7  9.7 0.9  9.7 1.3  9.9 1.5 /
\plot 6   2  7 2  7.5 1.5  9.5 1.5 /
\plot 6.2 0  7 0  7.5 0.5  9.3 0.5 /

\setquadratic
\plot 6 2      6.65  1.75     6.5 1.5  /
\plot 6.5 1.5  6.55  1.25    5 1  /
\plot 5 1    6.55  0.75   6.7 0.5 /
\plot 6.7 0.5  6.85  0.25   6.2 0 /

\plot 13 2    12.85 1.75     13.5 1.5  /
\plot 13.5 1.5  12.05 1.25  12 1  /
\plot 12 1    12.25 0.75    13.7 0.5 /
\plot 13.7 0.5  13.25 0.25   13.2 0 /

\setlinear
\plot 9.9 1.5  11.7 1.5  12.2 2  13   2 / 
\plot 9.7 0.5  11.7 0.5  12.2 0  13.2 0 /
\put{$\mod B$} at 4 1

\put{$\Cal Y'$} at 8 1
\put{$\Cal S$} at 9.5 0.3
\put{$\Cal X$} at 11 1

\endpicture}
$$
with all the indecomposables lying in one of the classes $\Cal Y', \Cal S, \Cal X$
and with no maps backwards (the only maps from $\Cal S$ to $\Cal Y'$, from $\Cal X$
to $\Cal S$, as well as from $\Cal X$ to $\Cal Y$ are the zero maps). 
Also note that any indecomposable
projective $B$-module belongs to $\Cal Y$ or $\Cal S$, any indecomposable injective
module to $\Cal S$ or $\Cal X.$ The module class $\Cal S$ is a slice
(as explained in Chapter 3 by Br\"ustle) and any slice is obtained in this way.
The modules in $\Cal Y'$ are those cogenerated by $\tau \Cal S$, the modules in
$\Cal X$ are those generated by $\tau^{-1}\Cal S$. 
	\medskip
Here is an example. Start with 
the path algebra $A$ of a quiver of Euclidean type $\widetilde {\Bbb A}_{22}$ having one
sink and one source. Let $B = \End(T)$, where 
$T$ is the direct sum of the simple projective, the simple injective
and the two indecomposable regular modules of length 3 (this is a tilting module),
then the quiver of $B$ is the same as the quiver of $A$, but $B$ is an algebra with
radical square zero. Thus $B$ is given by a square with two zero relations. 
$$
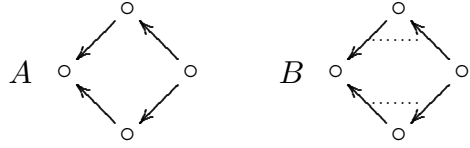
\hbox{\beginpicture
\setcoordinatesystem units <1cm,1cm>
\put{\beginpicture
\setcoordinatesystem units <.7cm,.7cm>
\multiput{$\circ$} at 0 1  1 0  1 2  2 1 /
\arr{0.8 1.8}{0.2 1.2}
\arr{0.8 0.2}{0.2 0.8}
\arr{1.8 1.2}{1.2 1.8}
\arr{1.8 0.8}{1.2 0.2}
\put{$A$} at -0.7 1
\endpicture} at -3 0
 
\put{\beginpicture
\setcoordinatesystem units <.7cm,.7cm>
\multiput{$\circ$} at 0 1  1 0  1 2  2 1 /
\arr{0.8 1.8}{0.2 1.2}
\arr{0.8 0.2}{0.2 0.8}
\arr{1.8 1.2}{1.2 1.8}
\arr{1.8 0.8}{1.2 0.2}
\setdots <.8mm>
\plot 0.6 1.5  1.4 1.5 /
\plot 0.6 0.5  1.4 0.5 /
\put{$B$} at -0.7 1
\endpicture} at 0 0 
\endpicture}
$$
The category $\mod B$ looks as follows: 
$$
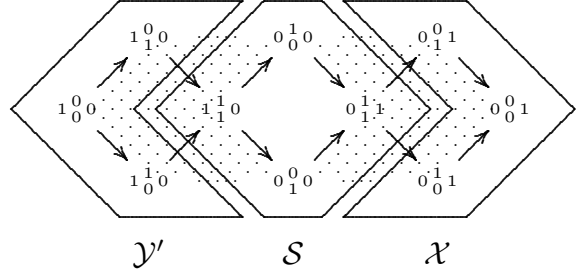
\hbox{\beginpicture
\setcoordinatesystem units <0.8cm,0.8cm>
\put{} at 0 0
\put{} at 0 3
\put{$\vectc 1100 $} at 1 1
\put{$\vectc 0010 $} at 3 1  
\put{$\vectc 0101 $} at 5 1

\put{$\vectc 1000 $} at 0 2
\put{$\vectc 1110 $} at 2 2 
\put{$\vectc 0111 $} at 4 2
\put{$\vectc 0001 $} at 6 2

\put{$\vectc 1010 $} at 1 3
\put{$\vectc 0100 $} at 3 3
\put{$\vectc 0011 $} at 5 3

\arr{0.3 1.7} {0.7 1.3} 
\arr{1.3 1.3} {1.7 1.7} 
\arr{2.3 1.7} {2.7 1.3} 
\arr{3.3 1.3} {3.7 1.7} 
\arr{4.3 1.7} {4.7 1.3} 
\arr{5.3 1.3} {5.7 1.7} 
\arr{0.3 2.3} {0.7 2.7} 
\arr{1.3 2.7} {1.7 2.3} 
\arr{2.3 2.3} {2.7 2.7} 
\arr{3.3 2.7} {3.7 2.3} 
\arr{4.3 2.3} {4.7 2.7} 
\arr{5.3 2.7} {5.7 2.3} 

\setdots <1mm>
\plot 1.7 1  2.3 1 /
\plot 3.7 1  4.3 1 /

\plot 1.7 3  2.3 3 /
\plot 3.7 3  4.3 3 /

\setsolid
\plot 1.1 2  2.6 3.5  3.4 3.5  4.9 2  3.4 0.5  2.6 0.5  1.1 2 /
\plot 0.8 2  2.3 3.5  .6 3.5  -.9 2  .6 0.5  2.3 0.5  0.8 2 /
\plot 5.2 2  3.7 0.5  5.4 0.5  6.9 2  5.4 3.5  3.7 3.5  5.2 2   /
\setshadegrid span <1mm>
\hshade 1 1 3 <,,,z> 2 0 2 <,,z,> 3 1 3 / 
\hshade 1 3 5 <,,,z> 2 4 6 <,,z,> 3 3 5 / 
\put{$\Cal Y'$} at 1 0
\put{$\Cal S$} at 3 0
\put{$\Cal X$} at 5 0
\endpicture} 
$$

The separation of $\mod B$ into the three classes $\Cal Y', \Cal S, \Cal X$ can be phrased
in the language of {\bf cotorsion pairs.} Cotorsion pairs are very well related to
tilting theory (see the Chapters 7 and 11 by Reiten and Trlifaj), 
but still have to be rated more as a sort of insider
tip. We recall the definition: the pair $(\Cal V, \Cal W)$ of full subcategories 
of $\mod A$ is said to be a {\it cotorsion pair} provided $\Cal V$ is the class of all
$A$-modules $V$ with $\Ext^1_A(V,W) = 0$ for all $W$ in $\Cal W$,
and $\Cal W$ is the class of all $A$-modules $W$
such that $\Ext^1_A(V,W) = 0$ for all $V$ in $\Cal V$. The cotorsion pair is said to be {\it split},
provided every indecomposable $A$-modules belongs to $\Cal V$ or $\Cal W$. Usually some
indecomposables will belong to both classes, they are said to form the {\it heart.} 
In our case the following holds: {\it The pair
$\bigl(\add(\Cal Y',\Cal S),\add(\Cal S,\Cal X)\bigr)$ forms a 
split cotorsion pair with heart $\Cal S$.} 

We also see that the modules in $\Cal Y'$ and in $\Cal S$ have projective dimension at most
1, those in $\Cal S$ and in $\Cal X$ have injective dimension at most 1. As a consequence,
if $X,Y$ are indecomposable modules with $\Ext^2_B(X,Y) \neq 0$, then $X$ belongs to $\Cal X$
and $Y$ belongs to $\Cal Y'$. 

Let me add a remark even if it may be considered to be superfluous --- 
its relevance should become clear in the last part of this appendix. If we feel that the
subcategory $\Cal Y'$ has the same importance as $\Cal X$ (thus that it is of interest), then we
should specify an equivalent subcategory, say $\Cal T'$ of $\mod A$ and an equivalence
$\Cal T' \to \Cal Y'$. Such an equivalence is given by the functor 
$$
 \Hom_A(\tau^{-1}T,-)\: \mod A \to \mod B
$$ or, equivalently, by $\Hom_A(T,\tau -)$, since
$\tau^{-1}$ is left adjoint to $\tau.$ This functor vanishes on $\Cal F$ as well as 
on $T$, and it yields an equivalence between the subcategory $\Cal T'$ of all $A$-modules
generated by $\tau^{-1}T$ and the subcategory $\Cal Y'$ of $\mod B$. Note that the
functor can also be written in the form $D\Ext^1_A(-,T)$, due to the
Auslander-Reiten formula $D\Ext^1(M,T) \simeq \Hom_A(T,\tau M).$ In this way, we see
that we deal with equivalences which are sort of dual to each other:
$$
 D\Ext^1_A(-,T)\: \Cal T' \to \Cal Y' \T{and} 
 \Ext^1_A(T,-)\:\Cal F \to \Cal X.
$$ 
	\medskip
It seems to be worthwhile to have a short look at the rather trivial case
when no modules are lost, so that
the tilting procedure is a kind of rearrangement of module classes.
{\it The following assertions are equivalent:}
\item{(i)} {\it The tilting module $T$ is a slice module.}
\item{(ii)} {\it The endomorphism ring $B$ is hereditary.}
\item{(iii)} {\it The torsion pair $(\Cal F,\Cal T)$ splits.}
\item{(iv)} {\it $\Ext^1_A(\tau T,T) = 0$.}
	\smallskip\noindent
The equivalence of these assertions are well-known, but not too easy to trace.
Some implications are quite obvious, for example that (ii) and (iii) are implied
by (i). Let us show that (ii) implies (i): 
Since $T$ is a tilting module, the $B$-module $T' = \Hom_A(T,DA) \simeq \Hom_k(T,k)$
is a slice module in $\mod B$. Since the $B$-module $T'$ is a tilting module and
$A = \End_B(T'),$ we can use this tilting module in order to tilt from 
$\mod B$ to $\mod A$. Since $B$ is hereditary, we obtain in
$\mod A$ the slice module $\Hom_B(T',DB) \simeq \Hom_k(T',k) \simeq T.$ 
This shows (i). The equivalences of (ii) and (iv), as well as of (iii) and (iv),
can be seen as consequences of more general considerations which will be presented later. 
	\medskip
If $T$ is not a slice module, so that
$B$ has global dimension equal to 2, then the algebras $A$ and $B$ 
play quite a different role: the first difference is of course the fact that
$A$ is hereditary, whereas $B$ is not. Second, there are the two torsion pairs
$(\Cal F,\Cal T)$ in $\mod A$ and $(\Cal Y,\Cal X)$ in $\mod B$ - the second one is
a split torsion pair, the first one not. This means that  we loose modules
going from $\mod A$ to $\mod B$ via tilting. 
Apparently, no one cared about the missing modules,
at least until quite recently. There are two reasons: First of all, we know 
(see Chapter 3 by Br\"ustle), that
the study of indecomposable modules over a representation-finite algebra is reduced
via covering theory to the study of representation-finite tilted algebras.
Such an algebra $B$ may be of the form $B = \End_A(T)$, where $T$ is a tilting $A$-module,
with $A$ representation-infinite. Here we describe the $B$-modules in terms of $A$ and
we are only interested in the finitely many indecomposable $A$-modules which belong to
$\Cal F$ or $\Cal T$, the remaining $A$-modules seem to be of no interest, we do not miss them.
But there is a second reason: the fashionable reference to derived categories is used to appease 
anyone, who still mourns about the missing modules. They are lost indeed as modules, 
but they survive as complexes: since the derived categories of $A$ and $B$ are equivalent, 
corresponding to any indecomposable $A$-module, there is an object in the derived category
which is given by a complex of $B$-modules. However, I have
to admit that I prefer modules to complexes, whenever possible --- thus I was delighted, 
when the lost modules were actually found, as described in Part III of this appendix. 
	\medskip
We will always denote by $n = n(A)$ the number of isomorphism classes of simple $A$-modules.
The interest in tilting $A$-modules directly leads to a corresponding interest in their 
direct summands. These are the modules without self-extensions and are 
called
{\it partial tilting modules.} In particular, one may consider the indecomposable ones: an indecomposable $A$-module without self-extensions is
said to be {\it exceptional} (or a ``stone'',  or  a ``brick without self-extensions'',
or a ``Schurian module without self-extensions''). But there is also an interest in the
partial tilting modules with precisely $n\!-\!1$ isomorphism classes of indecomposable direct
summands, the so-called {\it almost complete} partial tilting modules. If $\overline T$
is an almost complete partial tilting module and $X$ is indecomposable with $\overline T\oplus X$
a tilting module, then $X$ (or its isomorphism class) is called a 
{\it complement} for $\overline T$. It is of interest that any almost complete
partial tilting module $\overline T$ has either 1 or 2 complements, and it has 2 if and only if
$\overline T$ is sincere. Recall that a module is said to be 
{\it basic}, provided it is a direct
sum of pairwise non-isomorphic indecomposable modules. The isomorphism classes of basic
partial tilting modules form a simplicial complex $\Sigma_A$, with vertex set the set
of isomorphism classes of exceptional modules (the vertices of a simplex being its 
indecomposable direct summands), see Chapter 9 by Unger.
Note that this simplicial complex 
is of pure dimension $n-1$. The assertion concerning the complements shows that
it is a pseudomanifold with boundary. The boundary consists of all the non-sincere
almost complete partial tilting modules. 

As an example, consider the path algebra $A$
of the quiver $\circ\leftarrow \circ\leftarrow\circ$. The simplicial complex $\Sigma_A$ 
has the following shape:
$$
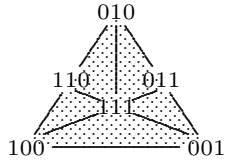
\hbox{\beginpicture
\setcoordinatesystem units <1cm,1.5cm>
\put{} at 1 1
\put{} at 3 2
\put{$\vectd 100 $} at 1 1
\put{$\vectd 110 $} at 1.5 1.5 
\put{$\vectd 010 $} at 2 2 
\put{$\vectd 111 $} at 2 1.3
\put{$\vectd 011 $} at 2.5 1.5
\put{$\vectd 001 $} at 3 1
\plot 1.1 1.1  1.4 1.4 /
\plot 1.6 1.6  1.9 1.9 /
\plot 2.1 1.9  2.4 1.6 /
\plot 2.6 1.4  2.9 1.1 /
\plot 1.3 1  2.7 1 /
\plot 2 1.9  2 1.4 /

\plot 1.55 1.42  1.82 1.35 /
\plot 2.2 1.25  2.8 1.1 /

\plot 2.45 1.42  2.18 1.35 /
\plot 1.8 1.25  1.2 1.1 /

\setshadegrid span <.5mm>
\hshade 1 1 3  2 2 2  / 
\endpicture}
$$

Some questions concerning the simplicial complex $\Sigma_A$ remained open:
What happens under a change of orientation? What happens under a tilting functor?
Is there a way to get rid of the boundary? Here we are again in a situation where 
a remedy is provided by the derived categories: 
If we construct the analogous simplicial complex of tilting complexes in $D^b(\mod A)$,
then one obtains a pseudo manifold without boundary, but this is quite a large 
simplicial complex!

If $\overline T$ is an almost complete partial tilting module, and $X$ and $Y$ are 
non-isomorphic complements for $\overline T$, then either $\Ext^1_A(Y,X) \neq 0$
or $\Ext^1_A(X,Y) \neq 0$ (but not both). If $\Ext^1_A(Y,X) \neq 0$ (what we may
assume), then there exists 
an exact sequence $0 \to X \to T' \to Y \to 0$ with $T' \in \add \overline T$,
and one may write $\overline T \oplus X < \overline T\oplus Y$.
In this way, one gets a partial ordering on the set of isomorphism classes of
basic tilting modules. One may consider the switch between the tilting modules
$\overline T \oplus X$ and $\overline T\oplus Y$ as an exchange process which stops
at the boundary. We will see in Part III that it is possible to define an exchange procedure across the boundary of $\Sigma_A$, and that this can be arranged in such a way that one
obtains an interesting small extension of the simplicial complex $\Sigma_A$.  

As long as reflection functors are defined only for sinks and for sources, this has to be 
considered as a real deficiency, since there is no similar restriction in Lie
theory. Indeed, there the use of reflections for all the vertices is an important tool.
A lot of efforts have been made in representation theory 
in order to overcome this deficiency, see for example
the work of Kac on the dimension vectors of the indecomposable representations of a quiver.  

A reflection functor furnishes a quite small change of the given module category. Let us
consider this in more detail. Let $Q$ be a quiver and $i$ a sink of $Q$. We denote by
$\sigma_iQ$ the quiver obtained from $Q$ by changing the orientation of all arrows ending in $i$,
thus $i$ becomes a source in $\sigma_iQ$.
Let $S(i)$ be the simple $kQ$-module corresponding to the vertex $i$, and
$S'(i)$ the simple $k\sigma_iQ$-module corresponding to the vertex $i$.
The reflection functor $\sigma_i\:\mod kQ \to \mod k\sigma_iQ$ provides an equivalence
between the categories
$$
 \mod kQ/\langle \add S(i)\rangle \to \mod k\sigma_iQ/\langle \add S'(i)\rangle.
$$
In general, given rings $R,R'$ one may look for a simple $R$-module $S$ and a
simple $R'$-module $S'$ such that the categories
$$
 \mod R/\langle \add S\rangle \to \mod R'/\langle \add S'\rangle
$$
are equivalent. In this case, let us say that $R,R'$ are {\it nearly Morita-equivalent}.
As we have seen, for $Q$ a quiver with a sink $i$, the path algebras
$kQ$ and $k\sigma_iQ$ are nearly Morita-equivalent. I am not aware that other pairs of
nearly Morita-equivalent rings have been considered until very recently, but 
Part III will provide a wealth of examples. 
	\medskip
A final question should be raised here. There is a very nice homological characterization
of the quasi-tilted algebras by Happel-Reiten-Smal{\o} [HRS]: these are the artin algebras
of global dimension at most 2, such that any indecomposable module has projective dimension
at most 1 or injective dimension at most 1. But it seems that a corresponding characterization
of the subclass of tilted algebras is still missing. Also, in case we consider tilted
$k$-algebras, where $k$ 
is an algebraically closed field, the possible quivers and their relations are not known.
	\bigskip
\centerline{\bf II.}
	\medskip
The relevance of tilting theory relies on the many different connections it has
not only to other areas of representation theory, but also to algebra and geometry in
general. 
Let me give some indications. If nothing else is said, $A$ will denote a
hereditary artin algebra, $T$ a tilting $A$-module and $B$ its endomorphism ring.
	\medskip
$\bullet$ {\bf Homology.} Already the definition (the vanishing of $\Ext^1$) refers
to homology. We have formulated above that the first feature which comes to mind
is the functor $\Hom_A(T,-).$ But actually all the tilting theory concerns the study
of the corresponding derived functors $\Ext^i_A(T,-)$, or better, of the right derived functor
$R\Hom_A(T,-)$. 

The best setting to deal with these functors are the corresponding derived
categories $D^b(\mod A)$ and $D^b(\mod B)$, they combine to the right derived functor
$R\Hom_A(T,-)$, and this functor is an equivalence, as Happel has shown. 
Tilting modules $T$ in general were defined in such a way that $R\Hom_A(T,-)$ is 
still an equivalence. The culmination of this development was Rickard's 
characterization of rings with equivalent derived categories: such equivalences are
always given by ``tilting complexes''. A detailed account can be found in Chapter 5
by Keller.

Tilting theory can be exhibited well by using spectral sequences. 
In Bongartz's presentation of tilting theory one finds the following formulation:
{\it Well-read mathematicians tend to understand tilting theory 
using spectral-sequences} (which is usually interpreted as a critical comment 
about the earlier papers). But it seems that the first general account 
of this approach is only now
available: the contribution of Brenner and Butler (see Chapter 4) in
this volume. A much earlier one by Vossieck should have been his
Bielefeld Ph.D.~thesis, but he never handed it in.
	\smallskip
$\bullet$ {\bf Geometry and Invariant Theory.} 
The Bielefeld interest in tilting modules
was first not motivated by homological, but by geometrical questions.
Happel's Ph.D.~thesis had focused the attention to quiver representations
with an open orbit (thus to all the partial tilting modules). 
In particular, he showed that the number $s(V)$ of isomorphism classes of indecomposable
direct summands of a representation $V$ with open orbit is bounded by the 
number of simple modules. In this way, the study of open orbits in quiver varieties 
was a (later hidden) step in the development of tilting theory.
When studying open orbits, we are in the setting of what Sato and Kimura
[SK] call {\bf prehomogeneous vector spaces}. On the one hand, the
geometry of the complement of the open orbit is of interest, on the other hand 
one is interested in the structure of the ring of semi-invariants.

Let $k$ be an algebraically closed field and $Q$ a finite quiver 
(with vertex set $Q_0$ and arrow set $Q_1$), and we may assume that $Q$ has no oriented cyclic path, thus the path algebra $kQ$ is just a basic hereditary finite-dimensional $k$-algebra. 
For any arrow $\alpha$ in $Q_1$, denote by $t\alpha$ its tail
and by $h\alpha$ its head, and fix some dimension vector $\bold d$.
Let us consider representations $V$ of $Q$ with a fixed dimension vector $\bold d$, 
we may assume $V(x) = k^{\bold d(x)}$; thus the set of these 
representations forms the affine space
$$
 \Cal R(Q,\bold d) = \bigoplus_{\alpha\in Q_1}\Hom_k(k^{\bold d(t\alpha)},k^{\bold d(h\alpha)}).
$$
The group $\GL(Q,\bold d) = \prod_{x\in Q_0} \GL(\bold d(x))$  
operates on this space via a sort of conjugation, and the orbits under this action
are just the isomorphism classes of representations. One of the results of Happel [H1]
asserts that given a sincere representation $V$ with open orbit, 
then $|Q_0| - s(V)$ is the number of isomorphism classes of representations $W$ with
$\bdim V = \bdim W$ and $\dim \Ext^1_A(W,W) = k$ (in particular, there are
only finitely many such isomorphism classes; we also see that $|Q_0| \ge s(V)$). 

Consider now the ring $\SI(Q,\bold d)$ of semi-invariants on
$\Cal R(Q,\bold d)$; by definition these are the invariants of the subgroup $\SL(Q,\bold d)
= \prod_{x\in Q_0} \SL(\bold d(x))$ of $\GL(Q,\bold d).$
Given two representations $V,W$ of $Q$, one may look at the map:
$$
 d^V_W\:\bigoplus_{x\in Q_0} \Hom_k(V(x),W(x)) \longrightarrow
 \bigoplus_{\alpha\in Q_1} \Hom_k(V(t\alpha),W(h\alpha)),
$$
sending $(f(x))_x$ to $(f(h\alpha)V(\alpha)-W(\alpha)f(t\alpha))_\alpha$. 
Its kernel is just
$\Hom_{kQ}(V,W)$, its cokernel $\Ext^1_{kQ}(V,W).$  In case $d^V_W$ is a square matrix, one may
consider its determinant. According to Schofield [Sc], 
this is a way of producing semi-invariants.
Namely,
the Grothendieck group $K_0(kQ)$ carries a (usually non-symmetric) bilinear form
$\langle-,-\rangle$ with $\langle \bdim V,\bdim W\rangle  = \dim_k 
\Hom_{kQ}(V,W) - 
\dim_k\Ext^1_{kQ}(V,W)$, thus $d^V_W$ is a square matrix if and only if 
$\langle \bdim V,\bdim W \rangle = 0.$ So, if $\bold d\in \Bbb N_0^{Q_0}$ and if
we select a representation $W$ such that $\langle \bold d,\bdim W \rangle = 0,$
then $c_W(V) = \det d^V_W$ yields a semi-invariant $c_W$ in $\SI(Q,\alpha)$. 
Derksen and Weyman (and also Schofield and Van den Bergh)
have shown that these semi-invariants 
form a generating set for $\SI(Q,\bold d)$.
In fact, it is sufficient to consider only indecomposable
representations $W$, thus exceptional $kQ$-modules.
	\smallskip
$\bullet$ {\bf Lie Theory.} It is a well-accepted fact that the representation theory
of hereditary artinian rings has a strong relation to Lie algebras and quantum groups 
(actually one should say: a strong relation to Lie algebras via quantum groups).
Such a relationship was first observed by Gabriel when he showed that the
representation-finite connected quivers are just those with underlying graph being
of the form $\Bbb A_n,\Bbb D_n,\Bbb E_6,\Bbb E_7,\Bbb E_8$ and that in these cases the
indecomposables correspond bijectively to the positive roots. According to Kac,
this extends to arbitrary finite quivers without oriented cycles: 
the dimension vectors of the indecomposable representations 
are just the positive roots of the corresponding Kac-Moody Lie-algebra. It is now known
that it is even possible to reconstruct these Lie algebra using the representation theory of hereditary artin algebras, via Hall algebras. Here one encounters 
the problem of specifying the subring of a Hall algebra, generated by the simple
modules. Schofield induction (to be discussed below) shows that all the exceptional modules belong to this subring. 

It seems to be appropriate to discuss the role of the necessary choices.
Let me start with a 
semisimple finite-dimensional complex Lie-algebra. The choice of a Cartan
subalgebra yields its root system, the choice of a root basis yields a triangular
decomposition (and the set of positive roots). Finally, the choice of a total
ordering of the root basis (or at least of an orientation of the edges)
allows to work with a Coxeter element in the Weyl group, to define a Borel
subalgebra and the corresponding category $\Cal O$. 
Of course, one knows that all these choices are inessential. The situation is
more subtle if we
deal with arbitrary Kac-Moody Lie-algebras: different orderings of the root basis
may yield Coxeter elements which are not conjugate -- the first case is 
$\widetilde {\Bbb A}_3$, where one has to distinguish between 
$\widetilde {\Bbb A}_{3,1}$ and
$\widetilde {\Bbb A}_{2,2}$. On the other hand, when we start
with a representation-finite hereditary artin algebra $A$, no choice at all is needed
in order to write down its Dynkin diagram: it is intrinsically given as the 
$\Ext$-quiver of the simple $A$-modules,
and we obtain in this way a Dynkin diagram with orientation.
A change of orientation corresponds to module categories with quite distinct
properties (as already the algebras of type $\Bbb A_3$ show). This difference is still
preserved when one looks at the corresponding Hall algebras, and it comes as
a big surprize that only a small twist of its multiplication 
is needed in order to get an algebra which is independent of the orientation.
	\smallskip
$\bullet$ {\bf The Combinatorics of Root Systems.} It is necessary to dig deeper
into root systems since they play an important role for dealing with 
$A$-modules. Of interest is the corresponding quadratic form,
and the reflections which preserve the root system (but not necessarily 
the positivity of roots), and their compositions, in particular the Coxeter transformations and the BGP-reflection functors.
Unfortunately, these reflection functors are only defined for sinks and for sources!
We will return to the reflection functors when we deal with generalizations of
Morita equivalences, but also in Part III. 

As we have mentioned the relationship between the representation theory of a hereditary
artin algebra $A$ and root systems is furnished by the dimension vectors: We consider
the Grothendieck group $K_0(A)$ (of all finite length $A$-modules modulo exact
sequences). Given an $A$-module $M$, we denote by $\bdim M$ the corresponding 
element in $K_0(A)$; this is what is called its {\it dimension vector}. The 
dimension vectors of the indecomposable $A$-modules are the positive roots of 
the root system in question.
A positive root $\bold d$ is said to be a {\it Schur root}
provided there exists an indecomposable $A$-module $M$ with $\bdim M = \bold d$
and $\End_A(M)$ a division ring. 
The dimension vectors of the exceptional modules are Schur roots, they
are just the real (or Weyl) Schur roots. In case $A$ is representation-finite
then all the positive roots are Schur roots, also for $n(A) = 2$ all the real
roots are Schur roots. But in all other cases, the set of real Schur roots depends
on the choice of orientation. For example, consider the following three orientations
of $\widetilde D_4$:
$$
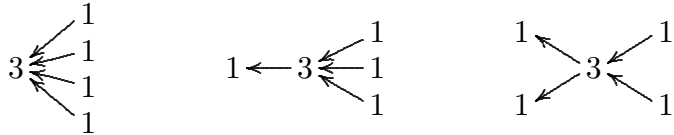
\hbox{\beginpicture
\setcoordinatesystem units <.8cm,.4cm>
\put{\beginpicture
\multiput{$1$} at 1 3  1 2  1 1  1 0 /
\put{$3$} at 0 1.5
\arr{0.8 0.2}{0.2 1.2}
\arr{0.8 1.1}{0.2 1.4}
\arr{0.8 1.9}{0.2 1.6}
\arr{0.8 2.8}{0.2 1.8}
\endpicture} at 0.5 0

\put{\beginpicture
\multiput{$1$} at 0 1  2 0  2 1  2 2 /
\put{$3$} at 1 1
\arr{0.8 1}{0.2 1}
\arr{1.8 1.8}{1.2 1.2}
\arr{1.8 1}{1.2 1}
\arr{1.8 0.1}{1.2 0.8}
\endpicture} at 4 0

\put{\beginpicture
\multiput{$1$} at 0 0  0 2  2 0  2 2 /
\put{$3$} at 1 1
\arr{0.8 1.15}{0.2 1.9}
\arr{0.8 0.85}{0.2 0.1}
\arr{1.8 1.9}{1.2 1.15}
\arr{1.8 0.1}{1.2 0.85}
\endpicture} at 8 0

\endpicture} 
$$
The two dimension vectors on the left are Schur roots, whereas the right one
is not a Schur root.

In order to present the dimension vectors of the indecomposable $A$-modules, one may depict the
Grothendieck group $K_0(A)$; a very convenient way seems to be to work with homogeneous
coordinates, say with the projective space of 
$K_0(A)\otimes_{\Bbb Z} \Bbb R$. It is the merit of Derksen and Weyman [DW3] of having 
popularized this presentation well: they managed to get it to a cover of the Notices of
the American Mathematical Society [DW2]. One such example has been shown in Part I, when we
presented the simplicial complex $\Sigma_A$, whith $A$ the path algebra of the linearly 
oriented quiver of type $\Bbb A_3$. In general, dealing with  hereditary $A$, one is interested
in the position of the Schur roots. Our main concern are 
the real Schur roots as the dimension vectors $\bdim E$ 
of the exceptional modules $E$. They are best
presented by marking the corresponding {\it exceptional lines:} look for
orthogonal exceptional pairs $E_1,E_2$ (this means: $E_1, E_2$ are exceptional modules with
 $\Hom_A(E_1,E_2) = \Hom_A(E_2,E_1) = \Ext_A^1(E_2,E_1) = 0$) and draw the line
segment from $\bdim E_1$ to $\bdim E_2$. The discussion of Schofield induction below will explain the importance of these exceptional lines.

As we know, a tilting $A$-module $T$ has precisely $n = n(A)$ isomorphism classes
of indecomposable direct summands, say $T_1,\dots,  T_n$ and the dimension vectors 
$\bdim T_1,\dots, \bdim T_n$ are linearly independent. We may consider the cone 
$C(T)$ in $K_0(A)\otimes \Bbb R$ generated by $\bdim T_1,\dots, \bdim T_n$. These
cones are of special interest. Namely, 
if a dimension vector $\bold d$ belongs to $C(T)$, then there is a unique 
isomorphism class of modules $M$ with $\bdim M = \bold d$ such that
$\End_A(M)$ is of minimal dimension, and such a module $M$ has no self-extensions. 
On the other hand, the dimension vector of any module
$M$ without self-extensions lies in such a cone (since these modules, the 
partial tilting modules, are the direct summands of tilting modules). The set of these cones forms a fan as they are considered in toric geometry (see for example
the books of Fulton and Oda).

As Hille [Hi] has pointed out, one should use the geometry of these cones in order to
introduce the following notion: If $T = \bigoplus_{i=1}^n T_i$ is a basic
tilting module with indecomposable modules $T_i$ of length $|T_i|$,
he calls $\prod_{i=1}^n |T_i|^{-1}$ the {\it volume}
of $T$. It follows that 
$$
 \sum\nolimits_T v(T) \le 1 
$$
(where the summation extends over all
isomorphism classes of basic tilting $A$-modules), 
with equality if and only if
$A$ is representation-finite or tame. This yields interesting equalities:
For example, the preprojective tilting modules of the Kronecker quiver yield
$$
       \frac11\cdot\frac13 + \frac13\cdot\frac15 + \frac15\cdot\frac17 + \cdots = \frac12.
$$
One may refine these considerations by replacing the length $|T_i|$ by the
$k$-dimension of $T_i$, at the same time replacing 
$A$ by all the Morita equivalent algebras. In this
way one produces power series identities in $n(A)$ variables which should be of 
general interest, for example
$$
       \frac1{x}\cdot\frac1{2x+y} + \frac1{2x+y}\cdot\frac1{3x+2y} + 
       \frac1{3x+2y}\cdot\frac1{4x+3y} + \cdots = \frac1{2xy}.
$$
Namely, let $x$ and $y$ denote the $k$-dimension of the simple projective, or simple
injective $A$-module, respectively. Then the indecomposable preprojective $A$-module
$P_t$ of length $2t\!-\!1$ is of dimension $tx+(t-1)y$.  This means that the
the tilting module $P_t\oplus P_{t+1}$ contributes the summand 
$\frac1{tx+(t-1)y)}\cdot \frac1{(t+1)x+ty}.$ The sequence of tilting modules
$P_1\oplus P_2,\ P_2\oplus P_3,\dots$ yields the various summands on the left side.
	\medskip
$\bullet$ {\bf Combinatorial Structure of Modules.} A lot of tilting theory is
devoted to combinatorial considerations. The combinatorial invariants just discussed 
concern the Grothendieck group. But also the exceptional modules themselves have
a combinatorial flavor: they are ``tree modules'' [R4]. 
As we have mentioned, the orbit of a tilting
module is open in the corresponding module variety, and this holds true with respect
to all the usual topologies, in particular, the Zariski topology, but also the usual
real topology in case the base field is $\Bbb R$ or $\Bbb C$. This means that a slight
change of the coefficients in any realization of $T$ using matrices will not change
the isomorphism class. Now in general to be able to change the coefficients slightly,
will not allow to prescribe a finite set (for example $\{0,1\}$)  of coefficients 
which one may like to use: the corresponding matrices may just belong to the complement
of the orbit. However, in case we deal with the path algebra of a quiver, the exceptional
modules have this nice property: there always exists a realization of $E$ using matrices with
coefficients only $0$ and $1$. A stronger statement holds true: If $E$ has dimension $d$,
then there is a matrix realization which uses precisely $d-1$ coefficients equal to $1$,
and all the remaining ones are $0$ (note that in order to be indecomposable, we need at least
$d-1$ non-zero coefficients; thus we assert that really the minimal possible number of
non-zero coefficients can be achieved). 
	\smallskip
$\bullet$ {\bf Numerical Linear Algebra.} Here we refer to the previous consideration:
The relevance of $0$-$1$-matrices in numerical linear algebra is well-known. 
Thus linear algebra problems, which
can be rewritten as dealing with partial tilting modules, are very suitable for
numerical algorithms, because of two reasons: one can restrict to $0$-$1$-matrices
and the matrices to be considered involve only very few non-zero entries.
	\smallskip
$\bullet$ {\bf Module Theory.} Of course, tilting theory is part of module theory. 
It provides a very useful collection of non-trivial examples for many 
central notions in ring and module theory. The importance of modules without self-extensions
has been realized a long time ago, for example one may refer to the lecture notes of Tachikawa
from 1973. Different names are in use for such modules such as ``splitters''.

It seems that the tilting theory exhibited for the first time a wide range of 
torsion pairs, with many different features: there are 
the splitting torsion pairs, which one finds in the module category of any tilted algebra,
as well as the various non-split torsion pairs in the category $\mod A$ itself.
As we have mentioned in Part I, tilting theory also gives rise to non-trivial 
examples of cotorsion pairs.
And there are corresponding approximations, but also
filtrations with prescribed factors. Questions concerning subcategories of module
categories are considered in many of the contributions in this Handbook, 
in particular in Chapter 7 by Reiten,
but also in the Chapters 8, 11 and 12 by Donkin, Trlifaj and Solberg, respectively.

We also should mention the use of {\bf perpendicular categories}. Starting
with an exceptional $A$-module $E$, the category $E^\perp$ of all $A$-modules
$M$ with $\Hom_A(E,M) = 0 = \Ext^1_A(E,M)$ is again a module category, say
$E^\perp \simeq \mod A'$, where $A'$ is a hereditary artin algebra with $n(A') = 
n(A) - 1.$ These perpendicular categories are an important tool for inductive
arguments and they can be considered as a kind of localisation.

Another notion should be illuminated here: recall that a left
$R$-module $M$ is said to have the {\bf double centralizer property} (or to be balanced), 
provided the following holds: If we denote by $S$ the endomorphism ring of ${}_RM$, say
operating on the right on $M$, we obtain a right $S$-module $M_S$, and we may now look
at the endomorphism ring $R'$ of $M_S$. Clearly, there is a canonical 
ring homomorphism  $R \to R'$ (sending $r\in R$ to the left multiplication by $r$ on $M$),
and now we require that this map is surjective (in case $M$ is a faithful $R$-module,
so that the map $R \to R'$ is injective, this means that we can identify $R$ and $R'$:
the ring $R$ is determined by the categorical properties of $M$, namely its endomorphism ring 
$S$, and the operation of $S$ on the underlying abelian group of $M$). 
Modules with the double centralizer property 
are very important in ring and module theory. Tilting modules satisfy the double centralizer property and this is used in many different ways. 

Of special interest is also the following {\bf subquotient realization} of $\mod A$. All the modules in $\Cal T$
are generated by $T$, all the modules in $\Cal F$ are cogenerated by $\tau T$.
It follows that {\it for any $A$-module $M$, there exists an $A$-module $X$ with
submodules $X''' \subseteq X'' \subseteq X'$ such that $X''$ is a direct sum
of copies of $T$, whereas $X/X''$ is a direct sum of copies of $\tau T$ and such that 
$M = X'/X'''$} (it then follows that $X''/X'''$ is the torsion submodule of $M$
and $X'/X''$ its torsion-free factor module). 
In particular, we see that $(\Cal F,\Cal T)$ is a split torsion pair
if and only if $\Ext^1_A(\tau T,T) = 0$ (the equivalence of conditions 
(iii) und (iv) in Part I). 
This is one of the results which stresses the importance of
the bimodule $\Ext^1_A(\tau T,T).$ Note that the extensions considered when we look at
$\Ext^1_A(\tau T,T)$ are 
opposite to those the Auslander-Reiten translation $\tau$ is famous for (namely the
Auslander-Reiten sequences, they correspond to elements of $\Ext^1_A(T_i,\tau T_i)$, 
where $T_i$ is
a non-projective indecomposable direct summand of $T$). We will return to the bimodule
$\Ext^1_A(\tau T,T)$ in Part III.
	\smallskip
$\bullet$ {\bf Morita equivalence.} Tilting theory is a powerful generalization of Morita equivalence. 
This can  already be demonstrated very well by the reflection functors. 
When Gabriel showed that the representations of a Dynkin quiver correspond to the positive
roots and thus only in an inessential way on the given orientation, this was considered
as a big surprise. The BGP-reflection functors explain in
which way the representation theory of a quiver is independent of the orientation: 
one can change the orientation of all the arrows in a sink or a source, and use
reflection functors in order to obtain a bijection between the indecomposables.  
Already for the quivers of type $\Bbb A_n$ with $n \ge 3,$ we get interesting examples, relating
say a serial algebra (using one of the two orientations with just one sink and one source)
to a non-serial one. 

The reflection functors are still near to classical Morita theory, since no modules are really 
lost: here, we only deal with a kind of rearrangement of the categories in question. We
deal with split torsion pairs $(\Cal F, \Cal T)$ in $\mod A$ and
$(\Cal Y,\Cal X)$ in $\mod B$, with $\Cal F$ equivalent to $\Cal X$ and $\Cal T$ equivalent
to $\Cal Y.$ Let us call two hereditary artin algebras {\it similar} provided they
can be obtained from each other by a sequence of reflection functors. In case we consider
the path algebra of a quiver which is a tree, then any change of orientation leads to
a similar algebra. But already for the cycle with 4 vertices and 4 arrows, there
are two similarity classes, namely the quiver $\widetilde{\Bbb A}_{3,1}$ with a path of length 3, 
and the quivers of type $\widetilde{\Bbb A}_{2,2}$.

One property of the reflection functors should be mentioned (since it will be used in Part III). Assume that $i$
is a sink for $A$ (this means that the corresponding simple $A$-module $S(i)$ is projective).
Let $S'(i)$ be the corresponding simple $\sigma_iA$-module (it is injective).
If $M$ is any $A$-module, then $S(i)$ is not a composition factor of $M$ if and only if
$\Ext^1_{\sigma_iA}(S'(i),\sigma_iM) = 0.$ This is a situation, where the reflection 
functor yields a universal extension; for similar situations, let me refer to [R1].

The general tilting process is further away from classical Morita theory,
due to the fact that the torsion pair $(\Cal F,\Cal T)$ in $\mod A$ is no longer split. 
	\smallskip
$\bullet$ {\bf Duality theory.} Tilting theory is usually formulated as dealing with
equivalences of subcategories (for example, that $\Hom_A(T,-)\:\Cal T \to \Cal Y$
is an equivalence). However, one may also consider it as a duality theory,
by composing the equivalences obtained with the duality functor $D$, thus obtaining
a duality between subcategories of the category $\mod A$ and subcategories of the
category $\mod B^{\text{op}}$. The new formulations obtained in this way actually look more
symmetrical, thus may be preferable. Of course, as long as we deal with finitely generated
modules, there is no mathematical difference. This changes, as soon as one takes into consideration also modules which are not of finite length.

But the interpretation of tilting processes as dualities is always of interest, 
also when dealing with modules of finite length:
In [A], Auslander considers (for $R$ a left and
right noetherian ring) the class $\Cal W(R)$ of all left $R$-modules of the form $\Ext^1_R(N_R,R_R)$,
where $N_R$ is a finitely generated right $R$-module, and he asserts that
it would be of interest to know whether this class is always closed under submodules.
A first example of a ring $R$ with $\Cal W(R)$ not being closed under submodules
has been exhibited by Huang [Hg], namely the path algebra $R = kQ$ 
of the quiver $Q$ of type $\Bbb A_3$ with 2 sources.
Let us consider the general case when $R = A$ is a hereditary artin algebra. 
The canonical injective cogenerator $T = DA$ is a tilting module, thus
$\Ext^1_A(T,-)$ is a full and dense functor
from $\mod A$ onto $\Cal X(T)$. The composition of functors
$$
 \mod A^{\text{op}} @>D>> \mod A @>{\Ext^1(T,-)}>> \mod A 
$$
is the functor $\Ext^1_A(-,A_A)$, thus we see that $\Cal W(A) = \Cal X(T)$. 
On the other hand, $\Cal T(T)$ are the injective $A$-modules, they are
mapped under $\Hom_A(T,-)$ to the class $\Cal Y(T)$, 
and these are the projective $A$-modules. It follows
that $\Cal X(T)$ is the class of all $A$-modules without an indecomposable projective direct summand. As a consequence, {\it $\Cal W(A) = \Cal X(T)$ is closed
under submodules if and only if $A$ is a Nakayama algebra.} (It is an easy exercise to 
show that $\Cal X(T)$ is closed under submodules if and only if the injective envelope
of any simple projective module is projective, thus if and only if $A$ is a Nakayama
algebra). 

It should be stressed that Morita himself seemed to be more interested in dualities than
in equivalences. What is called Morita theory was popularized by P.M.Cohn and 
H. Bass, but apparently was considered by Morita as a minor addition to his duality theory.
When Gabriel heard about tilting theory, he immediately interpreted it as a non-commutative
analog of Roos duality.

The use of general tilting modules as a source for dualities has been shown to be
very fruitful in the representation theory of algebraic groups, of Lie algebras and of
quantum groups. This is explained in detail in Chapter 8 by Donkin. As a typical
special case one should have the classical {\bf Schur-Weyl duality} in mind, which relates
the representation theory of the general linear groups and that of the symmetric groups,
see Chapter 8 by Donkin, but also [KSX].

In the realm of commutative complete local noetherian rings, 
Auslander and Reiten [AR2] considered {\bf Cohen-Macaulay rings} with dualizing module $W$.
They showed that {\it $W$ is the only basic cotilting module.} On the basis of this 
result, they introduced the notion of a dualizing module for arbitrary artin algebras. 
	\medskip
$\bullet$ {\bf Schofield induction.} This is an inductive procedure for constructing
all exceptional modules starting with the simple ones, by forming exact sequences of the
following kind: Assume we deal with a hereditary $k$-algebra, where $k$ is algebraically
closed, and let $E_1, E_2$ be orthogonal exceptional modules with $\dim\Ext^1_A(E_1,E_2) =
t$ and $\Ext^1_A(E_2,E_1) = 0.$ Then, for every pair $(a_1,a_2)$ of positive natural numbers satisfying $a_1^2 + a_2^2 - ta_1a_2 = 1$, there exists (up to equivalence) a unique non-split 
exact sequence of the form
$$
 0 \to E_2^{a_2} \to E \to E_1^{a_1} \to 0
$$
(call it a {\it Schofield sequence}). Note that the middle term of such a Schofield sequence
is exceptional again, and  
it is an amazing fact that starting with the
simple $A$-modules without self-extension, all the exceptional $A$-modules are obtained in
this way. Even a stronger assertion is true: If $E$ is an exceptional module with 
support of cardinality $s$ (this means that $E$ has precisely $s$ different composition factors),
then there are precisely $s-1$ Schofield sequences with $E$ as middle term.
What is the relation to tilting theory? Starting with $E$ one obtains the Schofield sequences
by using the various indecomposable direct summands of its Bongartz complement: the
$s-1$ summands yield the $s-1$ sequences [R3]. 
	\smallskip
$\bullet$ {\bf Exceptional sequences, mutations.} Note that a tilted algebra is always
directed: the indecomposable summands of a tilting module $E_1,\dots,E_m$ can be ordered 
in such a way that $\Hom_A(E_i,E_j) = 0$ for $i > j.$ We may call such a sequence
$(E_1,\dots,E_m)$ a tilting sequence, and there is the following generalization which
is of interest in its own (and which was considered by the Rudakov school [Ru]):
Call $(E_1,\dots,E_m)$ an {\it exceptional sequence} provided all the modules $E_i$ are
exceptional $A$-modules and $\Hom_A(E_i,E_j) = 0$ and $\Ext^1_A(E_i,E_j) = 0$ for $i > j.$
There are many obvious examples of exceptional sequences which are not tilting
sequences, the most important one being sequences of simple modules in case the
$\Ext$-quiver of the simple modules is directed. Now one may be afraid that this generalization
could yield too many additional sequences, but this is not the case. In general most of the
exceptional sequences are tilting sequences! 
An exceptional sequence $(E_1,\dots,E_m)$ is said to be complete provided
$m = n(A)$ (the number of simple $A$-modules). There is a braid group action on the set
of complete exceptional sequences, and this action is transitive [C,R2]. 
This means that 
all the exceptional sequences can be obtained from each other by what one calls ``mutations''.
As a consequence, one obtains the following: If $(E_1,\dots,E_n)$ is a complete
exceptional sequence, then there is a permutation $\pi$ such that $\End_A(E_i) = 
\End_A(S_{\pi(i)})$, where $S_1,\dots, S_n$ are the simple $A$-modules. In particular, this means
that for any tilted algebra $B$, the radical factor algebras of $A$ and of $B$ are Morita equivalent.

An exceptional module $E$ defines also {\bf partial reflection functors} [R1] as follows:
consider the following full subcategories of $\mod A$. Let $\Cal M^E$ be given by all
modules $M$ with $\Ext_A^1(E,M) = 0$ such that no non-zero direct summand of $M$ 
is cogenerated by $E$; dually, let $\Cal M_E$ be given by all modules $M$ with $\Ext_A^1(M,E) = 0$ such that no non-zero direct summand of $M$ is generated by $E$; also, 
let $\Cal M^{-E}$ be given by all $M$ with $\Hom_A(M,E) = 0$
and $\Cal M_{-E}$ by all $M$ with $\Hom_A(E,M) = 0$. For any module $M$, let 
$\sigma^{-E}(M)$ be the intersection of the kernels of maps $M \to E$ and
$\sigma_{-E}(M) = M/t_EM$, where 
$t_EM$ is the sum of the images of maps $M \to E.$ In this way, we
obtain equivalences
$$
 \sigma^{-E} \: \Cal M^E/\langle E\rangle \longrightarrow \Cal M^{-E}, \T{and}
 \sigma_{-E} \: \Cal M_E/\langle E\rangle \longrightarrow \Cal M_{-E}.
$$
Here $\langle E\rangle$ is the ideal of all maps which factor through $\add E$.
The reverse functors $\sigma^E$ and $\sigma_E$ 
are provided by forming universal extensions by copies of $E$ (from
above or below, respectively). Note that on the level of dimension
vectors these partial functors $\sigma = \sigma^E, \sigma^{-E}, \sigma_E, \sigma_{-E}$
yield the usual reflection formula:
$$
 \bdim \sigma(M) = \bdim M - \frac{2\langle \bdim M, \bdim E\rangle}
 {\langle \bdim E, \bdim E\rangle}\bdim E.
$$
	\medskip
$\bullet$ {\bf Slices.} An artin algebra $B$ is a tilted algebra if and only if 
$\mod B$ has a slice. Thus the existence of slices characterizes
the tilted algebras. The necessity to explain the importance of slices has to be mentioned as 
a (further) impetus for the development of tilting theory. 
In my 1979 Ottawa lectures, I tried to describe several module
categories explicitly. At that time, the knitting of preprojective components
was one of the main tools, and I used slices in such components in order to guess 
what later
turned out to be tilting functors, namely functorial constructions using pushouts and
pullbacks. The obvious question about a possible theoretical foundation was raised by several
participants, but it could be answered only a year later at the Puebla conference.
Under minor restrictions (for example, the existence of a sincere indecomposable module)
preprojective components will contain slice modules and these are tilting modules with
a hereditary endomorphism ring! This concerns the concealed algebras to be mentioned
below, but also all the representation-directed algebras. Namely, 
using covering theory, the problem of describing the
structure of the indecomposable modules over a representation-finite algebra is reduced to the  
representation-directed algebras with a sincere indecomposable module, and such an algebra is
a tilted algebra, since it obviously has a slice module.

In dealing with an artin algebra of finite representation
type, and looking at its Auslander-Reiten quiver, one may ask for sectional subquivers
say of Euclidean types. Given such a subquiver $\Gamma$, applying several times $\tau$ or 
$\tau^{-1}$ (and obtaining in this way ``parallel'' subquivers), 
one has to reach a projective, or an injective vertex, respectively.
Actually, Bautista and Brenner have shown that the number of parallel subquivers is bounded,
the bound is called the {\bf replication number.} If one is interested in algebras with optimal replication numbers, one only has to look at representation-finite tilted algebras of Euclidean type. Note that given a hereditary algebra $A$ of Euclidean type and a tilting $A$-module $T$,
then $B = \End_A(T)$ is representation-finite if and only if $T$ has both preprojective
and preinjective indecomposable direct summands. 

It is natural to look inside preprojective and preinjective components for
slices. 
In 1979 one did not envision that there could exist even regular components with 
a slice module. But any connected wild hereditary algebra with at least
three simple modules has a regular tilting module $T$, and then the {\bf connecting
component} of $B = \End_A(T)$ is regular. One should be aware that the category $\mod B$ 
looks quite amazing: the connecting component (which is a regular component in this case) 
connects two wild subcategories, like a tunnel between two busy regions. Inside the tunnel,
there are well-defined paths for the traffic, and the traffic goes in just one direction.
	\medskip
Tilting modules can be used to study {\bf specific classes of artin algebras.}
Some examples have been mentioned already. We have noted that 
all the representation finite $k$-algebras, with $k$ algebraically
closed, can be described using tilted algebras (the condition on $k$ is needed
in order to be able to use covering theory). We obtain in this way very detailed
information on the structure of the indecomposables. One of the first uses of tilting
theory concerned the representation-finite tree algebras, see Chapter 3 by Br\"ustle. 
	\smallskip
Other examples:
	\smallskip
$\bullet$ {\bf The Concealed Algebras.} This concerns again algebras $B$ with a slice module.
By definition, $B$ is a {\it concealed} algebra, provided $B = \End_A(T),$ where $T$
is a preprojective $A$-module with $A$ hereditary. The tame concealed $k$-algebras $B$
where $k$ is algebraically closed, have been classified
by Happel and Vossieck, and Bongartz has shown in which way they can be used in order
to determine whether a $k$-algebra is representation-finite.
	\smallskip
$\bullet$ {\bf Representations of Posets.} The representation theory of posets always has
been considered as an important tool when studying questions in representation theory in
general: there are quite a lot of reduction techniques which lead to a vector space
with a bunch of subspaces, but the study of a vector space with a bunch of subspaces with some
inclusions prescribed, really concerns the representation theory of the
corresponding poset. On the other hand, the representation theory of finite posets is very similar to the representation theory of some quite well-behaved algebras, and the
relationship is often given by tilting modules. For example, when dealing with 
a disjoint union of chains, then we deal with the subspace representations of
a star quiver $Q$ (the quiver $Q$ is obtained from a finite set of linearly oriented
quivers of type $\Bbb A$, with all the sinks identified to one vertex, the center of the
star). If $c$ is the center of the star quiver $Q$, then the subspace representations
are the torsion-free modules of the (split) 
torsion pair $(\Cal Y,\Cal X)$, with $\Cal X$ being the
representations $V$ of $Q$ such that $V_c = 0.$ We also may consider the opposite quiver
$Q^{\text{op}}$ and the (again split) torsion pair $(\Cal F,\Cal T)$, where now $\Cal F$
are the representations $V$ of $Q^{\text{op}}$ with $V_c = 0.$ The two orientations
used here are obtained by a sequence of reflections, and the two split torsion pairs
$(\Cal F,\Cal G), (\Cal Y,\Cal X)$ are given by a tilting module which is a slice module:

$$
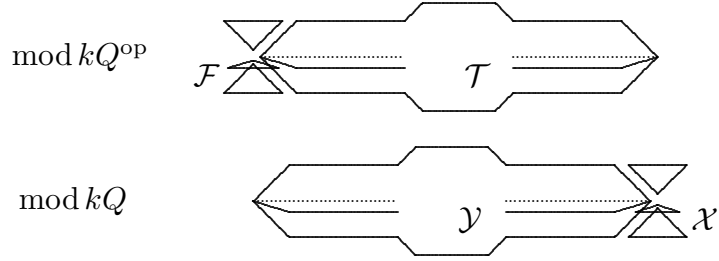
\hbox{\beginpicture
\setcoordinatesystem units <.4cm,.4cm>
\put{\beginpicture
\put{} at -6 -0.5
\put{} at 13 2.5
\plot 0 1  1 0  4 0  4.5 -0.5 6.5 -0.5 7 0  10 0  11 1  /
\plot 0 1  1 2  4 2  4.5  2.5 6.5  2.5 7 2  10 2  11 1  /
\setdots <.5mm>
\plot 0 1  4 1 /
\plot 7 1  11 1 /
\setsolid
\plot 0 1  1 0.7       4 0.7 /
\plot 7 0.7  10 0.7  11 1 /

\plot -0.2 1.2  0.6 2  -1 2  -0.2 1.2 /
\plot -0.2 0.8  0.6 0  -1 0  -0.2 0.8 /
\plot -0.2 0.9  0.5 0.7  -.9 0.7  -0.2 0.9 /

\put{$\mod kQ^{\text{op}}$} at -5 1
\put{$\Cal F$} at -1.5 0.5
\put{$\Cal T$} at  6 0.5
\endpicture} at 0 0 
\put{\beginpicture
\put{} at -6 -0.5
\put{} at 13 2.5
\plot 0 1  1 0  4 0  4.5 -0.5 6.5 -0.5 7 0  10 0  11 1  /
\plot 0 1  1 2  4 2  4.5  2.5 6.5  2.5 7 2  10 2  11 1  /
\setdots <.5mm>
\plot 0 1  4 1 /
\plot 7 1  11 1 /
\setsolid
\plot 0 1  1 0.7       4 0.7 /
\plot 7 0.7  10 0.7  11 1 /

\plot 11.2 1.2  10.4 2  12 2  11.2 1.2 /
\plot 11.2 0.8  10.4 0  12 0  11.2 0.8 /
\plot 11.2 0.9  10.6 0.7  11.8 0.7  11.2 0.9 /

\put{$\mod kQ$} at -5 1
\put{$\Cal X$} at 12.5 0.5
\put{$\Cal Y$} at  6 0.5
\endpicture} at 0 -4
\endpicture}
$$
	\smallskip
$\bullet$ {\bf The Crawley-Boevey-Kerner Functors.} If $R$ is an artin algebra and $W$
an $R$-module, let us write $\langle \tau^\bullet W\rangle $ for the ideal of
$\mod R$ of all maps which factor through a direct sum of modules of the form $\tau^z W$
with $z\in \Bbb Z.$ We say that the module categories $\mod R$ and $\mod R'$ are
{\it almost equivalent} provided there is an $R$-module $W$ and an $R'$-module $W'$
such that the categories $\mod R/\langle \tau^\bullet W\rangle$
and $\mod R'/\langle \tau^\bullet W'\rangle$ are equivalent. 
The Crawley-Boevey-Kerner functors were introduced in order to show the following:
{\it If $k$ is a field and $Q$ and $Q'$ are connected wild quivers, then the categories 
$\mod kQ$ and $\mod kQ'$ are almost equivalent.} The proof uses tilting modules, and 
the result may be rated as one of the most spectacular applications of tilting theory.
Thus it is worthwhile to outline the essential ingredients. This will be done below.

Here are some remarks concerning almost equivalent categories. 
It is trivial that the module categories of all representation-finite artin algebras are almost equivalent. 
If $k$ is a field, and $Q, Q'$ are tame connected quivers, then $\mod kQ$ and $\mod kQ'$ are 
almost equivalent only if $Q$ and $Q'$ have the same type ($\widetilde {\Bbb A}_{pq},
\widetilde {\Bbb D}_{n}, \widetilde {\Bbb E}_6,\widetilde {\Bbb E}_7,\widetilde {\Bbb E}_8$).
Let us return to wild quivers $Q, Q'$ and a Crawley-Boevey-Kerner  equivalence 
$$
 \eta\: \mod kQ/\langle\tau^\bullet W\rangle \longrightarrow 
 \mod kQ'/\langle\tau^\bullet W'\rangle,
$$
with finite length modules $W, W'$. 
Consider the case of an uncountable base field $k$,
so that there are uncountably many isomorphism classes of indecomposable modules for
$R = kQ$ as well as for $R' = kQ'$.
The ideals $\langle \tau^\bullet W\rangle$
and $\langle \tau^\bullet W'\rangle$ are given by the maps which factor through a countable set
of objects, thus nearly all the indecomposable modules remain indecomposable
in $\mod kQ/\langle \tau^\bullet W\rangle$
and $\mod kQ'/\langle \tau^\bullet W'\rangle$, and non-isomorphic ones (which are not sent to zero)
remain non-isomorphic.
In addition, one should note that the 
equivalence $\eta$ is really constructive 
(not set-theoretical rubbish), with no unfair choices whatsoever. 
This will be clear from the further discussion.

Nearly all quivers are wild. For example, if we consider the 
$m$-subspace quivers $Q(m)$, then one knows that $Q(m)$ is  wild  provided $m \ge 5$. 
Let us concentrate on a comparison of the wild quivers $Q(6)$ and $Q(5)$. 
To assert that $Q, Q'$ are wild quivers means
that there are full embeddings $\mod kQ \to \mod kQ'$ and $\mod kQ' \to \mod kQ$. But the
Crawley-Boevey-Kerner theorem  provides a completely
new interpretation of what ``wildness'' is about.
The definition of ``wildness'' itself is considered as quite odd, since it means
in particular that there is a full embedding of $\mod kQ(6)$ into $\mod kQ(5)$.
One may reformulate the wildness assertion as follows: 
any complication which occurs for $6$ subspaces can
be achieved (in some sense) already for $5$ subspaces. But similar results are known in  mathematics, since one is aware of other categories 
which allow to realize all kinds of categories as a subcategory. Also, ``wildness'' may be interpreted as
a kind of fractal behaviour: inside the category $\mod kQ(5)$ we find proper full subcategories
which are equivalent to $\mod kQ(5),$ again a quite frequent behaviour. 
These realization results are concerned with
small parts of say $\mod kQ(5)$; one looks at full subcategories of the category $\mod kQ$ 
which have desired properties, but one does not try to control the complement. This is in sharp
contrast to the Crawley-Boevey-Kerner property which 
provides a {\bf global} relation between $\mod kQ(5)$ and $\mod kQ(6)$, actually, between 
the module categories of any two wild connected quivers. 
In this way we see that there is a kind of homogeneity property of wild 
abelian length categories which had not been anticipated before.  

The Crawley-Boevey-Kerner result may be considered as a sort
of Schr{\"o}\-der-Bern\-stein property for abelian length categories. Recall that
the Schr{\"o}\-der-Bern\-stein theorem asserts that if two sets $S, S'$ can be embedded
into each other, then there is a bijection $S \to S'$. For any kind of mathematical
structure with a notion of embedding, one may ask whether two objects are isomorphic
in case they can be embedded into each other. Such a property is very rare, even if we 
replace the isomorphism requirement by some weaker requirement. But this is what
is asserted by the Crawley-Boevey-Kerner property.

Let us outline the construction of $\eta$. We start with a connected 
wild hereditary artin algebra $A$, and a regular exceptional module $E$ which is 
quasi-simple
(this means that the Auslander-Reiten sequence ending in $E$ has indecomposable middle 
term, call it $\mu(E)$), such a module exists provided $n(A) \ge 3.$
Denote by $E^\perp$ the category of all $A$-modules $M$ such that 
$\Hom_A(E,M) = 0 = \Ext^1_A(E,M).$ One knows (Geigle-Lenzing, Strau{\ss}) that
$E^\perp$ is equivalent to the category $\mod C$, where $C$ is a connected wild
hereditary algebra $C$ and $n(C) = n(A)-1$. The aim is to compare the
categories $\mod C$ and $\mod A$, they are shown to be almost equivalent. 

It is easy to see that the 
module $\mu(E)$ belongs to $E^\perp$, thus it can be regarded as a $C$-module.
Since $E^\perp = \mod C$, there is a projective generator $T'$ in $E^\perp$ with
$\End_A(T') = C.$ Claim: {\it $T'\oplus E$ is a tilting module.} 
For the proof we only have to check that
$\Ext_A^1(T',E) = 0.$ 
Since $T'$ is projective in $E^\perp,$ it follows that $\Ext^1(T',\mu(E)) = 0.$
However, there is a surjective map $\mu(E) \to E$ and this induces a surjective
map $\Ext^1(T',\mu(E)) \to \Ext^1(T',E).$ 

As we know, 
the tilting module $T = \overline T \oplus E$ defines a torsion pair $(\Cal F,\Cal T)$,
with $\Cal T$ the $A$-modules generated by $T$. Let us denote by $\tau_{\Cal T}M = t\tau_A M$
the torsion 
submodule\note{The notation shall indicate that this functor $\tau_{\Cal T}$ has to be considered
as an Auslander-Reiten translation: it is the relative Auslander-Reiten translation in the
subcategory $\Cal T$. And there is the equivalence $\Cal T \simeq \Cal Y$, where $\Cal Y$
is a full subcategory of $\mod B,$ with $B = \End_A(T)$. Since $\Cal Y$ is closed under $\tau$
in $\mod B$, the functor $\tau_{\Cal T}$ corresponds to the Auslander-Reiten translation 
$\tau_B$ in
$\mod B$.} of $\tau_A M$.
The functor $\eta$ is now defined as follows: 
$$
 \eta(M) = \lim_{t\to \infty} \tau_A^{-t}\tau_{\Cal T}^{2t}\tau_C^t(M).
$$
One has to observe that the limit actually stabilizes: 
for large $t$, there is no difference whether we
consider $t$ or $t+1$. The functor $\eta$ is full, the image is just the full subcategory 
of all regular $A$-modules. There is a non-trivial kernel: a map is send to zero if and only if
it belongs to $\langle \tau^\bullet W\rangle$, where 
$W = C \oplus \mu(E) \oplus DC$. Also, let $W' = A\oplus DA$. 
Then $\eta$ is an equivalence
$$
 \eta\: \mod C/\langle\tau^\bullet W\rangle \longrightarrow \mod A/\langle\tau^\bullet W'\rangle.
$$

One may wonder how special the assumptions on $A$ and $C$ are. 
Let us say that $A$ {\it dominates} $C$ provided there exists a 
regular exceptional module $E$ which is quasi-simple with $\mod C$ equivalent to $E^\perp$.
Given any two wild connected quivers $Q, Q'$, there is a sequence of wild connected
quivers $Q = Q_0,\dots,Q_t = Q'$ such that $kQ_i$ either dominates  or is dominated by 
$kQ_{i-1}$, for $1 \le i \le t.$ 
This implies that the module categories of all wild path algebras are almost
equivalent.
	\medskip
The equivalence $\eta$ can be constructed also in a different way [KT], 
using partial reflection functors. Let $E(i) = \tau^i E,$ for all $i\in \Bbb Z$.
Note that for any regular $A$-module $M$, one knows that
$$
 \Hom_A(M,E(-t)) = 0 = \Hom_A(E(t),M) \T{for} t \gg 0,
$$
according to Baer and Kerner. 
Thus, if we choose $t$ sufficiently large, we can apply the partial reflection functors
$\sigma^{E(-t)}$ and $\sigma_{E(t)}$ to $M$.
The module obtained from $M$ has the form
$$
\hbox{\beginpicture
\setcoordinatesystem units <1cm,.6cm>
\multiput{$E(-t)$} at 0 2  1 2  2.5 2 /
\put{$M$} at 1.25 1
\multiput{$E(t)$} at 0 0  1 0  2.5 0 /
\plot -0.5 1.5  3 1.5 /
\plot -0.5 0.5  3 0.5 /
\setdots <1mm>
\plot 1.6 2  2 2 /
\plot 1.6 0  2 0 /
\endpicture}
$$ 
and belongs to 
$$
 \Cal M^{E(-t)} \cap  \Cal M_{E(t)} \quad \subseteq \quad
   \Cal M^{-E(-t+1)} \cap \Cal M_{-E(t-1)}.
$$
Thus we can proceed, applying now $\sigma^{E(-t+1)}$ and 
$\sigma_{E(t-1)}.$  We use induction, the last partial
reflection functors to be applied are the functors $\sigma^{E(-1)}$, $\sigma_{E(1)}$, 
and then finally $\sigma^{E}$. In this way we obtain a module in
$$
 \Cal M_{E(1)} \cap \Cal M^{E} = E^\perp 
$$
as required. It has the following structure:
$$
\hbox{\beginpicture
\setcoordinatesystem units <1cm,.6cm>
\multiput{$E(-t)$} at 1 2  2.5 2 /
\put{$M$} at 1.75 1
\multiput{$E(t)$} at  1 0  2.5 0 /

\multiput{$E(-1)$} at 0 4  3.5 4 /
\multiput{$E(1)$} at  -0.2 -2  3.5 -2 /
\multiput{$E$} at  -0.5 5  4 5 /
\plot 0.5 1.5  3 1.5 /
\plot 0.5 0.5  3 0.5 /

\plot -0 2.5  3.5 2.5 /
\plot -0 -0.5 3.5 -0.5 /

\plot -.5 3.5  4 3.5 /
\plot -.4 -1.5 3.7 -1.5 /

\plot -.7 4.5 4.2 4.5 /

\setdots <1mm>
\plot 1.6 2  2 2 /
\plot 1.6 0  2 0 /

\plot 0.8 4  2.8 4 /
\plot 0.8 -2  2.8 -2 /

\plot 0.3 5  3.3 5 /

\plot 1.75 2.7  1.75 3.3 /
\plot 1.75 -.7  1.75 -1.3 /

\endpicture}
$$
	\smallskip
$\bullet$ {\bf The shrinking functors for the tubular algebras.} Again these are
tilting functors (here, $A$ no longer is a hereditary artin algebra, but say a canonical
algebra - we are still in the realm of the ``T'' displayed in Part I, 
now even in its center), 
and such functors belong to the origin of the development. If one
looks at the Brenner-Butler tilting paper, the main
examples considered there were of this kind. So one of the first applications of
tilting theory was to show the similarity of the module categories of various tubular algebras.
And this is also the setting which later helped to describe in detail the module category of a
tubular algebra: one uses the shrinking functors in order to construct all the
regular tubular families, as soon as one is known to exist.
	\smallskip
$\bullet$ {\bf Self-Injective Algebras.} Up to coverings and (in characteristic 2) deformations, the trivial
extensions of the tilted algebras of Dynkin type (those related to the left
arm of the ``T'' displayed in Part I) yield all the representation-finite self-injective algebras
(recall that the trivial extension of an algebra $R$ is the semi-direct product $R$\semi$DR$ 
of $R$
with the dual module $DR$). 
In private conversation, such a result was conjectured by Tachikawa already in 1978, and it was
the main force for the investigations of him and Wakamatsu, which he presented at
the Ottawa conference in 1979. There he also dealt with the trivial extension of 
a tilted algebra of Euclidean type (the module category has two tubular families).
This motivated Hughes-Waschb\"usch to introduce the concept of a repetitive algebra.
But it is also part of one of the typical quarrels between Z\"urich and the rest of the
world: with Gabriel hiding the Hughes-Waschb\"usch manuscript from Bretscher-L\"aser-Riedtmann 
(asking a secretary to seal the envelope with the manuscript 
and to open it only several months later...), so that they could proceed ``independently''.

The representation theory of artin algebras came into limelight when Dynkin diagrams
popped up for representation-finite algebras.
And this occurred twice, first for hereditary artin algebras in the work of Gabriel
(as the $\Ext$-quiver), 
but then also for self-injective algebras in the work of Riedtmann (as the tree class of
the stable Auslander-Reiten quiver). The link between these two classes of rings is
furnished by tilted algebras and their trivial extensions. As far as I know, it is
Tachikawa who deserves the credit for this important insight.

The reference to trivial extensions of tilted algebras actually closes the circle
of our considerations, due to another famous theorem of Happel. We have started with the 
fact that tilting functors provide derived equivalences. Thus the derived category
of a tilted algebra can be identified with the derived category of a hereditary artin
algebra. For all artin algebras $R$ of finite global dimension (in particular our
algebras $A$ and $B$), there is a an equivalence between $D^b(\mod R)$ and the stable module
category of the repetitive algebra $\widehat R$. But $\widehat R$ is just a $\Bbb Z$-covering
of the trivial extension of $R$.
	\medskip
{\bf Artin Algebras with Gorenstein Dimension at most 1.} We have mentioned that the two 
classes of algebras: the
selfinjective ones and the hereditary ones, look very different, but nevertheless they 
have some common behaviour.  
Auslander and Reiten [AR] have singled out an important
property which they share, they are Gorenstein algebras of Gorenstein dimension at most 1. 
An artin algebra $A$ is called  
{\it Gorenstein}\note{This definition is one of the many possibilities to
generalize the notion of a commutative Gorenstein ring to a non-commutative setting.
Note that a commutative artin algebra $R$ is a Gorenstein ring if and only if $R$
is selfinjective. Of course, a commutative connected artin algebra $R$ is a local ring,
and a local ring has a non-zero module of finite projective dimension only in case
$R$ is selfinjective.}
provided ${}_ADA$ has finite projective dimension and ${}_AA$ 
has finite injective dimension.
For Gorenstein algebras, $\projdim {}_ADA = \injdim {}_AA$ according to Happel, 
and this number is called the {\it Gorenstein dimension} of $A$. 
It is not known whether the finiteness of the projective dimension of $_ADA$ 
implies the finiteness of the injective dimension of ${}_AA$. It is conjectured
that this is the case: this is the Gorenstein symmetry conjecture, and this
conjecture is equivalent to the conjecture that the small finitistic dimension
of $A$ is finite [AR]. 
The artin algebras of Gorenstein dimension 0 
are the selfinjective algebras. An artin algebra has Gorenstein dimension at most 1 if and
only if $DA$ is a tilting module (of projective dimension at most 1). 

If $A$ is a Gorenstein algebra of Gorenstein dimension at most 1, then there is a strict
separation of the indecomposable modules: an $A$-module $M$ of finite projective dimension
or finite injective dimension satisfies both $\projdim M \le 1$, $\injdim M \le 1.$
(The proof is easy: Assume $\projdim  M \le m$, thus the $m$-th syzygy module 
  $\Omega_m(M)$ is projective.  Now for any short exact sequence
 $0 \to X \to Y \to Z$, it is clear that $\injdim X \le 1, \injdim Y \le 1$ imply
 $\injdim Z \le 1.$ One applies this inductively to the exact sequences
 $\Omega_i(M) \to P_i \to \Omega_{i-1}(M)$, where $P_i$ is projective, starting with 
 $i = m$ and ending with $i = 0.$ This shows that $\injdim(M) \le 1.$ The dual
 argument shows that a module of finite injective dimension has projective dimension at most 1.)
 As a consequence, if $A$ is not hereditary, then the global dimension of $A$ is infinite.
 Also, if $P$ is an indecomposable projective $A$-module, then
 either its radical is projective
 or else the top of $P$ is a simple module which has infinite projective and infinite
 injective dimension. 

Until very recently, the interest in artin algebras of Gorenstein dimension 
at most 1 has been
quite moderate, the main reason being a lack of tempting examples: of algebras which are
neither selfinjective nor hereditary. But now there is a wealth of such examples, as we will
see in Part III. 

	\medskip
We hope that we have convinced the reader that the use of tilting modules and 
tilted algebras lies at the heart of
nearly all the major developments in the representation theory of artin algebras in the last
25 years. In this report we usually restrict to tilting modules in the narrow sense (as
being finite length modules of projective dimension at most 1). In fact, most of the
topics mentioned are related to tilting $A$-modules $T$, where $A$ is a hereditary artin
algebra (so that there is no need to stress the condition $\projdim T \le 1).$ 
However, the following two sections will widen the viewpoint, taking into account 
also various generalizations. 
	\medskip
$\bullet$ {\bf Representations of semisimple complex Lie algebras and algebraic groups.}
The highest weight categories which arise in the representation theory of
semisimple complex Lie algebras and algebraic groups can be analyzed very well
using quasi-hereditary artin algebras as introduced by Cline-Parshal-Scott. One of the
main features of such a quasi-hereditary artin algebra is its characteristic module,
this is a tilting module (of finite projective dimension). Actually, the experts use
a different convention, calling its indecomposable direct summands ``tilting modules'',
see Chapter 8 by Donkin. If $T$ is the characteristic module, then $\add T$ consists
of the $A$-modules which have both a standard filtration and a costandard filtration,
and it leads to a duality theory which seems to be of great interest.
	\medskip
$\bullet$ {\bf The Homological Conjectures.} The homological conjectures are 
one of the central themes of module theory, so clearly they deserve special interest.
They go back to mathematicians like Nakayama, Eilenberg, Auslander, Bass, but
also Rosenberg, Zelinsky, Buchsbaum and Nunke should be mentioned,  and were formulated
between 1940 and 1960. Unfortunately, there are no written accounts about the origin,
but we may refer to surveys by Happel, Smal{\o} and Zimmermann-Huisgen. 
The modern development in representation theory of artin 
algebras was directed towards a solution of the Brauer-Thrall conjectures, and there was
for a long time a reluctance to work on the homological conjectures. The investigations
concerning the various representation types have produced a lot of information on special
classes of algebras, but for these algebras the homological conjectures are usually
true for trivial reasons. As Happel has pointed out, the lack of knowledge of non-trivial
examples may very well mean that  counter-examples could exist.
Here is a short discussion of this topic, in as far as modules without
self-extensions are concerned.

Let me start with the Nakayama conjecture which according to B.~M\"uller can be phrased as
follows: If $R$ is an artin algebra and
$M$ is a generator and cogenerator for $\mod R$ with $\Ext^i_R(M,M) = 0$ for all $i \ge 1$,
then $M$ has to be projective. Auslander and Reiten [AR1] proposed in 1975 that the same
conclusion should hold even if $M$ is not necessarily a cogenerator (this is called
the ``generalized Nakayama conjecture''). This incorporates a conjecture due to Tachikawa (1973):
If $R$ is self-injective and $M$ is an $R$-module with $\Ext^i_R(M,M) = 0$ for all $i > 0$, then
$M$ is projective. The relationship of the generalized Nakayama conjecture with tilting 
theory was noted in the Auslander-Reiten paper on
contravariantly finite subcategories [AR2]. 
Then there is the conjecture
on the finiteness of the number of complements of an almost complete partial tilting
module, due to Happel and Unger. And there
is a conjecture made by Beligiannis and Reiten [BR], 
called the Wakamatsu tilting conjecture (because it deals with Wakamatsu
tilting modules, see Chapter 7 by Reiten): 
If $T$ is a Wakamatsu tilting module of finite projective
dimension, then $T$ is a tilting module. 
The Wakamatsu tilting conjecture implies the generalized Nakayama conjecture
(apparently, this was first observed by Buan) and also the Gorenstein symmetry 
conjecture, see [BR]. 
In a joint paper, Mantese and Reiten [MR] showed that it is implied by the
finitistic dimension conjecture, and that it implies the conjecture on a finite number of complements, which according to Buan and Solberg is known to imply the 
generalized Nakayama conjecture. 
There is also the equivalence of the generalized Nakayama conjecture with
projective almost complete partial tilting modules having only a finite number of
complements (Happel-Unger, Buan-Solberg, both papers are in the Geiranger proceedings). 
Coelho, Happel and Unger proved that the finitistic dimension conjecture
implies the conjecture on a finite number of complements.

Further relationship of tilting theory with the finitistic dimension conjectures is discussed
in detail in Chapter 11 by Trlifaj and in Chapter 12 by Solberg. But also other results
presented in the Handbook have to be seen in this light. We know from 
Auslander and Reiten, that the finitistic dimension of an artin algebra $R$ is
finite, in case the subcategory  of all modules of finite projective dimension is
contravariantly finite in $\mod R$. This has been the motivation to look at the
latter condition carefully (see for example Chapter 9 by Unger). 
	\medskip
With respect to applications outside of ring and module theory, many more topics could
be mentioned. We have tried to stay on a basic
level, whereas there are a lot of mathematical objects which are derived from 
representation theoretical data and this leads to a fruitful interplay 
(dealing with questions on quantum groups, with the shellability of simplicial complexes,
or with continued fraction expansions of real numbers): There are many 
unexpected connections to analysis, to number theory, to combinatorics --- and again, it is usually the tilting theory which plays an important role.
	\bigskip
\centerline{\bf III.}
	\medskip
Let me repeat: at the time the Handbook was conceived, 
there was a common feeling that the tilted algebras (as the core of tilting theory)
were understood well and that this part of the theory had reached a sort of final shape.
But in the meantime this has turned out to be wrong: the tilted algebras have to be
seen as factor algebras of the so called cluster tilted algebras, and it may very well be, that 
in future the cluster tilted algebras and the cluster categories 
will topple the tilted algebras.
The impetus for introducing and studying cluster tilted algebras came from outside, in a
completely 
unexpected way. We will mention below some of the main steps of this development. 
But first let me jump directly to the relevant construction.
	\medskip
{\bf The cluster tilted algebras.}
We return to the basic setting, the hereditary artin algebra $A$, the tilting $A$-module
$T$ and its endomorphism ring $B$. Consider the semi-direct ring extension 
$$
 \widetilde B = B\,\semi\Ext^2_B(DB,B).
$$ 
This is called the {\it cluster tilted algebra} corresponding to $B$.
Since this is the relevant definition, let me say a little more about this 
construction\note
{One may wonder what properties the semi-direct product $R \semi \Ext^2_R(DR,R)$ for any
artin algebra $R$ has in general (at least in case $R$ has global dimension at most 2);
it seems that this question has not yet been studied.
}: 
$\widetilde B$ has $B$ as a subring, and there
is an ideal $J$ of $\widetilde B$ with $J^2 = 0$, 
such that $\widetilde B = B \oplus J$ as additive groups
and $J$ is as a $B$-$B$-bimodule isomorphic to $\Ext^2_B(DB,B)$; in order to construct
$\widetilde B$ one may take $B \oplus \Ext^2_B(DB,B)$, with 
componentwise addition, and one uses $(b,x)(b',x') = (bb',bx'+xb'),$ 
for $b,b'\in B$ and $x, x' \in \Ext^2_B(DB,B)$ as multiplication.

We consider again the example of $B$ given by a square with two zero relations. 
Here $\Ext^2_B(DB,B)$ is 8-dimensional
and $\widetilde B$ is a 16-dimensional algebra:
$$
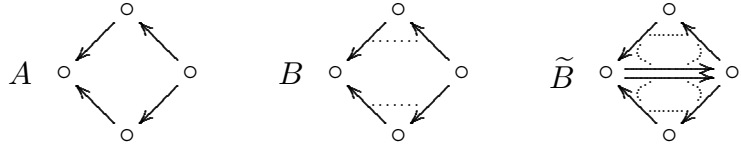
\hbox{\beginpicture
\setcoordinatesystem units <1cm,1cm>
\put{\beginpicture
\setcoordinatesystem units <.7cm,.7cm>
\multiput{$\circ$} at 0 1  1 0  1 2  2 1 /
\arr{0.8 1.8}{0.2 1.2}
\arr{0.8 0.2}{0.2 0.8}
\arr{1.8 1.2}{1.2 1.8}
\arr{1.8 0.8}{1.2 0.2}
\put{$A$} at -0.7 1
\endpicture} at -3 0
 
\put{\beginpicture
\setcoordinatesystem units <.7cm,.7cm>
\multiput{$\circ$} at 0 1  1 0  1 2  2 1 /
\arr{0.8 1.8}{0.2 1.2}
\arr{0.8 0.2}{0.2 0.8}
\arr{1.8 1.2}{1.2 1.8}
\arr{1.8 0.8}{1.2 0.2}
\setdots <.8mm>
\plot 0.6 1.5  1.4 1.5 /
\plot 0.6 0.5  1.4 0.5 /
\put{$B$} at -0.7 1
\endpicture} at 0 0 

\put{\beginpicture
\setcoordinatesystem units <.7cm,.7cm>
\multiput{$\circ$} at 0 1  1 0  1 2  2 1 /
\arr{0.8 1.8}{0.2 1.2}
\arr{0.8 0.2}{0.2 0.8}
\arr{1.8 1.2}{1.2 1.8}
\arr{1.8 0.8}{1.2 0.2}

\arr{0.3 0.93}{1.7 0.93}
\arr{0.3 1.07}{1.7 1.07}
\setdots <.5mm>
\plot 0.7 1.6  1.25 1.6 /
\plot 0.7 0.4  1.25 0.4 /
\setquadratic
\plot 0.7 1.15  0.5 1.3  0.6 1.5 /
\plot 1.3 1.15  1.5 1.3  1.4 1.5 /
\plot 0.7 0.85  0.5 0.7  0.6 0.5 /
\plot 1.3 0.85  1.5 0.7  1.4 0.5 /
\put{$\widetilde B$} at -0.7 1
\endpicture} at 3 0 

\endpicture}
$$
Non-isomorphic tilted algebras $B$ may yield isomorphic cluster tilted algebras $\widetilde B.$
Here are all the tilted algebras which lead to the cluster tilted algebra just considered:
$$
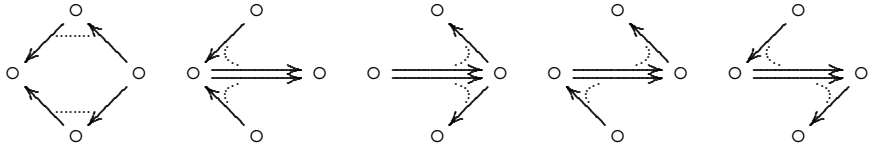
\hbox{\beginpicture
\setcoordinatesystem units <1cm,1cm>
\put{\beginpicture
\setcoordinatesystem units <.7cm,.7cm>
\multiput{$\circ$} at 0 1  1 0  1 2  2 1 /
\arr{0.8 1.8}{0.2 1.2}
\arr{0.8 0.2}{0.2 0.8}
\arr{1.8 1.2}{1.2 1.8}
\arr{1.8 0.8}{1.2 0.2}

\setdots <.5mm>
\plot 0.7 1.6  1.25 1.6 /
\plot 0.7 0.4  1.25 0.4 /
\endpicture} at 0 0 

\put{\beginpicture
\setcoordinatesystem units <.7cm,.7cm>
\multiput{$\circ$} at 0 1  1 0  1 2  2 1 /
\arr{0.8 1.8}{0.2 1.2}
\arr{0.8 0.2}{0.2 0.8}

\arr{0.3 0.93}{1.7 0.93}
\arr{0.3 1.07}{1.7 1.07}
\setdots <.5mm>
\setquadratic
\plot 0.7 1.15  0.5 1.3  0.6 1.5 /
\plot 0.7 0.85  0.5 0.7  0.6 0.5 /
\endpicture} at 2 0 

\put{\beginpicture
\setcoordinatesystem units <.7cm,.7cm>
\multiput{$\circ$} at 0 1  1 0  1 2  2 1 /
\arr{1.8 1.2}{1.2 1.8}
\arr{1.8 0.8}{1.2 0.2}

\arr{0.3 0.93}{1.7 0.93}
\arr{0.3 1.07}{1.7 1.07}
\setdots <.5mm>
\setquadratic
\plot 1.3 1.15  1.5 1.3  1.4 1.5 /
\plot 1.3 0.85  1.5 0.7  1.4 0.5 /
\endpicture} at 4 0 

\put{\beginpicture
\setcoordinatesystem units <.7cm,.7cm>
\multiput{$\circ$} at 0 1  1 0  1 2  2 1 /

\arr{0.8 0.2}{0.2 0.8}
\arr{1.8 1.2}{1.2 1.8}

\arr{0.3 0.93}{1.7 0.93}
\arr{0.3 1.07}{1.7 1.07}
\setdots <.5mm>
\setquadratic

\plot 1.3 1.15  1.5 1.3  1.4 1.5 /
\plot 0.7 0.85  0.5 0.7  0.6 0.5 /

\endpicture} at 6 0 

\put{\beginpicture
\setcoordinatesystem units <.7cm,.7cm>
\multiput{$\circ$} at 0 1  1 0  1 2  2 1 /
\arr{0.8 1.8}{0.2 1.2}
\arr{1.8 0.8}{1.2 0.2}

\arr{0.3 0.93}{1.7 0.93}
\arr{0.3 1.07}{1.7 1.07}
\setdots <.5mm>
\setquadratic
\plot 0.7 1.15  0.5 1.3  0.6 1.5 /
\plot 1.3 0.85  1.5 0.7  1.4 0.5 /
\endpicture} at 8 0 

\endpicture}
$$
	\medskip
It is quite easy to write down the {\bf quiver of a cluster-tilted algebra.} Here, we
assume that we deal
with $k$-algebras, where $k$ is an algebraically closed field. 
We get the quiver of $\widetilde B$ from the quiver with relations of $B$ by just replacing 
the dotted arrows\note{Actually, the usual convention for indicating relations
is to draw dotted lines, 
not dotted arrows. However, these dotted lines are to be seen
as being directed, since the corresponding relations are linear combinations of
paths with fixed starting point and fixed end point.}
by solid arrows in opposite direction [ABS]. 
The reason is the following: Let us denote by 
$N = \rad B$ the radical of $B$. Then $N\oplus J$
{\it is the radical of} 
$\widetilde B = B\, \semi\, J$, and $N^2 \oplus (NJ+JN)$ {\it is equal
to the square of the radical of $\widetilde B$.} This shows that the additional arrows for
$\widetilde B$ correspond to $J/(NJ+JN)$. Note that $J/(NJ+JN)$ is the top
of the $B$-$B$-bimodule $J$. Now the top of the bimodule $\Ext^2_B(DB,B)$ is
$\Ext^2_B(\soc {}_BDB,\top {}_BB)$, since $B$ has global dimension at most 2. It is well-known 
that $\Ext^2_B(\soc {}_BDB, \top {}_BB)$ describes the relations of the algebra $B$, and we
see in this way that relations for $B$ correspond to the additional arrows for $\widetilde B$. Since $\rad \widetilde B = \rad B \oplus J$ and $J$ is an
ideal of $\widetilde B$ with $J^2 = 0,$ 
we also see: {\it If $(\rad B)^t = 0$, then $(\rad \widetilde B)^{2t}
= 0$.} The quiver of any tilted algebra is directed, thus $(\rad B)^{n(B)} = 0,$
therefore $(\rad \widetilde B)^{2n(B)} = 0.$
	
The recipe for obtaining the quiver of $\widetilde B$ shows that there are always
oriented cyclic paths (unless $B$ is hereditary). However, such a path is always
of length at least 3. Namely, since the quiver of $B$ is 
directed, it follows that no relation of $B$ is a loop, thus
the quiver of $\widetilde B$ cannot have a loop [B-T]. Also, Happel ([H2], 
Lemma IV.1.11) 
has shown that for simple $B$-modules $S, S'$ with $\Ext^1_B(S,S') \neq 0$ one has
$\Ext^2_B(S,S') = 0.$ This means that 
the quiver of $\widetilde B$ cannot have a pair of arrows in opposite direction [BMR2].

It should be of interest whether knowledge about the quiver with relations of a cluster tilted algebra $\widetilde B$ can provide new insight into the structure of the tilted algebras themselves.
There is a lot of ongoing research on cluster tilted algebras, 
let us single out just one result. Assume that we deal with
$k$-algebras, where $k$ is algebraically closed. Then: 
{\it Any cluster tilted $k$-algebra 
of finite representation type is uniquely determined by its quiver} [BMR3].
This means: in the case of finite representation type, the quiver determines the relations! 
What happens in general is still under investigation.
	\bigskip
If $A$ is a hereditary artin algebra and $T$ a tilting $A$-module with endomorphism ring
$B$, we have introduced the corresponding cluster tilted algebras as the algebra
$\widetilde B = B \semi J$, with $B$-$B$-bimodule $J = \Ext^2_B(DB,B)$. The 
original definition of $\widetilde B$ by Buan, Marsh and Reiten [BMR1] used another
description of $J$, namely $J = \Ext^1_A(T,\tau^{-1}T)$, and 
it was observed by Assem, Br\"ustle and Schiffler [ABS] that the bimodules 
$\Ext^1_A(T,\tau^{-1}T)$ and $\Ext^2_B(DB,B)$ are 
isomorphic\note{In addition, we should remark that $\Ext_A^1(T,\tau^{-1}T)$ can be 
identified with $\Ext_A^1(\tau T,T)$ (as $B$-$B$-bimodules). The reason is the fact that
the functor $\tau^{-1}$ is left adjoint to $\tau$, for $A$ hereditary, thus 
$\Ext_A^1(T,\tau^{-1}T) 
 \simeq D\Hom_A(\tau^{-1}T,\tau T) 
 \simeq D\Hom_A(T,\tau^2 T) 
 \simeq \Ext_A^1(\tau T,T).$
The importance of the bimodule $\Ext_A^1(\tau T,T)$ has been stressed already in
Part II; I like to call it the ``magic'' bimodule for such a tilting process.
All the bimodule isomorphisms mentioned here should be of interest when dealing with
the magic bimodule $J$. In particular, when working with injective $\widetilde B$-modules, it seems to be convenient to know that
$DJ \simeq \Hom_A(T,\tau^2 T).$ }
(using this $\Ext^2$-bimodule has the advantage that it refers only to the
algebra $B$ itself, and not to $T$). It was Zhu Bin [Zh1] who stressed that 
cluster tilted algebras should be explored as semi-direct ring extensions.

Since this isomorphism is quite essential, let me sketch an elementary proof, without
reference to derived categories. Let $V$ be the universal extension of $\tau T$ by 
copies of $T$ from above, thus there is an exact sequence
$$
 0 \to \tau T \to V \to T^m \to 0 \tag{*}
$$
for some $m$, and $\Ext^1_A(T,V) = 0.$ Applying $\Hom_A(-,T)$ to $(*)$ shows that
$\Ext^1_A(V,T) \simeq \Ext^1_A(\tau T,T)$. Applying $\Hom_A(T,-)$ to $(*)$ yields the exact sequence
$$
 0 \to \Hom_A(T,V) \to \Hom_A(T,T^m) \to \Ext_A^1(T,\tau T) \to 0.
$$
This is an exact sequence of $B$-modules and $\Hom_A(T,T^m)$ is a free $B$-module, thus
we see that $\Hom_A(T,V)$ is a syzygy module for the $B$-module $\Ext^1_A(T,\tau T).$
But the latter means that
$$
 \Ext^2_B(\Ext_A^1(T,\tau T),{}_BB) \simeq \Ext_B^1(\Hom_A(T,V),{}_BB).
$$
The left hand side is nothing else than $\Ext^2_B(DB,B),$ since the $B$-module $DB$ and
$\Ext^1_A(T,\tau T)$ differ only by projective-injective direct summands. The right hand
side $\Ext^1_B(\Hom_A(T,V),\Hom_A(T,T))$ is the image of $\Ext_A^1(V,T)$ under the (exact) equivalence
$\Hom_A(T,-)\:\Cal T \to \Cal Y$ (here we use that $V$ belongs to $\Cal T$). This completes
the proof\note{Note that the isomorphy of $\Ext^2_B(DB,B)$ and $\Ext^1_A(\tau T,T)$
yields a proof for the implication (ii) $\implies$ (iv) mentioned in Part I. Since we know that $B$ has global dimension at most 2, the vanishing of $\Ext^2_B(DB,B)$
implies that $\Ext^2_B(X,Y) = 0$ for all $B$-modules $X,Y$, thus we also
see that (iv) $\implies$ (ii).}.
	\medskip
Now let us deal with the {\bf representations} of $\widetilde B$.
The $\widetilde B$-modules can be described as follows: they are pairs of the form
$(M,\gamma)$, where $M$ is a $B$-module, and $\gamma\:J\otimes_B M \to M$ is a $B$-linear map.
As we know, in $\mod B$ there is the splitting torsion pair $(\Cal Y, \Cal X)$
and it turns out that $J\otimes_B X = 0$ for $X\in \Cal X$, and that $J\otimes_B Y$ belongs to
$\Cal X$ for all $Y\in \Cal Y$. Let us consider a pair $(M,\gamma)$ in $\mod \widetilde B$
and write $M = X\oplus Y \oplus S$, with $X\in \Cal X,$ $Y\in \Cal Y'$ and $S\in \Cal S$.
Then the image of $\gamma$ is contained in $\Cal Y'$ and $Y\oplus S$ is contained in the 
kernel of $\gamma$ (in particular, $(S,0)$ is a direct summand of $(M, \gamma)$).

Note that $(\Cal Y, \Cal X)$ still is a torsion pair in $\mod \widetilde B$ (a 
module $(X\oplus Y,\gamma)$ with $X\in \Cal X$ and $Y\in \Cal Y$ has $(X,0)$ as torsion submodule,
has $(Y,0)$ as its torsion-free factor module, and the map $\gamma$ is the obstruction 
for the torsion submodule to split off). 
Let us draw the attention to a 
special feature of this torsion pair $(\Cal Y,\Cal X)$ in $\mod \widetilde B$: there 
exists an ideal, namely $J$, such that the modules annihilated by $J$
are just the modules in $\add(\Cal X, \Cal Y).$
	\medskip
Buan, Marsh and Reiten [BMR1] have shown that the category $\mod \widetilde B$ 
can be described in terms of $\mod A$ (via the corresponding cluster category). Let us
present such a description in detail.
We will use that $J = \Ext^1_A(T,\tau^{-1}T)$ (as explained above).
The algebra $\widetilde B$ has as $\Bbb Z$-covering the 
following (infinite dimensional) matrix algebra:
$$
\hbox{\beginpicture
\setcoordinatesystem units <0.4cm,0.4cm>
\put{} at 5 5
\put{} at 0 0
\plot 0.2 0  0 0  0 5  0.2 5 /
\plot 4.8 0  5 0  5 5  4.8 5 /
\multiput{$B$} at 1.5 3.5  2.5 2.5  3.5 1.5 /
\multiput{$J$} at 2.5 3.5  3.5 2.5 /
\setdots <1mm>
\plot 0.2 4.8 1 4 /
\plot 1.2 4.8 2 4 /
\plot 4 2  4.8 1.2 /
\plot 4 1  4.8 0.2 /
\put{$B_\infty =$} at -2 2.5
\endpicture},
$$
with $B$ on the main diagonal, $J$ directly above the main diagonal, and zeros elsewhere
(note that this algebra has no unit element in case $B \neq 0$). It turns out that
it is sufficient to determine the representations of the convex subalgebras of the form
$B_2 = \bmatrix B & J \cr
                0 & B \endbmatrix.$
We can write $B_2$-modules as columns 
$\left[\smallmatrix N \cr N' \endsmallmatrix\right]$
and use matrix multiplication, provided we have specified a map $\gamma\:J\otimes N' \to N$.
In the example considered ($B$ a square, with two zero relations), 
the algebras $B_\infty$ and $B_2$ are as follows:

$$
\hbox{\beginpicture
\setcoordinatesystem units <1cm,1cm>
\multiput{\beginpicture
\setcoordinatesystem units <.7cm,.7cm>
\put{} at 0 0
\put{} at 2 2 
\multiput{$\circ$} at 0 1  1 0  1 2  2 1 /
\arr{0.8 1.8}{0.2 1.2}
\arr{0.8 0.2}{0.2 0.8}
\arr{1.8 1.2}{1.2 1.8}
\arr{1.8 0.8}{1.2 0.2}

\setdots <.5mm>
\plot 0.7 1.6  1.25 1.6 /
\plot 0.7 0.4  1.25 0.4 /
\endpicture} at 0 0  2.6 0  5.2 0 /
\multiput{\beginpicture
\put{} at 0 0
\put{} at 2 2 
\arr{1.5 0.93}{0.5 0.93}
\arr{1.5 1.07}{0.5 1.07}
\setdots <.5mm>
\setquadratic
\plot 0.15 1.4  0.4 1.3  0.8 1.2 /
\plot 1.85 1.4  1.6 1.3  1.2 1.2 /
\plot 0.15 0.6  0.4 0.7  0.8 0.8 /
\plot 1.85 0.6  1.6 0.7  1.2 0.8 /
\endpicture} at  1.3 0  3.9 0 /

\put{\beginpicture
\put{} at 0 0
\put{} at 2 2 
\arr{1.5 0.93}{0.5 0.93}
\arr{1.5 1.07}{0.5 1.07}
\setdots <.5mm>
\setquadratic
\plot 1.85 1.4  1.6 1.3  1.2 1.2 /
\plot 1.85 0.6  1.6 0.7  1.2 0.8 /
\endpicture} at  -1.3 0  

\put{$B_\infty$} at -3.5 0
\plot -.7 -1  -.7 -1.2  3.3 -1.2 3.3 -1 /
\put{$B_2$} at 1.5 -1.5
\put{$\cdots$} at -2.2 0 
\put{$\cdots$} at 6.4 0 
\endpicture}
$$

In order to exhibit all the $B_2$-modules, we use the functor
$\Phi\:\mod A \to \mod B_2$ given by
$$
 \Phi(M) = \bmatrix \Ext^1_A(T,M)\cr
            \Hom_A(\tau^{-1}T,M) \endbmatrix, 
$$
with $\gamma\:\Ext^1_A(T,\tau^{-1}T)\otimes \Hom_A(\tau^{-1}T,M) \longrightarrow \Ext^1_A(T,M)$
being the canonical map of forming induced exact sequences (this is just the Yoneda multiplication)\note{The reader should recall that the 
functors $\Ext^1_A(T,-)$ and $\Hom_A(\tau^{-1}T,-)$ have been mentioned already in
Part I. These are the functors which provide the equivalences $\Cal F \simeq \Cal X$ and
$\Cal T' \simeq \Cal Y'$, respectively.
}.
Now $\Phi$ itself is
not faithful, since obviously $T$ is sent to 
zero\note{
The comparison with the Buan-Marsh-Reiten paper [BMR1] shows a slight deviation: 
The functor they use vanishes on the modules $\tau T$ and not on $T$ (and if we denote by
$T_i$ an indecomposable direct summand of $T$, then the image
of $T_i$ becomes an indecomposable projective $\widetilde B$-module). Instead of looking at
the functor $\Phi$, we could have worked with $\Phi'(M) = \bmatrix \Ext^1_A(\tau T,M)\cr
            \Hom_A(T,M) \endbmatrix,$ again taking for $\gamma$ the canonical map. 
This functor $\Phi'$ vanishes on $\tau T$. On the level
of cluster categories, the constructions corresponding to $\Phi$ and $\Phi'$ 
differ only by the Auslander-Reiten translation in the cluster category, and this is an
auto-equivalence of the cluster category. But as functors $\mod A \to \mod B_2$, 
the two functors $\Phi,\Phi'$ are quite different.
Our preference for the functor $\Phi$ has the following reason: the functor
$\Phi$ kills precisely $n = n(A)$ indecomposable $A$-modules, 
thus the number of indecomposable $\widetilde B$-modules which are not
contained in the image of $\Phi$ is also $n$, and these modules form a slice. This looks
quite pretty: the category $\mod \widetilde B$ is divided into the image of the functor
$\Phi$ and one additional slice.}. 
However, it induces a fully faithful
functor (which again will be denoted by $\Phi$):
$$
 \Phi\:\mod A/\langle T \rangle \longrightarrow \mod B_2,
$$
where $\mod A/\langle T \rangle$ denotes the factor category of $\mod A$ 
modulo the ideal of all maps which factor through $\add T$.
The image of the functor $\Phi$ is given by 
$$\left[\smallmatrix \Cal X\cr
               0 \endsmallmatrix\right] \under
\left[\smallmatrix 0\cr
               \Cal Y' \endsmallmatrix\right].
$$
In general, given module classes $\Cal K, \Cal L$ in $\mod R$, we write
$\Cal K\underkl \Cal L$ for the class of all $R$-modules $M$ with a submodule $K$ in 
$\Cal K$ such that $M/K$ belongs to $\Cal L$. Thus, we assert that the image of $\Phi$
is the class of the $\widetilde B$-modules 
$\left[\smallmatrix N\cr
               N' \endsmallmatrix\right]$ 
with $N \in \Cal X$ and $N'\in \Cal Y'$. (In order to see that 
$\Hom_A(\tau^{-1}T,M)\in \Cal Y'$, first note that $\Hom_A(\tau^{-1}T,M) =
\Hom_A(T,\tau M)$, thus this is a $B$-module in $\Cal Y.$ We further have $\Hom_A(T,\tau M)
= \Hom_A(T,t\tau M)$, where $t\tau M$ is the torsion submodule of $\tau M$.
If we assume that $\Hom_A(T,t\tau M)$ 
has an indecomposable submodule in $\Cal S$, say $\Hom_A(T,Q)$, where $Q$ is an indecomposable injective $A$-module, then we obtain a non-zero map $Q \to t\tau M
\subseteq \tau M$, since $\Hom_A(T,-)$ is fully faithful on $\Cal T$. However, the image of this
map is injective (since $A$ is hereditary) and $\tau M$ is indecomposable, thus $\tau M$
is injective, which is impossible). 

We want to draw a rough sketch of the shape of $\mod B_2$, in the same spirit as 
we have drawn a picture of $\mod B$ in Part I:
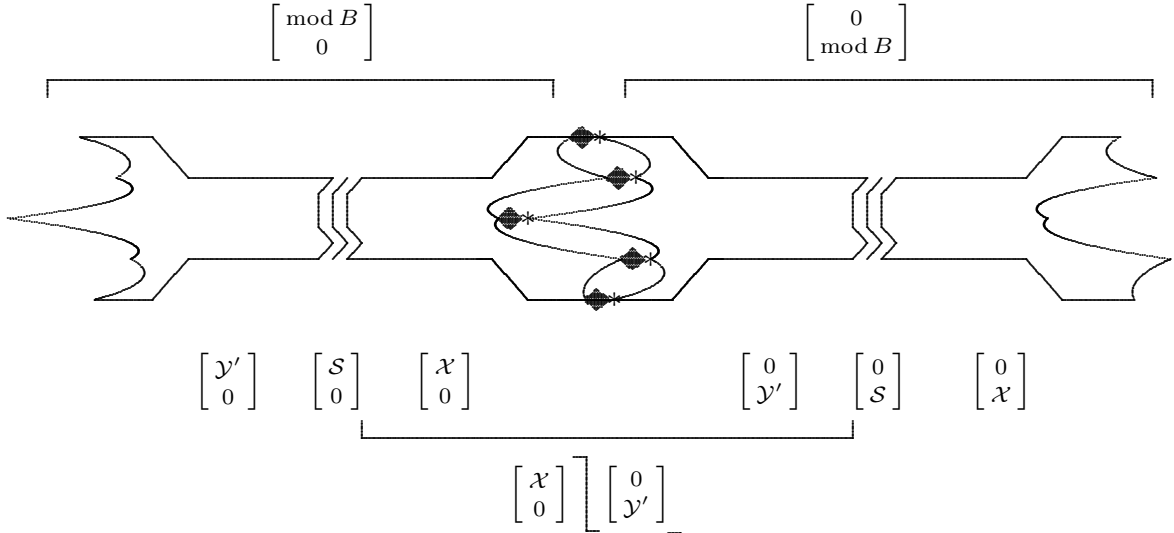
$$
\hbox{\beginpicture
\setcoordinatesystem units <.8cm,.9cm>
\put{\beginpicture 

\plot 2.5 4.5  2.7 4.7  2.5 4.9  2.5 5.3  2.7 5.5 /
\plot 2.5 4.5  4.5 4.5  5 4  7 4  7.5 4.5  9.5 4.5 / 
\plot 9.5 4.5  9.7 4.7  9.5 4.9  9.5 5.3  9.7 5.5 /
\plot 2.7 5.5  4.5 5.5  5 6  7 6  7.5 5.5  9.7 5.5 /
\plot 4.3 4.5  4.5 4.5  5 4  7 4  7.5 4.5  7.7 4.5 / 
\plot 4.3 5.5  4.5 5.5  5 6  7 6  7.5 5.5  7.7 5.5 /

\put{$*$} at 6 6
\put{$*$} at 6.5 5.5
\put{$*$} at 5 5
\put{$*$} at 6.7 4.5
\put{$*$} at 6.2 4
\setshadegrid span <.13mm>
\vshade 5.55 6 6 <,z,,> 5.75 5.85 6.15 <z,,,> 5.95 6 6 / 
\vshade 6.05 5.5 5.5 <,z,,> 6.25 5.35 5.65 <z,,,> 6.45 5.5 5.5 / 
\vshade 4.55 5 5 <,z,,> 4.75 4.85 5.15 <z,,,> 4.95 5 5 / 
\vshade 6.25 4.5 4.5 <,z,,> 6.45 4.35 4.65 <z,,,> 6.65 4.5 4.5 / 
\vshade 5.75 4 4 <,z,,> 5.95 3.85 4.15 <z,,,> 6.15 4 4 /

\setquadratic
\plot 6 6      6.65  5.75   6.5 5.5  /
\plot 6.5 5.5  6.55  5.25   5.0 5  /
\plot 5 5      6.55  4.75   6.7 4.5 /
\plot 6.7 4.5  6.85  4.25   6.2 4 /

\plot 5.6 6    5.45 5.75    6.1 5.5  /
\plot 6.1 5.5  4.65 5.25    4.6 5  /
\plot 4.6 5    4.85 4.75    6.3 4.5 /
\plot 6.3 4.5  5.85 4.25    5.8 4 /
\endpicture} at 0 0

\put{\beginpicture

\plot 9.3 0.5  9.5 0.7  9.3 0.9  9.3 1.3  9.5 1.5 /
\plot 9.5 0.5  9.7 0.7  9.5 0.9  9.5 1.3  9.7 1.5 /
\plot 6 2  7 2  7.5 1.5  9.5 1.5 /
\plot 6.2 0  7 0  7.5 0.5  9.3 0.5 /
\setquadratic

\plot 6 2      6.65  1.75     6.5 1.5  /
\plot 6.5 1.5  6.55  1.25    5 1  /
\plot 5 1    6.55  0.75   6.7 0.5 /
\plot 6.7 0.5  6.85  0.25   6.2 0 /

\setlinear
\put{} at 8 1
\put{} at 11 1

\endpicture} at -3.55 0

\put{\beginpicture

\plot 9.7 0.5  9.9 0.7  9.7 0.9  9.7 1.3  9.9 1.5 /
\put{} at 9.7 0.5
\put{} at 9.9 1.5

\endpicture} at 3.95 0
\put{\beginpicture

\plot 9.7 0.5  9.9 0.7  9.7 0.9  9.7 1.3  9.9 1.5 /

\setquadratic

\plot 13 2    12.85 1.75     13.5 1.5  /
\plot 13.5 1.5  12.05 1.25  12 1  /
\plot 12 1    12.25 0.75    13.7 0.5 /
\plot 13.7 0.5  13.25 0.25   13.2 0 /

\setlinear
\plot 9.9 1.5  11.7 1.5  12.2 2  13   2 / 
\plot 9.7 0.5  11.7 0.5  12.2 0  13.2 0 /
\put{} at 8 1
\put{} at 11 1

\endpicture} at 3.85 0
\put{$\left[\smallmatrix \Cal Y'\cr
               0 \endsmallmatrix\right]$} at -5 -2
\put{$\left[\smallmatrix \Cal S\cr
               0 \endsmallmatrix\right]$} at -3.5 -2
\put{$\left[\smallmatrix \Cal X\cr
               0 \endsmallmatrix\right]$} at -2 -2
\put{$\left[\smallmatrix 0\cr
               \Cal Y' \endsmallmatrix\right]$} at 2.5 -2
\put{$\left[\smallmatrix 0\cr
               \Cal S \endsmallmatrix\right]$} at 4 -2
\put{$\left[\smallmatrix 0\cr
               \Cal X \endsmallmatrix\right]$} at 5.7 -2
\plot -7.5 1.5  -7.5 1.7  -0.5 1.7   -0.5 1.5  /
\put{$\left[\smallmatrix \mod B \cr
               0 \endsmallmatrix\right]$} at -3.7 2.3
\plot   0.5 1.5   0.5 1.7  7.8  1.7  7.8 1.5 /
\put{$\left[\smallmatrix 0 \cr
               \mod B \endsmallmatrix\right]$} at 3.7 2.3
\plot   -3.15 -2.5   -3.15 -2.7  3.65 -2.7  3.65 -2.5 /
\put{$\left[\smallmatrix \Cal X\cr
               0 \endsmallmatrix\right] \under
\left[\smallmatrix 0\cr
               \Cal Y' \endsmallmatrix\right]$} at 0 -3.4

\endpicture}
$$
As we have mentioned, the middle part 
$\left[\smallmatrix \Cal X\cr
               0 \endsmallmatrix\right] \under
\left[\smallmatrix 0\cr
               \Cal Y' \endsmallmatrix\right]$
(starting with $\left[\smallmatrix \Cal X\cr
               0 \endsmallmatrix\right]$
and ending with $\left[\smallmatrix 0\cr
               \Cal Y' \endsmallmatrix\right]$)
is the image of the functor $\Phi$, thus this part of the category $\mod B_2$ is
equivalent to $\mod A/\langle T\rangle$. Note that this means that there are some 
small ``holes'' in this part, they are indicated by black lozenges; 
these holes correspond to the position $x$ in the Auslander-Reiten quiver of
$A$ which are given by the indecomposable direct summands $T_i$ of $T$ (and are
directly to the left of the small stars). 
	\medskip
It follows that $\mod \widetilde B$ has the form:
$$
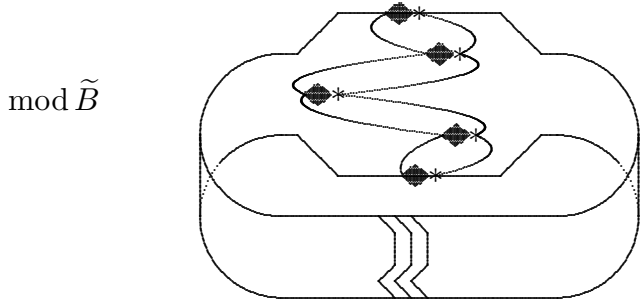
\hbox{\beginpicture
\setcoordinatesystem units <.9cm,.9cm>
\put{} at  0 2.5
\put{} at 13.2 6

\plot 4.3 4.5  4.5 4.5  5 4  7 4  7.5 4.5  7.7 4.5 / 
\plot 4.3 5.5  4.5 5.5  5 6  7 6  7.5 5.5  7.7 5.5 /

\put{$*$} at 6 6
\put{$*$} at 6.5 5.5
\put{$*$} at 5 5
\put{$*$} at 6.7 4.5
\put{$*$} at 6.2 4
\setshadegrid span <.13mm>
\vshade 5.55 6 6 <,z,,> 5.75 5.85 6.15 <z,,,> 5.95 6 6 / 
\vshade 6.05 5.5 5.5 <,z,,> 6.25 5.35 5.65 <z,,,> 6.45 5.5 5.5 / 
\vshade 4.55 5 5 <,z,,> 4.75 4.85 5.15 <z,,,> 4.95 5 5 / 
\vshade 6.25 4.5 4.5 <,z,,> 6.45 4.35 4.65 <z,,,> 6.65 4.5 4.5 / 
\vshade 5.75 4 4 <,z,,> 5.95 3.85 4.15 <z,,,> 6.15 4 4 / 

\setquadratic
\plot 6 6      6.65  5.75   6.5 5.5  /
\plot 6.5 5.5  6.55  5.25   5.0 5  /
\plot 5 5      6.55  4.75   6.7 4.5 /
\plot 6.7 4.5  6.85  4.25   6.2 4 /

\plot 5.6 6    5.45 5.75    6.1 5.5  /
\plot 6.1 5.5  4.65 5.25    4.6 5  /
\plot 4.6 5    4.85 4.75    6.3 4.5 /
\plot 6.3 4.5  5.85 4.25    5.8 4 /

\put{$\mod \widetilde B$} at 1.5 5

\setlinear
\circulararc  90 degrees from 3.3 4.5 center at 4.3 4.5
\circulararc -90 degrees from 3.3 4.5 center at 4.3 4.5
\circulararc  90 degrees from 3.3 3.5 center at 4.3 3.5
	\setdots <.5mm>
\circulararc -40 degrees from 3.3 3.5 center at 4.3 3.5
	\setsolid
\circulararc  50 degrees from 4.3 4.5 center at 4.3 3.5

\circulararc -90 degrees from 8.7 4.5 center at 7.7 4.5
\circulararc  90 degrees from 8.7 4.5 center at 7.7 4.5
\circulararc -90 degrees from 8.7 3.5 center at 7.7 3.5
	\setdots <.5mm>
\circulararc  40 degrees from 8.7 3.5 center at 7.7 3.5
	\setsolid
\circulararc -50 degrees from 7.7 4.5 center at 7.7 3.5

	\setsolid
\setlinear
\plot 3.3 4.5  3.3 3.5 /
\plot 8.7 4.5  8.7 3.5 /

\plot 4.3 2.5  7.7 2.5 /
\plot 4.3 3.5  7.7 3.5 /

\plot 5.7 2.5  5.5 2.7  5.7 2.9  5.7 3.3  5.5 3.5 /
\plot 5.9 2.5  5.7 2.7  5.9 2.9  5.9 3.3  5.7 3.5 /
\plot 6.1 2.5  5.9 2.7  6.1 2.9  6.1 3.3  5.9 3.5 /

\endpicture}
$$
Here, we have used the covering functor $\Pi\: \mod B_\infty \to \mod \widetilde B$ (or
better its restriction to $\mod B_2$): under this functor the subcategories
$\left[\smallmatrix \mod B\cr
               0 \endsmallmatrix\right]$ 
and $\left[\smallmatrix 0\cr
               \mod B \endsmallmatrix\right]$
are canonically identified. In particular, a fundamental domain for the covering
functor is given by the module classes 
$\left[\smallmatrix \Cal X\cr
               0 \endsmallmatrix\right] \under
\left[\smallmatrix 0\cr
               \Cal Y' \endsmallmatrix\right]$
and $\left[\smallmatrix 0 \cr
               \Cal S \endsmallmatrix\right]$.

This shows that $\mod \widetilde B$ decomposes into the modules in $\Cal X\underkl \Cal Y'$
(these are the $\widetilde B$-modules $N$ with a submodule $X \subseteq N$ in $\Cal X$, such
that $N/X$ belongs to $\Cal Y'$) on the one hand, and the modules in $\Cal S$ on the other hand. 
Under the functor $\Phi$, 
$\mod A/\langle T\rangle$ is embedded into $\mod \widetilde B$ with image the module class
$\Cal X\underkl \Cal Y'$. This is a controlled embedding (as defined in [R5]), with control class $\Cal S.$

The functor
$$
 \mod A @>\Phi>> \mod B_2 @>\Pi>> \mod \widetilde B
$$
has the following interesting property: only finitely many indecomposables are killed
by the functor (the indecomposable direct summands of $T$)
and there are only finitely many indecomposables (actually, the same number)
which are not in the image of the functor (the indecomposable modules in $\Cal S$).
Otherwise, it yields a bijection between 
indecomposables.

It should be noted that some of the strange phenomena of tilted algebras disappear
when passing to cluster tilted algebras. For example, the tunnel effect mentioned
above changes as follows: there still is the tunnel, but no longer does it connect two
separate regions; it now is a sort of by-pass for a single region. 
On the other hand, we should stress that the pictures which we have presented and which 
emphasise the existence of cyclic paths in $\mod \widetilde B$ are misleading in the
special case when $T$ is a slice module: in this case, $J =
\Ext_B^2(DB,B) = 0$, thus $\widetilde B = B$ is again hereditary. 
	\medskip
The cluster tilting theory has produced a lot of surprising results --- it even answered
some question which one did not dare to ask. For example, dealing with certain classes of algebras
such as special biserial ones, one observes that sometimes there do exist
indecomposable direct summands $X$ of the radical of  an indecomposable projective module
$P$, such that the Auslander-Reiten translate $\tau X$ is a direct summand 
of the socle factor module
of an indecomposable injective module $I$. Thus, in the Auslander-Reiten quiver of $\widetilde B$,
there are non-sectional paths of length 4 from $I$ to $P$
$$
 I \to \tau X \to X' \to X \to P.
$$
Is this configuration of interest? I did not think so,
but according to [BMR1], 
this configuration is a very typical one when dealing with cluster tilted algebras.

As an illustration, we show what happens in the non-regular components of our example $B_\infty$
(where $B$ is the square with two zero relations).
The upper line exhibits the part of the
quiver of $B_\infty$ which is needed as support for the modules shown below:
$$
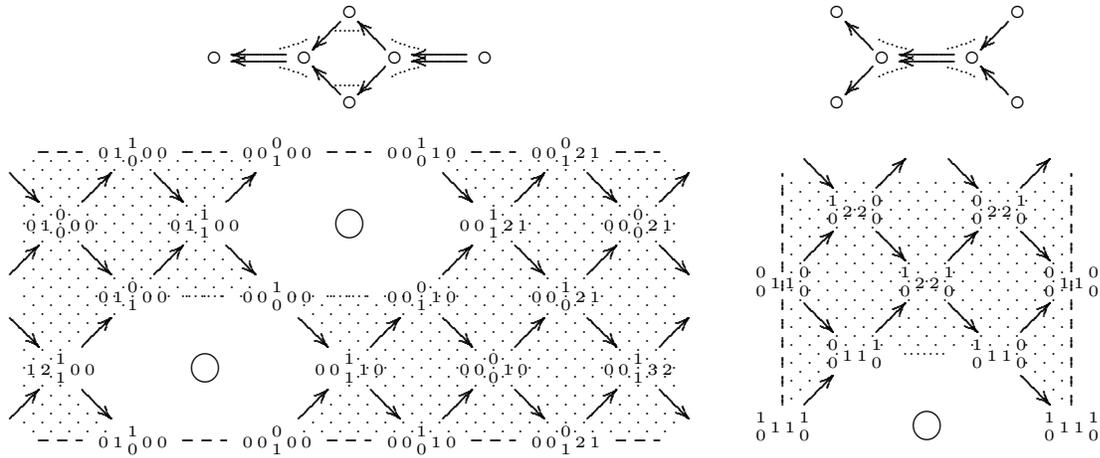
\hbox{\beginpicture
\setcoordinatesystem units <0.8cm,0.8cm>

\put{\beginpicture
\setcoordinatesystem units <.5cm,.5cm>
\put{} at 0 0
\put{} at 2 2 
\multiput{$\circ$} at 0 0  0 2  1 1  3 1  4 0  4 2 /
\arr{3.8 1.8}{3.2 1.2}
\arr{3.8 0.2}{3.2 0.8}
\arr{0.8 1.2}{0.2 1.8}
\arr{0.8 0.8}{0.2 0.2}

\arr{2.6 0.93}{1.4 0.93}
\arr{2.6 1.07}{1.4 1.07}

\setdots <.5mm>
\setquadratic
\plot 0.95 1.4  1.2 1.3  1.6 1.2 /
\plot 0.95 0.6  1.2 0.7  1.6 0.8 /

\plot 3.05 1.4  2.8 1.3  2.4 1.2 /
\plot 3.05 0.6  2.8 0.7  2.4 0.8 /

\endpicture} at 8 3.7

\put{\beginpicture
\setcoordinatesystem units <.5cm,.5cm>
\put{} at 0 0
\put{} at 2 2 
\multiput{$\circ$} at 0 1  1 0  1 2  2 1  -2 1  4 1 /
\arr{0.8 1.8}{0.2 1.2}
\arr{0.8 0.2}{0.2 0.8}
\arr{1.8 1.2}{1.2 1.8}
\arr{1.8 0.8}{1.2 0.2}

\setdots <.5mm>
\plot 0.7 1.6  1.25 1.6 /
\plot 0.7 0.4  1.25 0.4 /
\setsolid
\arr{3.6 0.93}{2.4 0.93}
\arr{3.6 1.07}{2.4 1.07}
\arr{-.4 0.93}{-1.6 0.93}
\arr{-.4 1.07}{-1.6 1.07}
\setdots <.5mm>
\setquadratic
\plot 1.95 1.4  2.2 1.3  2.6 1.2 /
\plot 1.95 0.6  2.2 0.7  2.6 0.8 /

\plot 0.05 1.4  -0.2 1.3  -.6 1.2 /
\plot 0.05 0.6  -0.2 0.7  -.6 0.8 /
\endpicture} at 0 3.7

\put{\beginpicture
\setcoordinatesystem units <0.8cm,0.8cm>
\put{} at 0 0
\put{} at 0 3
\put{$\vecta 011000 $} at 1 1
\put{$\vecta 000100 $} at 3 1  
\put{$\vecta 001010 $} at 5 1
\put{$\vecta 000121 $} at 7 1

\put{$\vecta 121100 $} at 0 2
\put{$\bigcirc$} at 2 2 
\put{$\vecta 001110 $} at 4 2 
\put{$\vecta 000010 $} at 6 2
\put{$\vecta 001132 $} at 8 2

\put{$\vecta 010100 $} at 1 3
\put{$\vecta 001000 $} at 3 3
\put{$\vecta 000110 $} at 5 3
\put{$\vecta 001021 $} at 7 3

\put{$\vecta 010000 $} at 0 4
\put{$\vecta 011100 $} at 2 4
\put{$\bigcirc$} at 4 4
\put{$\vecta 001121 $} at 6 4
\put{$\vecta 000021 $} at 8 4

\put{$\vecta 011000 $} at 1 5
\put{$\vecta 000100 $} at 3 5 
\put{$\vecta 001010 $} at 5 5
\put{$\vecta 000121 $} at 7 5
\arr{-.7 1.3} {-0.3 1.7} 
\arr{0.3 1.7} {0.7 1.3} 
\arr{3.3 1.3} {3.7 1.7} 
\arr{4.3 1.7} {4.7 1.3} 
\arr{5.3 1.3} {5.7 1.7} 
\arr{6.3 1.7} {6.7 1.3} 
\arr{7.3 1.3} {7.7 1.7} 
\arr{8.3 1.7} {8.7 1.3} 
\arr{-.7 2.7} {-0.3 2.3} 
\arr{0.3 2.3} {0.7 2.7} 
\arr{3.3 2.7} {3.7 2.3} 
\arr{4.3 2.3} {4.7 2.7} 
\arr{5.3 2.7} {5.7 2.3} 
\arr{6.3 2.3} {6.7 2.7}
\arr{7.3 2.7} {7.7 2.3} 
\arr{8.3 2.3} {8.7 2.7}
\arr{-.7 3.3} {-0.3 3.7} 
\arr{0.3 3.7} {0.7 3.3} 
\arr{1.3 3.3} {1.7 3.7} 
\arr{2.3 3.7} {2.7 3.3} 
\arr{5.3 3.3} {5.7 3.7} 
\arr{6.3 3.7} {6.7 3.3} 
\arr{7.3 3.3} {7.7 3.7} 
\arr{8.3 3.7} {8.7 3.3} 
\arr{-.7 4.7} {-0.3 4.3} 
\arr{0.3 4.3} {0.7 4.7} 
\arr{1.3 4.7} {1.7 4.3} 
\arr{2.3 4.3} {2.7 4.7} 
\arr{5.3 4.7} {5.7 4.3} 
\arr{6.3 4.3} {6.7 4.7} 
\arr{7.3 4.7} {7.7 4.3} 
\arr{8.3 4.3} {8.7 4.7} 

\put{} at 0 2
\setdashes<3pt>
\plot -.3 1  0.3 1 /
\plot 1.7 1  2.3 1 /
\plot 3.7 1  4.3 1 /
\plot 5.7 1  6.3 1 /
\plot 7.7 1  8.3 1 /

\plot -0.3 5  0.3 5 /
\plot 1.7 5  2.3 5 /
\plot 3.7 5  4.3 5 /
\plot 5.7 5  6.3 5 /
\plot 7.7 5  8.3 5 /

\setdots<2pt>
\plot 1.7 3  2.3 3 /
\plot 3.7 3  4.3 3 /
\setshadegrid span <1mm>
\hshade 1 -.6 1 <,,,z> 2 -.6 0   <,,z,z> 2.9 -.6 1 
   <,,z,z>  3 -.6 3 <,,z,z> 4 -.6 2 <,,z,> 5 -.6 3 / 
\hshade  1 3 8.6  <,,,z> 2 4 8.6   <,,z,z> 3 3 8.6 
   <,,z,z>  3.1 5 8.6 <,,z,z> 4 6 8.6 <,,z,> 5 5 8.6 / 
\endpicture} at 0 0

\put{\beginpicture
\setcoordinatesystem units <0.8cm,0.8cm>
\put{} at 0 0
\put{} at 0 3
\put{$\vectb 101 110 $} at 0 0
\put{$\vectb 101 110 $} at 4 0
\put{$\vectb 001 110 $} at 1 1
\put{$\vectb 101 100 $} at 3 1
\put{$\vectb 001 100 $} at 0 2
\put{$\vectb 102 210 $} at 2 2
\put{$\vectb 001 100 $} at 4 2 
\put{$\vectb 102 200 $} at 1 3
\put{$\vectb 002 210 $} at 3 3

\arr{0.3 0.3} {0.7 0.7} 
\arr{3.3 0.7} {3.7 0.3} 
\arr{0.3 1.7} {0.7 1.3} 
\arr{1.3 1.3} {1.7 1.7} 
\arr{2.3 1.7} {2.7 1.3} 
\arr{3.3 1.3} {3.7 1.7} 
\arr{0.3 2.3} {0.7 2.7} 
\arr{1.3 2.7} {1.7 2.3} 
\arr{2.3 2.3} {2.7 2.7} 
\arr{3.3 2.7} {3.7 2.3} 
\arr{0.3 3.7} {0.7 3.3} 
\arr{1.3 3.3} {1.7 3.7} 
\arr{2.3 3.7} {2.7 3.3} 
\arr{3.3 3.3} {3.7 3.7} 
\put{$\bigcirc$} at 2 0 
\put{} at 0 2
\setdashes<3pt>
\plot 0 0.3  0 3.5 /
\plot 4 0.3  4 3.5 /
\setdots<2pt>
\plot 1.7 1  2.3 1 /
\setshadegrid span <1mm>
\vshade 0 0 3.5 <,z,,> 1 1 3.5 <z,z,,> 3 1 3.5 <z,,,>  4 0 3.5 / 

\endpicture} at 8 0.1
\endpicture}
$$
For both components, the dashed boundary 
lines have to be identified. In this way, the right picture with the vertical
identification yields what is called a tube, the left picture gives a kind of horizontal hose.
In contrast to the tube with its mouth, the hose extends in both directions indefinitely.
The big circles
indicate the position of the modules $T_i$ in the corresponding components of $\mod A$,
these are the modules which are killed by the functor $\Phi.$ In both components we
find non-sectional paths of length 4 from an indecomposable injective $B_\infty$-module
$I$ to an indecomposable projective $B_\infty$-module $P$ such that the simple modules
$\soc I$ and $\top P$ 
 are identified under the covering functor $\Pi$.  
	\medskip
We also want to use this example in order to illustrate the fact that 
the image of $\Phi$ in $\mod \widetilde B$ is complemented by a slice $\Cal S$:
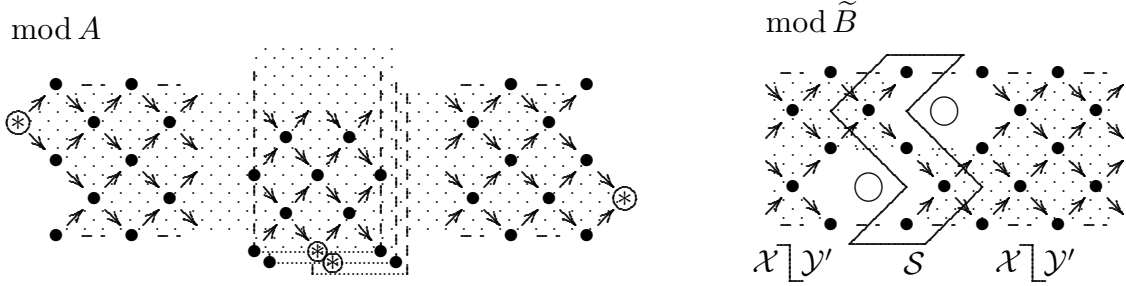
$$
\hbox{\beginpicture
\setcoordinatesystem units <0.42cm,0.42cm>

\put{\beginpicture
\put{} at 0 0
\put{} at 14 6 
\put{$\bullet$} at -1 1
\put{$\bullet$} at 1 1

\put{$\bullet$} at 11 1  
\put{$\bullet$} at 13 1

\put{$\bullet$} at 0 2
\put{$\bullet$} at 2 2

\put{$\bullet$} at 10 2 
\put{$\bullet$} at 12 2
\put{$\ast$} at 14 2
\circulararc 360 degrees from 14.3 2 center at 14 2

\put{$\bullet$} at -1 3
\put{$\bullet$} at 1 3

\put{$\bullet$} at 11 3
\put{$\bullet$} at 13 3

\put{$\ast$} at -2 4
\circulararc 360 degrees from -2.3 4 center at -2 4  

\put{$\bullet$} at 0 4
\put{$\bullet$} at 2 4

\put{$\bullet$} at 10 4
\put{$\bullet$} at 12 4

\put{$\bullet$} at -1 5
\put{$\bullet$} at 1 5

\put{$\bullet$} at 11 5
\put{$\bullet$} at 13 5
\arr{-.7 1.3} {-0.3 1.7} 
\arr{0.3 1.7} {0.7 1.3} 
\arr{1.3 1.3} {1.7 1.7} 
\arr{2.3 1.7} {2.7 1.3}
 
\arr{9.3 1.3} {9.7 1.7} 
\arr{10.3 1.7} {10.7 1.3} 
\arr{11.3 1.3} {11.7 1.7} 
\arr{12.3 1.7} {12.7 1.3} 
\arr{13.3 1.3} {13.7 1.7} 
\arr{-.7 2.7} {-0.3 2.3} 
\arr{0.3 2.3} {0.7 2.7} 
\arr{1.3 2.7} {1.7 2.3} 
\arr{2.3 2.3} {2.7 2.7} 

\arr{9.3 2.7} {9.7 2.3} 
\arr{10.3 2.3} {10.7 2.7} 
\arr{11.3 2.7} {11.7 2.3} 
\arr{12.3 2.3} {12.7 2.7}
\arr{13.3 2.7} {13.7 2.3} 
\arr{-1.7 3.7} {-1.3 3.3} 
\arr{-.7 3.3} {-0.3 3.7} 
\arr{0.3 3.7} {0.7 3.3} 
\arr{1.3 3.3} {1.7 3.7} 
\arr{2.3 3.7} {2.7 3.3} 

\arr{9.3 3.3} {9.7 3.7} 
\arr{10.3 3.7} {10.7 3.3} 
\arr{11.3 3.3} {11.7 3.7} 
\arr{12.3 3.7} {12.7 3.3} 
\arr{-1.7 4.3} {-1.3 4.7} 
\arr{-.7 4.7} {-0.3 4.3} 
\arr{0.3 4.3} {0.7 4.7} 
\arr{1.3 4.7} {1.7 4.3} 
\arr{2.3 4.3} {2.7 4.7} 

\arr{9.3 4.7} {9.7 4.3} 
\arr{10.3 4.3} {10.7 4.7} 
\arr{11.3 4.7} {11.7 4.3} 
\arr{12.3 4.3} {12.7 4.7} 

\setdashes<3pt>
\plot -.3 1  0.3 1 /
\plot 1.7 1  2.3 1 /

\plot 9.7 1  10.3 1 /
\plot 11.7 1  12.3 1 /

\plot -0.3 5  0.3 5 /
\plot 1.7 5  2.3 5 /

\plot 9.7 5  10.3 5 /
\plot 11.7 5  12.3 5 /

\put{\beginpicture
\setcoordinatesystem units <0.35cm,0.42cm>

\put{} at 0 0
\put{} at 4 6 
\setdashes<3pt>
\plot 0 0  0 5 /
\plot 4 0  4 5 /
\setsolid
\put{$\ast$} at 2 0
\circulararc 360 degrees from 2.3 0 center at  2 0 

\multiput{$\bullet$} at 0 0       4 0 
                        1 1  3 1
                        0 2  2 2  4 2
                        1 3  3 3 /
\arr{0.3 0.3} {0.7 0.7} 
\arr{1.3 0.7} {1.7 0.3} 
\arr{2.3 0.3} {2.7 0.7} 
\arr{3.3 0.7} {3.7 0.3} 
\arr{0.3 1.7} {0.7 1.3} 
\arr{1.3 1.3} {1.7 1.7} 
\arr{2.3 1.7} {2.7 1.3}
\arr{3.3 1.3} {3.7 1.7} 
\arr{0.3 2.3} {0.7 2.7} 
\arr{1.3 2.7} {1.7 2.3} 
\arr{2.3 2.3} {2.7 2.7} 
\arr{3.3 2.7} {3.7 2.3} 
\arr{0.3 3.7} {0.7 3.3} 
\arr{1.3 3.3} {1.7 3.7} 
\arr{2.3 3.7} {2.7 3.3}
\arr{3.3 3.3} {3.7 3.7} 
\setdots <.5mm>
\plot 0 0  4 0 /
\endpicture} at 5.9 3.5 

\put{\beginpicture
\setcoordinatesystem units <0.35cm,0.42cm>
\put{} at 0 0
\put{} at 4 6 
\setdashes<3pt>
\plot 4 0  4 5 /
\plot 0 0  0 .5 /
\setsolid
\multiput{$\bullet$} at 0 0    4 0 /
\put{$\ast$} at 2 0
\circulararc 360 degrees from 2.3 0 center at 2 0


\setdots <.5mm>
\plot 0 0  4 0 /
\endpicture} at 6.3 3.2 

\put{\beginpicture
\setcoordinatesystem units <0.35cm,0.42cm>
\put{} at 0 0
\put{} at 4 6 
\setdashes<3pt>
\plot 4 0  4 5 /
\plot 1 0  1 .5 /
\setsolid

\setdots <.5mm>
\plot 1 0  4 0 /
\endpicture} at 6.6 3



\setshadegrid span <1mm>
\hshade 1 -1  13 <,,,z> 2 0 14  <,,z,z> 4 -2  12 <,,z,>  5 -1 13 /
\hshade .5 4 8  1 4 8 /
\hshade 4.8 4 8  6 4 8 /

\put{$\mod A$} at -1 6.5

\endpicture} at 0 0  

\put{\beginpicture
\put{} at 0 0
\put{} at 8 0
\put{$\bullet$} at 1 1
\put{$\bullet$} at 3 1  
\put{$\bullet$} at 5 1
\put{$\bullet$} at 7 1

\put{$\bullet$} at 0 2
\put{$\bigcirc$} at 2 2 
\put{$\bullet$} at 4 2 
\put{$\bullet$} at 6 2
\put{$\bullet$} at 8 2

\put{$\bullet$} at 1 3
\put{$\bullet$} at 3 3
\put{$\bullet$} at 5 3
\put{$\bullet$} at 7 3

\put{$\bullet$} at 0 4
\put{$\bullet$} at 2 4
\put{$\bigcirc$} at 4 4
\put{$\bullet$} at 6 4
\put{$\bullet$} at 8 4

\put{$\bullet$} at 1 5
\put{$\bullet$} at 3 5 
\put{$\bullet$} at 5 5
\put{$\bullet$} at 7 5
\arr{-.7 1.3} {-0.3 1.7} 
\arr{0.3 1.7} {0.7 1.3} 
\arr{3.3 1.3} {3.7 1.7} 
\arr{4.3 1.7} {4.7 1.3} 
\arr{5.3 1.3} {5.7 1.7} 
\arr{6.3 1.7} {6.7 1.3} 
\arr{7.3 1.3} {7.7 1.7} 
\arr{8.3 1.7} {8.7 1.3} 
\arr{-.7 2.7} {-0.3 2.3} 
\arr{0.3 2.3} {0.7 2.7} 
\arr{3.3 2.7} {3.7 2.3} 
\arr{4.3 2.3} {4.7 2.7} 
\arr{5.3 2.7} {5.7 2.3} 
\arr{6.3 2.3} {6.7 2.7}
\arr{7.3 2.7} {7.7 2.3} 
\arr{8.3 2.3} {8.7 2.7}
\arr{-.7 3.3} {-0.3 3.7} 
\arr{0.3 3.7} {0.7 3.3} 
\arr{1.3 3.3} {1.7 3.7} 
\arr{2.3 3.7} {2.7 3.3} 
\arr{5.3 3.3} {5.7 3.7} 
\arr{6.3 3.7} {6.7 3.3} 
\arr{7.3 3.3} {7.7 3.7} 
\arr{8.3 3.7} {8.7 3.3} 
\arr{-.7 4.7} {-0.3 4.3} 
\arr{0.3 4.3} {0.7 4.7} 
\arr{1.3 4.7} {1.7 4.3} 
\arr{2.3 4.3} {2.7 4.7} 
\arr{5.3 4.7} {5.7 4.3} 
\arr{6.3 4.3} {6.7 4.7} 
\arr{7.3 4.7} {7.7 4.3} 
\arr{8.3 4.3} {8.7 4.7} 

\put{} at 0 2
\setdashes<3pt>
\plot -.3 1  0.3 1 /
\plot 1.7 1  2.3 1 /
\plot 3.7 1  4.3 1 /
\plot 5.7 1  6.3 1 /
\plot 7.7 1  8.3 1 /

\plot -0.3 5  0.3 5 /
\plot 1.7 5  2.3 5 /
\plot 3.7 5  4.3 5 /
\plot 5.7 5  6.3 5 /
\plot 7.7 5  8.3 5 /

\setdots<2pt>
\plot 1.7 3  2.3 3 /
\plot 3.7 3  4.3 3 /
\setshadegrid span <1mm>
\hshade 1 -.6 1 <,,,z> 2 -.6 0   <,,z,z> 2.9 -.6 1 
   <,,z,z>  3 -.6 3 <,,z,z> 4 -.6 2 <,,z,> 5 -.6 3 / 
\hshade  1 3 8.6  <,,,z> 2 4 8.6   <,,z,z> 3 3 8.6 
   <,,z,z>  3.1 5 8.6 <,,z,z> 4 6 8.6 <,,z,> 5 5 8.6 / 
\setsolid
\plot 1.5 .5  3.5 .5  5 2  3 4  4.5 5.5 /
\plot         1.5 .5  3 2  1 4  2.5 5.5  4.5 5.5 /
\put{$\Cal X \underkl \Cal Y'$} at   0 0
\put{$\Cal S$} at 3.2 0
\put{$\Cal X \underkl \Cal Y'$} at 6.4 0

\put{$\mod \widetilde B$} at 0.5 6.5

\endpicture} at 16 0
\endpicture}
$$

When looking at the non-sectional paths from $I$ to $P$ of length 4, 
where $I$ is an indecomposable injective $\widetilde B$-module, $P$ an indecomposable 
projective $\widetilde B$-module such that $S = \soc I \simeq \top P$, 
one should
be aware that the usual interest lies in paths from $P$ to $I$. Namely, there is the
so called ``hammock'' for the simple module $S$, 
dealing with pairs of maps of the form
$P \to M \to I$ with composition having image $S$
(and $M$ indecomposable). 
$$
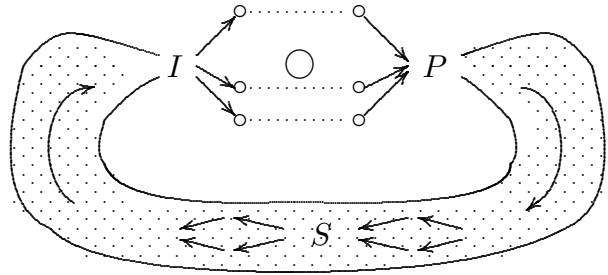
\hbox{\beginpicture
\setcoordinatesystem units <.6cm,.6cm>

\put{} at 0 0
\put{} at 0 -4
\put{$I$} at -0.2 0
\put{$P$} at 4.6 0
\multiput{$\circ$} at 1 -1  1 -0.4 1 1  3.2 -1 3.2 -0.4 3.2 1 /
\arr{0.2 0.2}{0.9 0.9}
\arr{0.2 0}{0.9 -0.4}
\arr{0.2 -0.2}{0.9 -0.9}

\arr{3.3 0.9}{4.1 0.1}
\arr{3.3 -0.4}{4.1 0}
\arr{3.3 -0.9}{4.1 -0.1}
\setdots <1mm>
\plot 1.2 1     3 1 /
\plot 1.2 -0.4  3 -0.4 /
\plot 1.2 -1    3 -1 /
\setshadegrid span <1mm>
\hshade -3.8  2 2 <,,,z> -3.5 -2 6 <,,z,> -2.5 -3.3 7.5 / 
\hshade -2.6 -3.3 -1 <,,,z> -1 -3.3 -1.7 <,,z,z> 0 -3 -1 <,,z,> 0.5 -2.5 -1.7 / 
\hshade -2.6  5.5 7.8 <,,,z>  -1 6.5 7.8 <,,z,z> 0 5 7.2 <,,z,> 0.5 5 7 / 

\setsolid

\setquadratic
\plot  5.1 0.2   6 .5     6.5  .6   7. .5    7.5 0    7.7  -2  7 -3  5 -3.6  2 -3.8 /
\plot    -.5 0.2   -1.5 .5     -2 .6  -2.5  .5   -3 0   -3.2 -2  -2.5 -3  -1 -3.6 2 -3.8 /

\plot  5.1 -0.2  5.5 -.5   6 -1    6.1 -1.5       6 -1.8    5 -2.4  2 -2.5   /
\plot -.5 -0.2   -1 -.5   -1.5 -1  -1.6 -1.5       -1.5 -1.8   -.5 -2.4  2 -2.5   /

\ellipticalarc axes ratio 1:1.5 -160 degrees from 6.2 -.4 center at 6.2 -1.5  
\ellipticalarc axes ratio 1:1.5  160 degrees from -1.8 -.4 center at -1.8 -1.5  
\arr{-1.8 -.4}{-1.7 -.38}
\arr{6.5 -2.5}{6.3 -2.65}

\arr{4.1 -3.4}{3.2 -3.2}
\arr{4.1 -2.8}{3.2 -3.}
\arr{5.1 -3.2}{4.3 -3.4}
\arr{5.1 -3.0}{4.3 -2.8}

\arr{1.8 -3.2}{.9 -3.4}
\arr{1.8 -3.0}{.9 -2.8}
\arr{.7 -3.4}{-.1 -3.2}
\arr{.7 -2.8}{-.1 -3.}

\put{$S$} at 2.5 -3.1 

\put{$\bigcirc$} at 2.1 0.05

\endpicture}
$$
Taking into account not only the hammock, but also the non-sectional 
paths of length 4 from $I$ to $P$ leads to a kind of organized round trip. 
Since the simple module $S$ has no self-extension, it is the only indecomposable
module $M$ such that $\Hom_{\widetilde B}(P',M) = 0$, for any indecomposable
projective $\widetilde B$-module $P'\not\simeq P.$ We will return to this hammock
configuration $(P,S,I)$ later. 
	\medskip

Readers familiar with the literature will agree that despite
of the large number of papers devoted to questions in the representation theory of
artin algebras, 
only few classes of artin algebras are known where there is a
clear description of the module categories\note{Say 
in the same way as the module categories of hereditary artin
algebras are described. We consider here algebras which may be wild, thus
we have to be cautious of what to expect from a ``clear description''.}. 
The new developments outlined here show that the cluster tilted algebras 
are such a class: As for the hereditary artin algebras, the description of the
module category is again given by the root system of a Kac-Moody
Lie algebra. 

Keller and Reiten [KR] have shown, that {\it cluster tilted algebras
are Gorenstein algebras of Gorenstein dimension at most 1.} This is a very
remarkable assertion! The proof uses in an
essential way (generalizations of) cluster categories, and provides further
classes of Gorenstein algebras of Gorenstein dimension at most 1. Note that
the cluster tilted algebra $\widetilde B$ is hereditary only in case $\widetilde B =
B$, thus only for $T$ a slice module. There are examples where $\widetilde B$
is self-injective (for example for $B = kQ/\langle\rho\rangle$ with $Q$ the
linearly directed $\Bbb A_3$-quiver and $\rho$ the path of length 2). In general,
$\widetilde B$ will be neither hereditary nor self-injective.
	\bigskip
{\bf The Complex $\Sigma'_A$.} We have mentioned in part I that the simplicial complex 
$\Sigma_A$ of tilting modules always has a non-empty boundary (for $n(A) \ge 2$). 
Now the cluster theory
provides a recipe for embedding this simplicial complex in a slightly larger one without
boundary. Let me introduce here this complex $\Sigma'_A$ directly in terms of $\mod A$,
using a variation of the work of Marsh, Reineke und Zelevinsky 
[MRZ]\note{The title of the paper refers to ``associahedra'': in the case of the path algebra 
of a quiver of type $\Bbb A_n$, the dual of the simplicial complex $\Sigma'_A$ is an
associahedron (or Stasheff polytope). For quivers of type $\Bbb B_n$ and $\Bbb C_n$ one 
obtains a Bott-Taubes cyclohedron.}.
It is obtained from $\Sigma_A$ by just adding 
$n = n(A)$ vertices, and of course further simplices. 
Recall that a {\it Serre 
subcategory}\note{The Serre subcategories are nothing else then the subcategories of the form
$\mod A/AeA$, where $e$ in an idempotent of $A$. 
} 
$\Cal U$ of an abelian category is a subcategory which is closed under submodules, factor modules and extensions; 
thus in case we deal with a length category such as $\mod A$, then $\Cal U$ is specified by 
the simple modules contained in $\Cal U$ (an object belongs to $\Cal U$ if and only if its
composition factors lie in $\Cal U$). In particular, for a simple $A$-module $S$, let us denote
by $(-S)$ the subcategory of all $A$-modules which do not have $S$ as a composition
factor. Any Serre subcategory is the intersection of such subcategories.

Here is the definition of $\Sigma'_A$: As simplices take the pairs $(M,\Cal U)$
where $\Cal U$ is a Serre subcategory of $\mod A$ and $M$ is (the isomorphism class of)
a basic module in $\Cal U$
without self-extensions; write $(M',\Cal U') \le (M,\Cal U)$ provided $M'$ is a direct
summand of $M$ and $\Cal U' \supseteq \Cal U$ (note the reversed order!).
Clearly\note{In the same way, we may identify the set of simplicies of the form $(M,\Cal U)$
with $\Cal U$ fixed, as $\Sigma_{A/AeA}$, where $\Cal U = \mod A/AeA.$ In this way, we
see that $\Sigma'_A$ can be considered as a union of all the simplicial
complexes $\Sigma_{A/AeA}$. 
}, 
$\Sigma_A$ can be considered as a subcomplex of $\Sigma'_A$, namely as the set
of all pairs $(M,\mod A)$. 

There are two kinds of vertices of $\Sigma'_A$,  namely those of the form
$(E,\mod A)$ with $E$ an exceptional $A$-module (these are the vertices belonging to
$\Sigma_A)$, and those of the form $(0,(-S))$ with $S$ simple. It is fair to say 
that the latter ones are indexed by the ``negative 
 simple roots''; of course these
are the vertices which do not belong to $\Sigma_A$. 
Given a simplex $(M,\Cal U)$, its vertices 
are the elements $(E,\mod A)$,
where $E$ is an indecomposable direct summand of $M$, and the elements $(0,(-S))$,
where $\Cal U \subseteq (-S).$ 
The $(n-1)$-simplices are those of the form $(M,\Cal U),$ where $M$ is
a basic tilting module in $\Cal U.$ 
The vertices outside $\Sigma_A$ belong to one $(n-1)$-simplex, namely to $(0,\{0\})$. 
The $(n-2)$-simplices are of the form $(M,\Cal U)$, where $M$ is an
almost complete partial tilting module for $\Cal U.$ If it is sincere in $\Cal U,$
there are precisely two complements in $\Cal U$. If it is not sincere in $\Cal U$, then
there is only one complement in $\Cal U$, but there also is a simple module $S$ such that
$X$ belongs to $(-S)$, thus $X$ is a tilting module for $\Cal U\cap (-S).$ 
This shows that any $(n-2)$-simplex belongs to precisely two
$(n-1)$-simplices. 

As an example, we consider again the path
algebra $A$ of the quiver $\circ \leftarrow \circ \leftarrow \circ$.
The simplicial complex $\Sigma'_A$ is a 2-sphere and looks as follows
(considering the 2-sphere as the 1-point compactification of the real plane):
$$
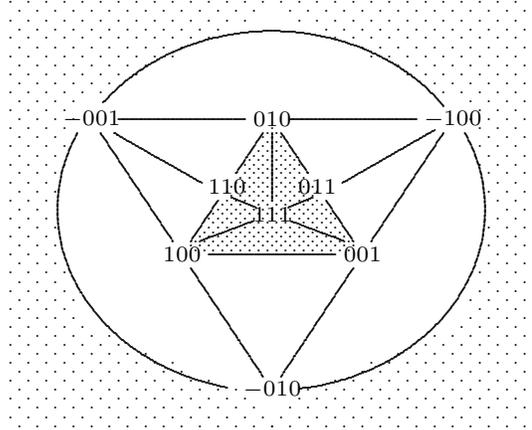
\hbox{\beginpicture
\setcoordinatesystem units <1cm,1.5cm>

\put{} at 0 0
\put{} at 0 2.7
\put{$\vectd 100 $} at 1 1
\put{$\vectd 110 $} at 1.5 1.5 
\put{$\vectd 010 $} at 2 2 
\put{$\vectd 111 $} at 2 1.3
\put{$\vectd 011 $} at 2.5 1.5
\put{$\vectd 001 $} at 3 1
\plot 1.1 1.1  1.4 1.4 /
\plot 1.6 1.6  1.9 1.9 /
\plot 2.1 1.9  2.4 1.6 /
\plot 2.6 1.4  2.9 1.1 /
\plot 1.3 1  2.7 1 /
\plot 2 1.9  2 1.4 /

\plot 1.55 1.42  1.82 1.35 /
\plot 2.2 1.25  2.8 1.1 /

\plot 2.45 1.42  2.18 1.35 /
\plot 1.8 1.25  1.2 1.1 /

\setshadegrid span <.5mm>
\hshade 1 1 3  2 2 2  / 

\put{$\vectd {-1}00 $} at 4 2
\put{$\vectd {-0}10 $} at 2 0
\put{$\vectd {-0}01 $} at 0 2 
\plot 0.1 1.9  0.9 1.1 /
\plot 1.1 0.9  1.9 0.1 /
\plot 2.1 0.1  2.9 0.9 /
\plot 3.1 1.1  3.9 1.9 /
\plot 0.3 2  1.7 2 /
\plot 2.2 2  3.6 2 /

\plot 0.25 1.9 1.22 1.53 / 
\plot 3.8 1.92 2.78 1.53 /

\ellipticalarc axes ratio 1.18:1 108 degrees from 2.3 0.01 center at 2 1.333
\ellipticalarc axes ratio 1.18:1 108 degrees from 3.9 2.1  center at 2 1.333
\ellipticalarc axes ratio 1.18:1 -103 degrees from 1.5 0.02 center at 2 1.333
\setshadegrid span <1mm>
\hshade 
-0.3 -1 2.1   <,,,z>  
-0.1 -1 2   <,,z,z>  
0.5  -1 -.1 <,,z,z>
1.3  -1 -.4 <,,z,z>
2.1  -1 -.1 <,,z,z>
2.8  -1 2   <,,z,> 
2.9  -1 2 / 

\hshade 
-0.3  1.9    5 <,,,z>  
-0.1 2    5 <,,z,z>  
0.5  4.1  5 <,,z,z>
1.3  4.4  5 <,,z,z>
2.1  4.1  5 <,,z,z>
2.8  2    5 <,,z,> 
2.9  2    5 / 

\endpicture}
$$
Here, the vertex $(0,(-S))$ is labeled as $-\bdim S$.
We have shaded the subcomplex $\Sigma_A$ (the triangle in the middle)
as well as the $(n-1)$-simplex $(0,\{0\})$ (the outside). 
	
Consider now a reflection functor $\sigma_i$, where $i$ is a sink, say. We obtain an
embedding of $\Sigma_{\sigma_iA}$ into $\Sigma'_A$ as follows: 
There are the exceptional $\sigma_iA$-modules of the form $\sigma_iE$ 
with $E$ an exceptional $A$-module, different from the simple $A$-module $S(i)$ concentrated
at the vertex $i$, and 
in between these modules $\sigma_iE$ the simplex structure
is the same as in between the modules $E$. In addition, there is the simple $\sigma_iA$-module
$S'(i)$ again concentrated at $i$. Now we know that $E$ has no composition factor
$S(i)$ if and only if 
$\Ext^1_{\sigma_iA}(S'(i),\sigma_iM) = 0.$ This shows that the simplex structure of 
$\Sigma_A$ involving $(0,(-S(i))$ and  vertices of the form $(E,\mod A)$ is the same
as the simplex structure of $\Sigma_{\sigma_iA}$ in the vicinity of $(S'(i),\mod\sigma_iA)$.
	\medskip
We may consider the simplicial complex $\Sigma'_A$ as a subset of the real $n$-dimensional
space $K_0(A)\otimes \Bbb R$, where $n = n(A)$,
namely as a part of the corresponding unit $(n-1)$-sphere,
with all the $(n-1)$-simplices defined by $n$ linear inequalities. In case $A$ is
representation-finite, we deal with the $(n-1)$-sphere itself, otherwise with a proper
subset. For example, in the case of the path algebra $A$ of the quiver 
$\circ \leftarrow \circ \leftarrow \circ$, the inequalities are $\phi_1 \ge 0,\ 
\phi_2 \ge 0,\  \phi_3 \ge 0$, where $\phi_1,\phi_2,\phi_3$ 
are the linear forms inserted in the
corresponding triangle:
$$
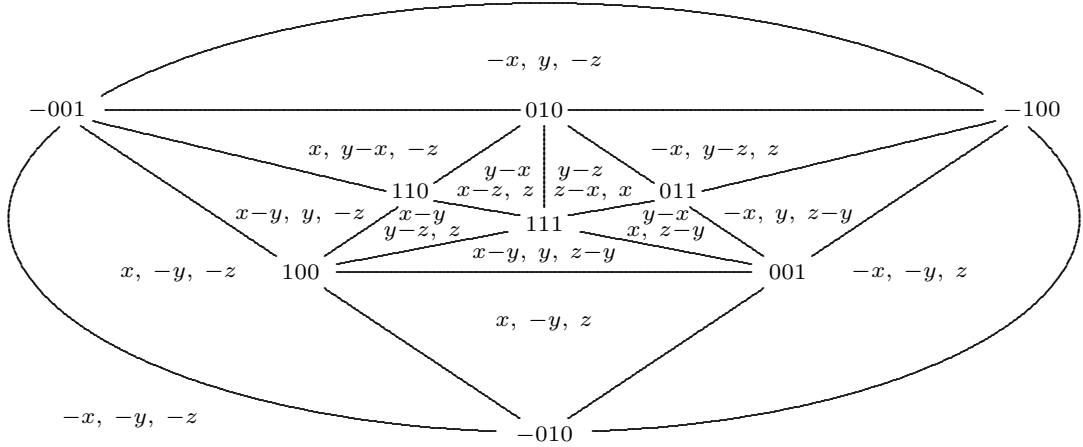
\hbox{\beginpicture
\setcoordinatesystem units <2.7cm,1.8cm>

\put{} at 0 0
\put{} at 0 2.7
\put{$\vectd 100 $} at 1 1
\put{$\vectd 110 $} at 1.45 1.5 
\put{$\vectd 010 $} at 2 2 
\put{$\vectd 111 $} at 2 1.3
\put{$\vectd 011 $} at 2.55 1.5
\put{$\vectd 001 $} at 3 1
\plot 1.1 1.1  1.4 1.4 /
\plot 1.55 1.55  1.9 1.9 /

\plot 2.1 1.9  2.45 1.55 /
\plot 2.6 1.4  2.9 1.1 /
\plot 1.15 1  2.85 1 /
\plot 2 1.9  2 1.4 /

\plot 1.55 1.44  1.9 1.35 /
\plot 2.15 1.25  2.85 1.05 /

\plot 2.45 1.44  2.1 1.35 /
\plot 1.85 1.25  1.15 1.05 /

\put{$\vectd {-1}00 $} at 4 2
\put{$\vectd {-0}10 $} at 2 0
\put{$\vectd {-0}01 $} at 0 2 
\put{$\ssize -x,\ y,\ -z$} at  2  2.3
\put{$\ssize  x,\ y-x,\ -z$} at 1.3 1.75
\put{$\ssize -x,\ y-z,\ z$} at 2.7 1.75
\put{$\ssize y-x$} at 1.85 1.62
\put{$\ssize x-z,\ z$} at 1.8 1.5

\put{$\ssize y-z$} at 2.15 1.62
\put{$\ssize z-x,\ x$} at 2.2 1.5

\put{$\ssize x-y,\ y,\ -z$} at 1 1.36
\put{$\ssize -x,\ y,\ z-y$} at 3 1.36

\put{$\ssize -x,\ -y,\ -z$} at 0.3 0.1 
\put{$\ssize x,\ -y,\ -z$} at 0.5 1 
\put{$\ssize -x,\ -y,\ z$} at 3.5 1 
\put{$\ssize x-y$} at 1.5 1.35
\put{$\ssize y-z,\ z$} at 1.5 1.25 

\put{$\ssize y-x$} at 2.5 1.35
\put{$\ssize x,\ z-y$} at 2.5 1.25 

\put{$\ssize x-y,\ y,\ z-y$} at 2 1.12 
\put{$\ssize x,\ -y,\ z$}    at 2 0.7

\plot 0.1 1.9  0.9 1.1 /
\plot 1.1 0.9  1.9 0.1 /
\plot 2.1 0.1  2.9 0.9 /
\plot 3.1 1.1  3.9 1.9 /

\plot 0.2 2  1.9 2 /
\plot 2.1 2  3.8 2 /

\plot 0.15 1.93 1.35 1.5 / 
\plot 3.85 1.93 2.65 1.5 /

\ellipticalarc axes ratio 2.5:1  110 degrees from 3.8 2.1  center at 2 1.333
\ellipticalarc axes ratio 2.5:1 -110 degrees from 1.8 0.02 center at 2 1.333
\ellipticalarc axes ratio 2.5:1  110 degrees from 2.2 0.02 center at 2 1.333

\endpicture}
$$
In general, any $(n-1)$-simplex $(M,\Cal U)$ is equipped with $n$ linear forms
$\phi_1,\dots,\phi_n$ on $K_0(A)$ such that the following holds: an $A$-module $N$ without
self-extensions belongs to  $\add M$ if and only if $\phi_i(\bdim N) \ge 0.$

In the same way as $\Sigma_A$, also $\Sigma'_A$ can be identified with a fan
in $K_0(A)\otimes \Bbb R$. For any simplex $(M,\Cal U)$ with vertices
$(E,\mod A)$ and $(0,(-S))$, where $E$ are indecomposable 
direct summands of $M$, and $S$ simple modules which do not belong to $\Cal U$,
take the cone $C(M,\Cal U)$ generated by the vectors $\bdim E$ and
$-\bdim S$.
	\bigskip
{\bf The cluster categories.} We have exhibited the cluster tilted algebras 
without reference to cluster categories, in order to show
the elementary nature of these concepts. But a genuine understanding
of cluster tilted algebras as well as of $\Sigma'_A$ 
is not possible in this way.
Starting with a hereditary artin algebra $A$, let us introduce now the corresponding
cluster category $\Cal C_A.$ We have to stress that this 
procedure reverses the historical 
development\note{In the words of Fomin and Zelevinsky [FZ4], this Part III is
completely revisionistic.}: 
the cluster categories were introduced first, 
and the cluster tilted algebras only later. The aim of the 
definition of cluster categories was to illuminate the combinatorics behind the so called
cluster algebras, in particular the combinatorics of the cluster complex. 

Let me say a little how cluster tilted algebras were found. 
Everything started with the introduction of ``cluster algebras'' 
by Fomin and Zelevinksy [FZ1]:
these are certain subrings of rational function fields, thus commutative integral domains.
At first sight, one would not guess any substantial relationship to 
non-commutative artin algebras. 
But it turned out
that the Dynkin diagrams, as well as the general Cartan data, 
play an important role for cluster algebras too. As it holds true for the 
hereditary artin algebras, it is the corresponding root system, which is of interest. 
This is a parallel situation, although not completely. For the cluster
algebras one needs to understand not only the positive roots, but the {\it almost positive} roots:
this set includes besides the
positive roots also the negative simple roots. As far as we know, 
the set of almost positive roots had not been considered
before\note{Lie theory is based on the existence of perfect symmetries --- 
partial structures (such as the set of positive roots) 
which allow only broken symmetries tend to be accepted just 
as necessary working tools. The set of almost positive roots seems to be
as odd as that of the positive ones: it depends on the same choices, but does not 
even enjoy the plus-minus merit of being half of a neat entity.
This must have been the mental reasons that the intrinsic beauty of the cluster
complex was realized only very recently. But let me stress here that the
cluster complex seems to depend not only on the choice of a root basis, but
on the ordering of the basis (or better, on the similarity class): with a
difference already for the types $\widetilde{\Bbb A}_{2,2}$ and 
$\widetilde{\Bbb A}_{3,1}$.}. 
The first link between clustor theory and tilting theory was given
by Marsh, Reineke, Zelevinksy in [MRZ] when they constructed the complex $\Sigma'_A$.
Buan, Marsh, Reineke, Reiten, Todorov [B-T]
have shown in which way the 
representation theory of hereditary artin algebras can be used in order
to construct a category $\Cal C_A$ (the cluster category) which is related to the
set of almost positive roots\note{A slight unease 
should be mentioned: as we will see,
there is an embedding of 
$\mod A$ into the cluster category which preserves indecomposability and
reflects isomorphy (but it is not a full embedding), thus this part of the cluster
category corresponds to the positive roots. There are precisely 
$n = n(A)$ additional indecomposable objects: 
they should correspond to the negative simple roots, but
actually the construction relates them to the negative of the dimension
vectors of the indecomposable projectives. Thus, 
the number of additional objects is correct, and there is even a natural bijection between
the additional indecomposable objects and the simple modules, thus the simple
roots. But in this interpretation one may hesitate to say that ``one has added
the negative simple roots'' (except in case any vertex is a sink or a source).
On the other hand, in our presentation of the cluster
complex we have used as additional vertices the elements $(0,(-S))$, and 
they really look like ``negative simple roots''. Thus, we hope
that this provides a better feeling.}
in the same way as the module category of a hereditary artin algebra
is related to the corresponding set of positive roots. 

As we have seen, a tilted algebra $B$ should be
regarded as the factor algebra of its cluster tilted algebra $\widetilde B$, 
if we want to take into account also the missing modules. 
But $\mod \widetilde B$ has to be considered as the factor category of some
triangulated category $\Cal C_A$, the corresponding cluster category. 
Looking at $\Cal C_A$, 
we obtain a common ancestor of all the algebras tilted from algebras in the
similarity class of $A$. In the setting of the pictures shown, the corresponding
cluster category has the form
$$
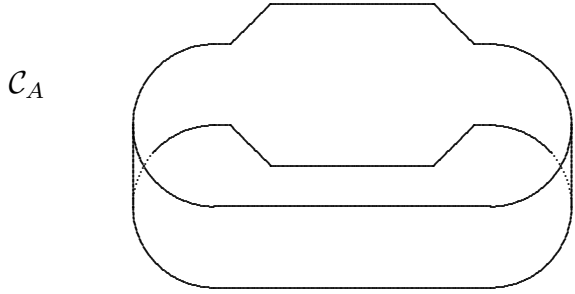
\hbox{\beginpicture
\setcoordinatesystem units <.9cm,.9cm>
\put{} at  0 2.5
\put{} at 13.2 6

\plot 4.3 4.5  4.5 4.5  5 4  7 4  7.5 4.5  7.7 4.5 / 
\plot 4.3 5.5  4.5 5.5  5 6  7 6  7.5 5.5  7.7 5.5 /

\put{$\Cal C_A$} at 2 5

\setlinear
\circulararc  90 degrees from 3.3 4.5 center at 4.3 4.5
\circulararc -90 degrees from 3.3 4.5 center at 4.3 4.5
\circulararc  90 degrees from 3.3 3.5 center at 4.3 3.5
	\setdots <.5mm>
\circulararc -40 degrees from 3.3 3.5 center at 4.3 3.5
	\setsolid
\circulararc  50 degrees from 4.3 4.5 center at 4.3 3.5

\circulararc -90 degrees from 8.7 4.5 center at 7.7 4.5
\circulararc  90 degrees from 8.7 4.5 center at 7.7 4.5
\circulararc -90 degrees from 8.7 3.5 center at 7.7 3.5
	\setdots <.5mm>
\circulararc  40 degrees from 8.7 3.5 center at 7.7 3.5
	\setsolid
\circulararc -50 degrees from 7.7 4.5 center at 7.7 3.5

	\setsolid
\setlinear
\plot 3.3 4.5  3.3 3.5 /
\plot 8.7 4.5  8.7 3.5 /

\plot 4.3 2.5  7.7 2.5 /
\plot 4.3 3.5  7.7 3.5 /

\endpicture}
$$
The cluster category $\Cal C_A$ should be considered as a universal kind of
category belonging to the similarity class of the hereditary artin algebra $A$ 
in order to obtain all the module categories $\mod \widetilde B,$ where $\widetilde B$
is a cluster tilted algebra of type similar to $A$.

What one does is the following:
start with the derived category $D^b(\mod A)$ of the hereditary artin algebra $A$, with
shift functor $[1]$, 
and take as $\Cal C_A$ the orbit
category with respect to the functor $\tau_D^{-1}[1]$ (we write $\tau_D$ for the
Auslander-Reiten translation in the derived category, and $\tau_c$ for the 
Auslander-Reiten translation in $\Cal C_A$).
As a fundamental domain for the
action of this functor one can take the disjoint union of $\mod A$ (this yields all the
positive roots) and the shifts of the projective $A$-modules by [1] (this yields
$n = n(A)$ additional indecomposable objects).
It should be mentioned that Keller [K] has shown that
$\Cal C_A$ is a triangulated category. Now if we take a tilting module $T$ in $\mod A$,
we may look at the endomorphism ring $\widetilde B$ of $T$ in $\Cal C_A$ 
(or better: the endomorphism ring of the image of $T$ under the canonical functors
$\mod A \subseteq D^b(\mod A) \to \Cal C_A$), and obtain a cluster
tilted algebra\note{This is the way, the cluster tilted algebras were introduced and studied by
Buan, Marsh and Reiten [BMR1].} as considered above.
The definition immediately yields that $\widetilde B = B\semi J$, 
where $J = \Hom_{D^b(\mod A)}(T,\tau_D^{-1}T[1]) = \Ext^1_A(T,\tau^{-1}T).$ The decisive
property is that there is a canonical equivalence of 
categories\note
{Instead of $\Cal C_A/\langle T\rangle$, one may also take the equivalent category 
$\Cal C_A/\langle \tau_c T\rangle$. The latter is of interest 
if one wants  the indecomposable summands
of $T$ in $\Cal C_A$ to become indecomposable projective objects.}
$$
 \Cal C_A/\langle T\rangle  \longrightarrow \mod \widetilde B.
$$ 
 In particular, we see
that the triangulated category $\Cal C_A$ has many factor categories which are 
abelian\note
{We have mentioned that the cluster theory brought many surprises. Here is another one:
One knows for a long time 
many examples of abelian categories $\Cal A$ with an object $M$ such that
the category $\Cal A/\langle M\rangle$ (obtained by setting zero all maps which factor
through $\add M$) becomes a triangulated category: just take $\Cal A = \mod R$, where $R$
is a self-injective artin algebra $R$ and $M = {}_RR$. The category
$\mod R/\langle{}_RR\rangle = \underline{\mod} R$ is the stable module category of $R$.
But we are not aware that non-trivial 
examples where known of a triangulated category $\Cal D$ with an object $N$ such that
$\Cal D/\langle N\rangle$ becomes abelian. Cluster tilting theory is just about this!}.
	\medskip
What happens when we form the factor category $\Cal C_A/\langle T\rangle$? 
Consider an indecomposable direct summand $E$ of the tilting $A$-module $T$
as an object in the cluster category $\Cal C_A$ and the meshes 
starting and ending in $E$:
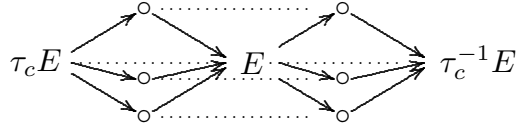
$$
\hbox{\beginpicture
\setcoordinatesystem units <1cm,.6cm>
\put{} at 0 0
\put{} at 0 0
\put{$\tau_c E$} at -.2 0
\put{$E$} at 2.2 0
\put{$\tau_c^{-1}E$} at 4.7 0
\multiput{$\circ$} at 1 -1  1 -0.3 1 1  3.2 -1 3.2 -0.3 3.2 1 /
\arr{0.2 0.2}{0.9 0.9}
\arr{0.2 0}{0.9 -0.3}
\arr{0.2 -0.2}{0.9 -0.9}
\arr{1.1 0.9}{1.9 0.1}
\arr{1.1 -0.3}{1.9 0}
\arr{1.1 -0.9}{1.9 -0.1}

\arr{2.5 0.3}{3.1 0.9}
\arr{2.5 0}{3.1 -0.3}
\arr{2.5 -0.3}{3.1 -0.9}
\arr{3.3 0.9}{4.1 0.1}
\arr{3.3 -0.3}{4.1 0}
\arr{3.3 -0.9}{4.1 -0.1}
\setdots <1mm>
\plot 0.2 0  1.9 0 /
\plot 2.4 0  3.8 0 /
\plot 1.2 1  2.8 1 /
\plot 1.2 -0.3  2.8 -0.3 /
\plot 1.2 -1  2.8 -1 /
\endpicture}
$$
In the category $\Cal C_A/\langle T\rangle$, the object $E$ becomes zero, whereas both
$\tau_c E$ and $\tau_c^{-1}E$ remain non-zero. In fact, $\tau_c^{-1}E$ becomes a
projective object and $\tau_c E$ becomes an injective object: 
We obtain in this way in $\mod \widetilde B =
\Cal C_A/\langle T\rangle$ an indecomposable projective module $P = \tau_c^{-1}E$ and an
indecomposable injective module $I = \tau_c E$, such that $\top P \simeq \soc I.$ 
This explains the round trip phenomenon for $\widetilde B$ mentioned above:  
there is the hammock corresponding to the simple $\widetilde B$-module 
$\top P \simeq \soc I$, starting from $I = \tau_c E$,
and ending in $P = \tau_c^{-1}E$. And either $\rad P$ is projective
(and $I/\soc I$ injective) or else there are non-sectional paths of length 4 from
$I$ to $P$.
	\medskip
There is a decisive symmetry 
condition\note{If we write $\Ext^1(X,Y) = \Hom_{\Cal C}(X,Y[1])$, then this symmetry condition reads
that $\Ext^1(X,Y)$ and $\Ext^1(Y,X)$ are dual to each other, in particular they have
the same dimension.} in the cluster category $\Cal C = \Cal C_A$: 
$$
 \Hom_{\Cal C}(X,Y[1]) \simeq D\Hom_{\Cal C}(Y,X[1]).
$$
This is easy to see: since we form the orbit category with respect to $\tau_D^{-1}[1]$,
this functor becomes the identity functor in $\Cal C$, and therefore the 
Auslander-Reiten functor $\tau_c$ and the shift functor $[1]$ in $\Cal C$
coincide. On the other hand, the Auslander-Reiten (or Serre duality) formula for $\Cal C$
asserts that $\Hom_{\Cal C}(X,Y[1]) \simeq D\Hom_{\Cal C}(Y,\tau_c X).$
A triangulated category is said to be {\it $d$-Calabi-Yau} provided the shift functor $[d]$
is a Serre (or Nakayama) functor, thus provided there is a functorial isomorphism
$$
 \Hom(X,-) \simeq D\Hom(-,X[d])
$$
(for a discussion of this property, see for example [K]). 
As we see, {\it the cluster category is $2$-Calabi-Yau.}
	\medskip
The cluster category has Auslander-Reiten sequences. One component
$\Gamma_0$ of the Auslander-Reiten quiver of $\Cal C_A$ has only
finitely many $\tau_{\Cal C}$-orbits, namely the component
containing the indecomposable projective (as well as the indecomposable
injective) $A$-modules. The remaining components of the Auslander-Reiten
quiver of $\Cal C_A$ have tree class $\Bbb A_\infty.$
	\medskip
In a cluster category $\Cal C = \Cal C_A$, an object is said to be a 
{\it cluster-tilting 
object}\note{It has to be stressed that the notion of a ``cluster-tilting object'' in a cluster 
category does not conform to the tilting notions used otherwise in this Handbook!
If $T$ is such a cluster-tilting 
object, then it may be that $\Hom_{\Cal C}(T,T[i]) \neq 0$ in 
$\Cal C = \Cal C_A$ for some $i \ge 2.$}
 provided first
$\Hom_{\Cal C}(T,T[1]) = 0$, and second, that $T$ is maximal with this property in the following sense:
if $\Hom_{\Cal C}(T\oplus X,(T\oplus X)[1]) = 0$, then $X$ is in $\add T$. If $T$ is a tilting $A$-module,
then one can show that $T$, considered as an object of $\Cal C_A$ is a cluster-tilting object.

Let us consider the hereditary artin algebra in one
similarity class and the reflection functors between them. One may 
identify the corresponding cluster categories using the reflection functors,
as was pointed out by Bin Zhu [Zh2].
In this way, one can compare the tilting modules of all the
hereditary artin algebras in one similarity class. And it turns out 
that {\it the cluster-tilting objects in $\Cal C_A$ are just the tilting modules for
the various artin algebras obtained from $A$ by using reflection functors} [B-T]. 
In order to see this, let $T$ be a cluster-tilting
object in $\Cal C_A$. Let $\Gamma_0$ be the component of the
Auslander-Reiten quiver of $\Cal C_A$ which contains the indecomposable
projective $A$-modules. If no indecomposable direct summand of $T$ belongs
to $\Gamma_0,$ then $T$ can be considered as an $A$-module, and it is
a regular tilting $A$-module. On the other hand, if there is an indecomposable
direct summand of $T$, say $T_1$, which belongs to $\Gamma_0$, then let
$\Cal S$ be the class of all indecomposable objects $X$ in $\Gamma_0$ with
a path from $X$ to $T_1$ in $\Gamma_0$, and such that any path from $X$ to $T_1$
in $\Gamma_0$ is sectional. Then no indecomposable direct summand of $T$
belongs to $\tau_{\Cal C}\Cal S$. We may identify the factor category
$\Cal C_A/\langle\tau_{Cal C}\Cal S\rangle$ with $\mod A'$ for some
hereditary artin algebra $A'$, and consider $T$ as an $A'$-module (the object
$T_1$, considered as an $A'$-module, is projective and faithful). Clearly,
$A'$ is obtained from $A$ by a sequence of BGP-reflection functors.

Also, the usual procedure of going from a tilting
module to another one by exchanging just one indecomposable direct summand gets more regular.
Of course, there is the notion of an almost complete partial 
cluster-tilting object and of a complement,
parallel to the corresponding notions of an almost complete partial tilting module and its
complements. Here we get:
{\it Any almost complete partial cluster-tilting object $\overline T$ has precisely two complements} 
[B-T]. 
We indicate the proof: 
We can assume that $\overline T$ is an $A$-module. If $\overline T$ is sincere, then we
know that there are two complements for $\overline T$ considered as an almost complete
partial tilting $A$-module. If $\overline T$ is not sincere, then there is only one
complement for $\overline T$ considered as an almost complete partial tilting $A$-module.
But there is also one (and 
obviously only one) indecomposable projective module $P$ with $\Hom_A(P,\overline T) = 0,$
and the $\tau_c$-shift of $P$ in the cluster category is the second complement we are
looking for!
	\medskip
The main 
point seems to be the following: {\it The simplicial 
complex of partial cluster-tilting
objects in the cluster category $\Cal C_A$ is nothing else than $\Sigma'_A$}, with
the following identification: If $T$  is a 
basic partial cluster-tilting object in $\Cal C_A$, we can write $T$ as the direct sum of a module
$M$ in $\mod A$ and objects of the form $\tau_cP(i)$, with $P(i)$ 
indecomposable projective in $\mod A$, and $i$ in some index set $\Theta$. 
Then $M$ corresponds in $\Sigma'_A$ to the pair
$(M,\Cal U)$, where $\Cal U = \bigcap_{i\in \Theta} (-S(i)).$ The reason is very simple:
$\Hom_{\Cal C}(\tau_cP(i),M[1]) \simeq \Hom_{\Cal C}(P(i),M) = \Hom_A(P(i),M),$
with $\Cal C = \Cal C_A.$ 
	
The complex $\Sigma'_A$ should be viewed as a convenient index 
scheme\note{But we should also mention the following:  
The set of isomorphism classes of basic cluster-tilting objects in $\Cal C_A$ 
is no longer partially ordered.
In fact, given an almost complete partial cluster-tilting object $\overline T$ and its two
complements $X$ and $Y$, there are triangles $X \to T' \to Y \to$ and
$Y\to T'' \to X \to$ with $T',T''\in \add \overline T$.} 
for
the set of cluster tilted algebras obtained from the hereditary artin algebras  
in the similarity class of $A$. Any maximal simplex of $\sigma'_A$ 
is a cluster-tilting object in $\Cal C_A$,
and thus we can attach to it its endomorphism ring. Let us redraw the complex $\Sigma'_A$ 
for the path algebra $A$ of the quiver $\circ \leftarrow \circ\leftarrow \circ,$ 
so that the different vertices and triangles are better seen: 
$$
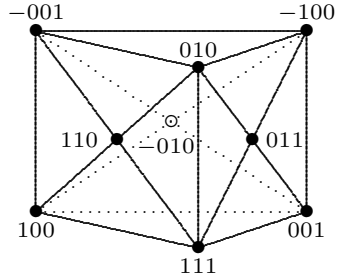
\hbox{\beginpicture
\setcoordinatesystem units <.6cm,.4cm>
\put{} at 0 -1
\put{} at 5 5
\multiput{$\bullet$} at 0 0  0 5  5 0  5 5  3 -1  3 4  1.5 2  4 2  /
\put{$\circ$} at 2.5 2.5 
\plot 0 0  3 -1  5 0  5 5  0 5  0 0 /
\plot 0 0  3 4  5 0 /
\plot 0 5  3 -1  5 5 /
\plot 0 5  3 4  5 5 /
\plot  3 4  3 -1 /
\setdots <1mm>
\plot 0 0  5 5  /
\plot 0 5  5 0  0 0 /
\put{$\ssize -001$} at 0 5.5
\put{$\ssize -010$} at 2.4 1.8
\put{$\ssize -100$} at 5 5.5
\put{$\ssize 100$} at 0 -.5
\put{$\ssize 001$} at 5 -.5
\put{$\ssize 111$} at 3 -1.5
\put{$\ssize 010$} at 3 4.5
\put{$\ssize 110$} at .8 2
\put{$\ssize 011$} at 4.6 2
\endpicture}
$$
There are two
kinds of vertices, having either 4 or 5 neighbours. The vertices with 5 neighbours form two
triangles (the bottom and the top triangle), 
and these are the cluster-tilting objects with endomorphism ring of infinite global dimension.
The remaining triangles yield hereditary endomorphism rings and again, there
are two kinds: The quiver of the endomorphism ring may have one sink and one source, these
rings are given by the six triangles which have an edge in common with the bottom or the top
triangle. Else, the endomorphism ring is hereditary and the radical square is zero:
these rings correspond to the remaining six triangles:
$$
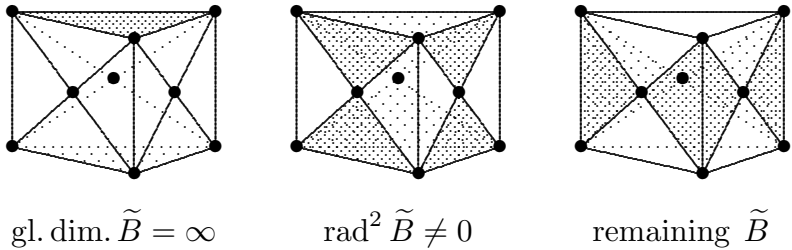
\hbox{\beginpicture
\setcoordinatesystem units <.45cm,.3cm>
\put{\beginpicture
\setcoordinatesystem units <.45cm,.3cm>
\put{} at 0 -1
\put{} at 5 5
\multiput{$\bullet$} at 0 0  0 5  5 0  5 5  3 -1  3 4  1.5 2  4 2  2.5 2.5 /
\plot 0 0  3 -1  5 0  5 5  0 5  0 0 /
\plot 0 0  3 4  5 0 /
\plot 0 5  3 -1  5 5 /
\plot 0 5  3 4  5 5 /
\plot  3 4  3 -1 /
\setdots <1mm>
\plot 0 0  5 5  /
\plot 0 5  5 0  0 0 /
\setshadegrid span <.4mm>
\hshade 4 3 3  5 0 5 /

\setshadegrid span <.8mm>
\hshade -1  3 3  0 0 5  /

\endpicture} at 0 0
\put{\beginpicture
\setcoordinatesystem units <.45cm,.3cm>
\put{} at 0 -1
\put{} at 5 5
\multiput{$\bullet$} at 0 0  0 5  5 0  5 5  3 -1  3 4  1.5 2  4 2  2.5 2.5 /
\plot 0 0  3 -1  5 0  5 5  0 5  0 0 /
\plot 0 0  3 4  5 0 /
\plot 0 5  3 -1  5 5 /
\plot 0 5  3 4  5 5 /
\plot  3 4  3 -1 /
\setdots <1mm>
\plot 0 0  5 5  /
\plot 0 5  5 0  0 0 /
\setshadegrid span <.45mm>
\vshade 0 5 5 <,z,,> 1.5 2 4.5 <z,,,> 3 4 4  /
\vshade 3 4 4 <,z,,> 4   2 4.5 <z,,,> 5 5 5  /

\vshade 0 0 0 <,z,,> 1.5 -.5  2  <z,,,> 3 -1 -1 /
\vshade 3 -1 -1 <,z,,> 4 -.5  2  <z,,,> 5 0 0  /

\setshadegrid span <.8mm>
\vshade 0 5 5  <,z,,> 2.5 2.5 5 <z,,,>  5 5 5  /
\vshade 0 0 0  <,z,,> 2.5 0 2.5 <z,,,>  5 0 0   /

\endpicture} at 7 0
\put{\beginpicture
\setcoordinatesystem units <.45cm,.3cm>
\put{} at 0 -1
\put{} at 5 5
\multiput{$\bullet$} at 0 0  0 5  5 0  5 5  3 -1  3 4  1.5 2  4 2  2.5 2.5 /
\plot 0 0  3 -1  5 0  5 5  0 5  0 0 /
\plot 0 0  3 4  5 0 /
\plot 0 5  3 -1  5 5 /
\plot 0 5  3 4  5 5 /
\plot  3 4  3 -1 /
\setdots <1mm>
\plot 0 0  5 5  /
\plot 0 5  5 0  0 0 /
\setshadegrid span <.45mm>
\vshade 0 0 5    1.5 2 2 /
\vshade 1.5 2 2  3 -1 4 /
\vshade 3 -1 4  4 2 2 /
\vshade 4 2 2  5 0 5 /
\setshadegrid span <.8mm>
\vshade 0 0 5  2.5 2.5 2.5  5 0 5  /

\endpicture} at 14 0
\put{$\gldim \widetilde B = \infty$} at 0 -5
\put{$\rad^2 \widetilde B \neq 0$} at  7 -5

\put{remaining $\widetilde B$} at  14 -5
\endpicture}
$$

Consider an almost complete partial cluster-tilting object $\overline T$ in
$\Cal C = \Cal C_A$.
As we have mentioned, there are precisely two complements for $\overline T$, say $E$ and $E'$.
Let $T = \overline T\oplus E,$ and 
$T' = \overline T \oplus E' .$ Thus, there are given two cluster-tilted algebras 
$\widetilde B = \End_{\Cal C}(T)$, and $\widetilde {B'} = \End_{\Cal C}(T')$, 
we may call them {\it adjacent,} this corresponds to the position of $T$ and $T'$ in the
complex $\Sigma'_A$. 
We can identify $\Cal C_A/\langle T \rangle$ with $\mod \widetilde B,$ 
and $\Cal C_A/\langle T' \rangle$ with $\mod \widetilde {B'}$ 
We saw that $E$ as an indecomposable direct summand of $T$ yields an indecomposable projective
$\widetilde B$-module $P = \tau^{-1}_cE$ and an indecomposable injective 
$\widetilde B$-module $I = \tau_cE$, such that $\soc I \simeq \top P.$ 
Since $\overline T \oplus E'$ is a cluster-tilting object, it is not difficult to show that
$\Hom_{\widetilde B}(P',E') = 0$ for any indecomposable projective $\widetilde B$-module $P'
\not\simeq P$. But this implies that $E'$ is identified under the equivalence of
$\Cal C_A/\langle T \rangle$ and $\mod \widetilde B$ with 
the simple $\widetilde B$-module which is
the socle of $I$ and the top of $P$. In the same way, 
we see that $E$ is a simple $\widetilde{B'}$-module, namely the top of the
$\widetilde {B'}$-module $P = \tau^{-1}_cE'$ and the socle of the indecomposable injective 
$\widetilde {B'}$-module $I = \tau_cE'$.
Thus, there is the following sequence of
identifications:
$$
 \mod \widetilde B/\langle \add E'\rangle \simeq \Cal C_A/ \langle \add(T\oplus E')\rangle =
 \Cal C_A/ \langle \add(T'\oplus E)\rangle \simeq
\mod \widetilde {B'}/\langle \add E\rangle.
$$
Altogether this means that artin algebras $\widetilde B$ and $\widetilde{B'}$ which are
adjacent, are
nearly Morita equivalent [BMR2]. We had promised to the reader, that we will return
to the hammock configuration $(P,S,I)$, 
where $S$ is a simple $\widetilde B$-module, $P = P(S)$ its
projective cover, and $I = I(S)$ its injective envelope: but this is the present setting. 
Using the cluster category
notation, we can write $P = \tau^{-1}_cE,$  $I = \tau_cE$, and then $S = E'$, where 
$E, E'$ are
complements to an almost complete partial cluster-tilting object $\overline T$. 
When we form the category $\mod \widetilde B/\langle \add S\rangle$, the killing of the
simple $\widetilde B$-module $S$ creates a hole in $\mod \widetilde B$. From the
hammock $\Hom(P,-)$ in $\mod \widetilde B$ the following parts survive:
$$
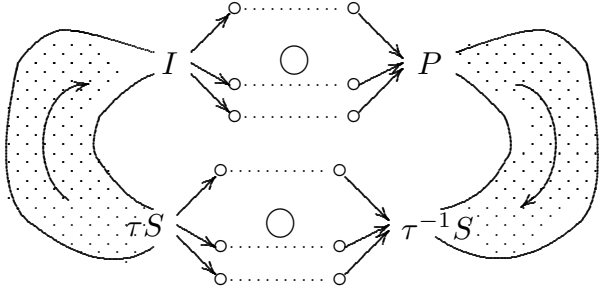
\hbox{\beginpicture
\setcoordinatesystem units <.6cm,.6cm>

\put{\beginpicture
\put{$\bigcirc$} at 2.1 0.05

\put{$\tau S$} at -0.4 0
\put{$\tau^{-1}S$} at 5 0
\multiput{$\circ$} at 1 -1  1 -0.4 1 1  3.2 -1 3.2 -0.4 3.2 1 /
\arr{0.2 0.2}{0.9 0.9}
\arr{0.2 0}{0.9 -0.4}
\arr{0.2 -0.2}{0.9 -0.9}

\arr{3.3 0.9}{4.1 0.1}
\arr{3.3 -0.4}{4.1 0}
\arr{3.3 -0.9}{4.1 -0.1}
\setdots <1mm>
\plot 1.2 1     3 1 /
\plot 1.2 -0.4  3 -0.4 /
\plot 1.2 -1    3 -1 /

\endpicture} at 2.2 -3

\put{$\bigcirc$} at 2.1 0.05 

\put{$I$} at -0.2 0
\put{$P$} at 4.6 0
\multiput{$\circ$} at 1 -1  1 -0.4 1 1  3.2 -1 3.2 -0.4 3.2 1 /
\arr{0.2 0.2}{0.9 0.9}
\arr{0.2 0}{0.9 -0.4}
\arr{0.2 -0.2}{0.9 -0.9}

\arr{3.3 0.9}{4.1 0.1}
\arr{3.3 -0.4}{4.1 0}
\arr{3.3 -0.9}{4.1 -0.1}
\setdots <1mm>
\plot 1.2 1     3 1 /
\plot 1.2 -0.4  3 -0.4 /
\plot 1.2 -1    3 -1 /

\setshadegrid span <1mm>
\hshade -3.8  -.4 -.6 <,,,z> -3.5 -2 -.7 <,,z,> -2.5 -3.3 -0.4 / 
\hshade -2.6 -3.3 -1 <,,,z> -1 -3.3 -1.7 <,,z,z> 0 -3 -1 <,,z,> 0.5 -2.5 -1.7 / 
\hshade -2.6  5.5 7.8 <,,,z>  -1 6.5 7.8 <,,z,z> 0 5 7.2 <,,z,> 0.5 5 7 / 

\hshade -3.6  5.5 5.5 <,,,z> -3.5 5.5 6.5 <,,z,> -2.5 5 7.5 / 
\hshade -2.6 -3.3 -1 <,,,z> -1 -3.3 -1.7 <,,z,z> 0 -3 -1 <,,z,> 0.5 -2.5 -1.7 / 
\hshade -2.6  5.5 7.8 <,,,z>  -1 6.5 7.8 <,,z,z> 0 5 7.2 <,,z,> 0.5 5 7 / 

\setsolid

\setquadratic
\plot  5.1 0.2   
       6 .5     
       6.5  .6   
       7. .5    
       7.5 0    
       7.7  -2  
       7.3 -3  
       6  -3.6  
       4.7 -3.2 /
\plot   -.5 0.2   
        -1.5 .5    
        -2 .6 
        -2.5  .5  
        -3 0  
        -3.2 -2
        -2.7 -3  
        -1.5 -3.6 
         -.5 -3.4 /

\plot  5.1 -0.2  
       5.5 -.5  
        6 -1  
        6.1 -1.5 
        6 -1.8   
       5.5 -2.4  
       4.8 -2.7   /
\plot -.5 -0.2 
        -1 -.5  
        -1.5 -1  
        -1.6 -1.5  
        -1.5 -1.8
        -1 -2.4
         -.5 -2.7   /

\ellipticalarc axes ratio 1:1.5 -160 degrees from 6.2 -.4 center at 6.2 -1.5  
\ellipticalarc axes ratio 1:1.5  160 degrees from -1.8 -.4 center at -1.8 -1.5  
\arr{-1.8 -.4}{-1.7 -.38}
\arr{6.5 -2.5}{6.3 -2.65}
\put{} at 9 0 

\endpicture}
$$
Note that the 
new hole is of the same nature as the hole between $I$ and $P$ (which was created when
we started with the cluster category $\Cal C$ by killing the object $E$).
Indeed,
one may fill alternatively one of the two holes and obtains $\mod \widetilde B$, or
$\mod \widetilde{B'}$, respectively. 

Altogether, we see: 
A cluster category $\Cal C = \Cal C_A$ has a lot of nice factor categories 
which are abelian (the module categories $\mod \widetilde B$), and one should regard 
$\Cal C$ as being obtained from patching together the various factor categories 
in the same way as manifolds are built up from open subsets by specifying the identification
maps of two such subsets along what will become their intersection. The patching process for the
categories $\mod \widetilde B$ is done via the nearly Morita equivalences
for adjacent tilting 
objects\note{It seems that there is not yet any kind of axiomatic approach to
this kind of patching process.
}. 

The reader will have noticed that this exchange process for adjacent algebras
generalizes the BGP-reflection functors (and the APR-tilting functors) to vertices
which are not sinks or sources. Indeed, for 
$\widetilde B = \End_{\Cal C}(\overline T\oplus E)$, and $\widetilde {B'} = \End_{\Cal C}(\overline T\oplus E')$, the indecomposable direct summand $E$ of 
$\overline T\oplus E$ corresponds to a vertex of the quiver of $\widetilde B$,
and similarly, $E'$  corresponds to a vertex of the quiver of $\widetilde {B'}$.
In the BGP and the APR setting, one of the modules $E,E'$ is simple projective, the
other one is simple injective --- here now $E$ and $E'$ are arbitrary exceptional
modules\note{A direct description of this reflection process seems to be still 
missing. It will require a proper understanding of all the cluster tilted algebras
$\widetilde B$ with $n(\widetilde B) = 3.$ A lot is already known about such algebras, 
see [BMR2].}. 
	\medskip
This concludes our attempt to report about some of the new results in tilting theory 
which are based on cluster categories. Let us summerize the importance of this
development. First of all, the cluster tilted algebras provide a nice
depository for storing the modules which are lost when we pass from hereditary artin
algebras to tilted algebras; there is a magic bimodule which controls the situation.
We obtain in this way a wealth of algebras whose module categories are described
by the root system of a Kac-Moody Lie-algebra. These new algebras are no longer
hereditary, but are still of Gorenstein dimension at most 1. For the class of
cluster tilted algebras, there is a reflection process at any vertex of the quiver,
not only at sinks and sources. This is a powerful generalization of the
APR-tilting functors (thus also of the BGP-reflection functors), and 
adjacent cluster tilted algebras are nearly Morita equivalent. The index set for this
reflection process is the simplicial complex $\Sigma'_A$ and the introduction
of this simplicial complex solved also another riddle of tilting theory: it provides
a neat way of enlarging the simplicial complex of tilting $A$-modules in order
to get rid of its boundary. We have mentioned in Part I that both the missing
modules problem as well as the boundary problem concern the module category,
but disappear on the level of derived categories. Thus it is not too
surprising that derived categories play a role: as it has turned out,
the cluster categories, as suitable orbit categories of the 
corresponding derived categories, are the decisive new objects. These are
again triangulated categories, and are to be considered as the universal
structure behind all the tilted and cluster tilted algebras obtained from 
a single hereditary artin algebra. 
	\bigskip
\centerline{\bf Appendix: Cluster Algebras.}
	\medskip
Finally we should speak about the source of all these developments,
the introduction of cluster algebras by Fomin and Zelevinsky. But we are hesitant,
for two reasons: first, there is our complete lack of proper expertise, but also
it means that we leave the playground of tilting theory. Thus this will be just
an appendix to the appendix. 
The relationship between cluster algebras on the one hand,
and the representation theory of hereditary artin algebras and cluster tilted algebras
on the other hand is fascinating, but also very 
subtle\note{Since this report is written for the Handbook of Tilting Theory, 
we are only concerned with the relationship of the cluster algebras to tilting theory. 
There is 
a second relationship to the representation theory of artin algebras,
namely to Hall algebras, as found by Caldero and Chapoton [CC], and Caldero-Keller
[CK1,CK2] , see also Hubery [Hu]. 
And there are numerous interactions between cluster theory and many different parts of mathematics. But all this lies beyond the scope of this volume.
}. 
At first, one
observed certain analogies and coincidences. Then there was an experimental
period, with many surprising findings (for example, that the Happel-Vossieck list of
tame concealed algebras corresponds perfectly to the Seven list of minimal infinite cluster
algebras [S], as explained in [BRS]). In the meantime,
many applications of cluster-tilted algebras to cluster algebras have been found [BR, BMRT],
and the use of Hall algebra methods provides a conceptual
understanding of this relationship [CC, CK1, CK2, Hu].

Here is at least a short indication what cluster algebras are.
As we said already, the cluster algebras are (commutative) integral domains. The 
cluster algebras we
are interested in (those related to hereditary artin 
algebras)\note{these are the so-called acylic cluster algebras [FM3].}
are finitely generated (this means finitely generated ``over nothing'', say over $\Bbb Z$), 
thus they can be considered as subrings of a finitely generated function field 
$\Bbb Q(x_1,\dots,x_n)$ over the rational numbers $\Bbb Q$. 
This is the way they usually are presented in the literature (but the finite generation  
is often not stressed). In fact, one of the main theorems of cluster theory asserts
that we deal with subrings of the ring of Laurent polynomials $\Bbb Z[x_1^{\pm 1},\dots,
x_n^{\pm 1}]$ (this is the subring of all elements of the form $\frac pq$ where $p$
is in  the polynomial ring $\Bbb Z[x_1,\dots,x_n]$ and $q$ is a monomial in the variables
$x_1,\dots,x_n$).

Since we deal with a noetherian integral domain,
the reader may expect to be confronted with problems in algebraic geometry, or, since
we work over $\Bbb Z$ with those of arithmetical geometry. But this was not the primary
interest. Instead, the cluster theory belongs in some sense to algebraic combinatorics,
and the starting question concerns the existence of a nice $\Bbb Z$-basis of such a cluster
algebra, say similar to all the assertions about canonical bases in Lie theory. 

What are clusters? Recall that a cluster algebra is a subring of 
$\Bbb Z[x_1^{\pm 1},\dots,x_n^{\pm 1}]$. A cluster algebra always 
has a $\Bbb Z$-basis 
consisting of elements of the form $\frac pq$, where $p \in \Bbb Z[x_1,\dots,x_n]$ 
is not divisible by the variables $x_1,\dots,x_n$ and
$q = x_1^{d_1}\cdots x_n^{d_n}$ with exponents $d_i \in \Bbb Z.$ There is
an inductive procedure to produce such a basis and the elements $\frac pq$ obtained in
this way are called the cluster variables. The Laurent monomial $q$ is said to be the 
{\it denominator}\note{The variables $x_1,\dots,x_n$ are cluster variables.
Such an element $x_i$ is considered as being written in the form $1/(x_i^{-1})$, 
its denominator is $q = x_i^{-1}.$} of $\frac pq$ and 
$\bdim q = (d_1,\dots,d_n)$ its dimension vector. 
One of the main topics discussed recently 
in cluster theory concerns the denominators of the cluster variables.

Consider now the case of the path algebra of a finite quiver
$Q$ without oriented cycles. According to Caldero-Keller [CK2], the simplicial
complex $\Sigma'_A$ with $A = kQ$ 
can be identified with the cluster complex corresponding to $Q$.
Under this correspondence, the cluster variables correspond to the 
exceptional $A$-modules and the elements of the form $(0,(-S))$.
When we introduced the simplicial complex $\Sigma'_A$, the
maximal simplices were labeled $(M,\Cal U)$ with $M$ a basic tilting module in a
Serre subcategory $\Cal U$ of $\mod A$. 
Recall that such an $(n-1)$-simplex $(M,\Cal U)$ in $\Sigma'_A$ 
is equipped with $n$ linear forms
$p_1,\dots,p_n$ on $K_0(A)$ such that an $A$-module $N$ without
self-extensions belongs to $\add M$ if and only if $\phi_i(\bdim N) \ge 0,$
for $1\le i \le n$. And there is the parallel assertion: 
A cluster monomial with denominator $q$ belongs to the
cluster corresponding to $(M,\Cal U)$ 
if and only if $\phi_i(\bdim q) \ge 0,$ for
$1 \le i \le n$. 

Here are the cluster variables for the cluster algebra of type $\Bbb A_3$,
inserted as the vertices of the cluster complex $\Sigma'_A$:
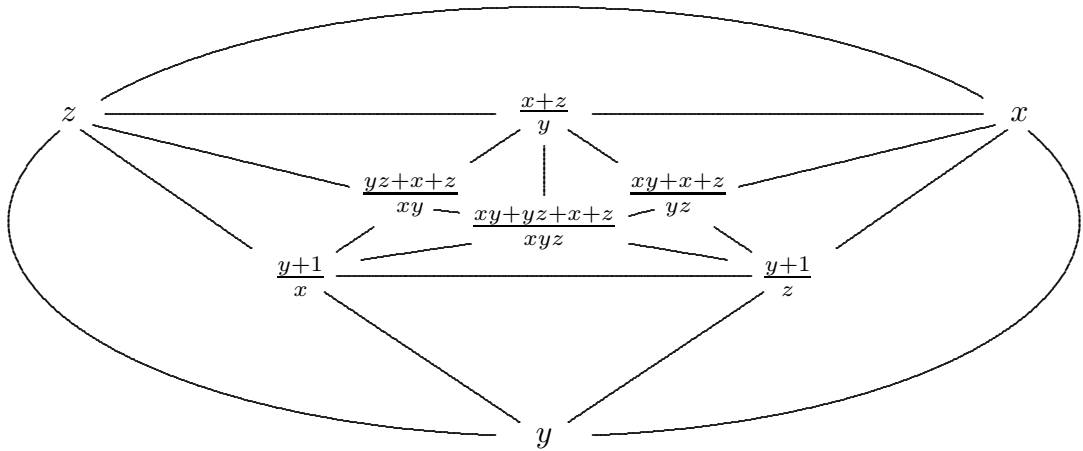
$$
\hbox{\beginpicture
\setcoordinatesystem units <2.7cm,1.8cm>

\put{} at 0 0
\put{} at 0 2.7
\put{$\frac{y+1}x$} at 1 1
\put{$\frac{yz+x+z}{xy}$} at 1.45 1.5 
\put{$\frac{x+z}y$} at 2 2 
\put{$\frac{xy+yz+x+z}{xyz}$} at 2 1.3
\put{$\frac{xy+x+z}{yz}$} at 2.55 1.5
\put{$\frac{y+1}z $} at 3 1
\plot 1.15 1.15  1.3 1.3 /
\plot 1.7 1.7  1.9 1.9 /

\plot 2.1 1.9  2.3 1.7 /
\plot 2.7 1.3  2.85 1.15 /

\plot 1.15 1  2.85 1 /
\plot 2 1.8  2 1.5 /

\plot 1.55 1.41  1.65 1.39 /
\plot 2.35 1.2  2.75 1.1 /

\plot 2.45 1.41  2.35 1.37 /
\plot 1.7 1.2  1.25 1.1 /

\put{$x $} at 3.95 2
\put{$y $} at 2 0
\put{$z $} at 0.05 2 

\plot 0.1 1.9  0.8 1.17 / 
\plot 1.1 0.9  1.9 0.1 /
\plot 2.1 0.1  2.9 0.9 / 
\plot 3.2 1.17  3.9 1.9 / 

\plot 0.2 2  1.8 2 /
\plot 2.2 2  3.8 2 /

\plot 0.15 1.93 1.2 1.55 / 
\plot 3.85 1.93 2.8 1.55 /

\ellipticalarc axes ratio 2.5:1  110 degrees from 3.8 2.1  center at 2 1.333
\ellipticalarc axes ratio 2.5:1 -110 degrees from 1.8 0.02 center at 2 1.333
\ellipticalarc axes ratio 2.5:1  110 degrees from 2.2 0.02 center at 2 1.333

\endpicture}
$$
	\bigskip
{\bf Acknowledgment.} 
We were able to provide only a
small glimpse of what cluster tilting theory is about --- it is a
theory in fast progress, but also with many problems still open.
The report would not have been possible without the constant
support of many experts. In particular, we have to mention Idun Reiten who
went through numerous raw versions. The author has also
found a lot of help by Thomas Br\"ustle, Philippe Caldero,
Dieter Happel, Otto Kerner and Robert Marsh, who answered all
kinds of questions. And he is indebted to Aslak Buan, Philipp
Fahr, Rolf Farnsteiner, Lutz Hille, Bernhard Keller and Markus
Reineke for helpful comments concerning an early draft of this
report.
\par\bigskip\bigskip 

{\bf References.}
	\medskip
In order to avoid a too long list of references, we tried to restrict to the new
developments. We hope that further papers mentioned throughout our presentation
can be identified well using the appropriate chapters of this Handbook as well as
standard lists of references. But also concerning the cluster approach,
there are many more papers of interest and most are still preprints (see the arXiv).
Of special interest should be the survey by Buan and Marsh [BM].

	\medskip
\item{[ABS]} Assem, Br\"ustle and Schiffler: Cluster-tilted
algebras
  as trivial extensions. arXiv: math.RT/0601537
\item{[A]} Auslander: Comments on the functor $\Ext$. Topology 8 (1969), 151-166.
\item{[AR1]} Auslander, Reiten: On a generalized version of the Nakayama conjecture.
  Proc. Amer.Math.Soc. 52 (1975), 69-74.
\item{[AR2]}  Auslander, Reiten: Applications of contravariantly finite subcategories.
  Advances Math. 86 (1991), 111-152.
\item{[BeR]} Beligiannis, Reiten: Homological and homotopical aspects of
   torsion theories. Memoirs Amer.Math.Soc. (To appear).
\item{[BFZ]} Berenstein, Fomin, Zelevinsky:
  Cluster algebras III: Upper bounds and double Bruhat cells.
 arXiv: math.RT/0305434
\item{[BM]} Buan, Marsh: Cluster-tilting theory. Proceedings of the ICRA meeting, Mexico 2004.
  Contemporary Mathematics. (To appear).
\item{[B-T]} Buan, Marsh, Reineke, Reiten, Todorov: Tilting theory and cluster combinatorics.
   Advances Math. (To appear). arXiv: math.RT/0402054
\item{[BMRT]} Buan, Marsh, Reiten, Todorov: Clusters and seeds in acyclic cluster algebras.
  arXiv: math.RT/0510359
\item{[BMR]} Buan, Marsh, Reiten: Cluster tilted algebras. Trans. Amer. Math. Soc. (To appear). arXiv: math.RT/0402075
\item{[BMR2]} Buan, Marsh, Reiten: Cluster mutation via quiver representations.
  arXiv: math.RT/0412077.
\item{[BMR3]} Buan, Marsh, Reiten: Cluster-tilted algebras of finite representation type.
  arXiv: math.RT/0509198
\item{[BR]} Buan, Reiten: Cluster algebras associated with extended Dynkin quivers.
 arXiv: math.RT/0507113
\item{[BRS]} Buan, Reiten, Seven: Tame concealed algebras and cluster algebras of minimal
    infinite type. arXiv: math.RT/0512137
\item{[CC]} Caldero, Chapoton:
   Cluster algebras as Hall algebras of quiver representations. arXiv:
   math.RT/0410187
\item{[CK1]}Caldero, Keller: From triangulated categories to cluster algebras. I.
   arXiv: math.RT/0506018

\item{[CK2]}Caldero, Keller: From triangulated categories to cluster algebras. II.
   arXiv: math.RT/0510251.
\item{[C]} Crawley-Boevey: Exceptional sequences of representations of quivers
   in 'Representations of algebras', Proc. Ottawa 1992, eds V. Dlab and H. Lenzing,
  Canadian Math. Soc. Conf. Proc. 14 (Amer. Math. Soc., 1993), 117-124.
\item{[CKe]} Crawley-Boevey, Kerner:
  A functor between categories of regular modules for wild hereditary algebras.
  Math. Ann., 298 (1994), 481-487.

\item{[DW1]} Derksen, Weyman: Semi-invariants of quivers and saturation for Littlewood-Richardson   coefficients, J. Amer. Math. Soc. 13 (2000), no. 3, 456-479.
\item{[DW2]} Derksen, Weyman: Quiver Representations. Notices of the AMS 52, no. 2.
\item{[DW3]} Derksen, Weyman:
  On the canonical decomposition of quiver representations, Compositio Math. (To appear).
\item{[FZ1]} Fomin, Zelevinsky:
  Cluster algebras I: Foundations.
  J. Amer. Math. Soc. 15 (2002), 497-529.
\item{[FZ2]} Fomin, Zelevinsky:
  Cluster algebras II. Finite type classification. Invent. Math. 154 (2003), 63-121.

\item{[FZ3]} Fomin, Zelevinsky: Cluster algebras. Notes for the CDM-03 conference.
  arXiv: math.RT/0311493
\item{[H1]}  Happel: Relative invariants and subgeneric orbits of quivers of
finite and tame type. Spinger LNM 903 (1981), 116-124.
\item{[H2]} Happel:
  Triangulated categories in the representation theory of
 finite dimensional algebras. London Math.Soc. Lecture Note Series 119 (1988)
\item{[HRS]} Happel, Reiten, Smal{\o}: Tilting in abelian categories and quasitilted algebras.
 Memoirs Amer. Math. Soc. 575 (1995).
\item{[Hi]} Hille: The volume of a tilting module. Lecture. Bielefeld 2005, see \newline
  www.math.uni-bielefeld.de/birep/select/hille\_6.pdf.
\item{[Hg]} Huang: Generalized tilting modules with finite injective dimension.
   arXiv: math.RA/0602572
\item{[Hu]} Hubery: Acyclic cluster algebras via Ringel-Hall algebras. 
  See wwwmath.uni-paderborn.de/\~hubery/Cluster.pdf
\item{[HW]} Hughes, Waschb\"usch: Trivial extensions of tilted algebras. Proc. London Math. Soc.
   46 (1983), 347-364.
\item{[K]} Keller: On triangulated orbit categories.
   arXiv: math.RT/0503240.
\item{[KR]} Keller, Reiten: Cluster-tilted algebras are Gorenstein and stably Calabi-Yau.
   arXiv: math.RT/0512471
\item{[KT]} Kerner, Takane: Universal filtrations for modules in perpendicular categories. In
   Algebras and Modules II (ed. Reiten, Smal{\o}, Solberg). Amer.Math. Soc. 1998. 347-364.
\item{[KSX]} K\"onig, Slungard, Xi: Double centralizer properties, dominant dimension,
 and tilting modules. J. Algebra 240 (2001), 393-412.
\item{[MR]} Mantese, Reiten: Wakamatsu tilting modules. J Algebra 278 (2004), 532-552.
\item{[MRZ]} Marsh, Reineke, Zelevinsky:
   Generalized associahedra via quiver representations. Trans. Amer. Math. Soc.
   355 (2003), 4171-4186
\item{[R1]} Ringel: Reflection functors for hereditary algebras. J. London Math. Soc. (2) 21    (1980), 465-479.
\item{[R2]} Ringel: The braid group action on the set of exceptional sequences of a hereditary algebra. In: Abelian Group Theory and Related Topics. Contemp. Math. 171 (1994), 339-352.
\item{[R3]} Ringel: Exceptional objects in hereditary categories.
Proceedings Constantza Conference. An. St. Univ. Ovidius Constantza Vol. 4 (1996), f. 2, 150-158.
\item{[R4]} Ringel: Exceptional modules are tree modules.
   Lin. Alg. Appl. 275-276 (1998) 471-493.
\item{[R5]} Ringel: Combinatorial Representation Theory: History and Future.
  In: Representations of Algebras. Vol. I (ed. D.Happel, Y.B.Zhang). BNU Press. (2002) 122-144.
\item{[Ru]} Rudakov: Helices and vector bundles. London Math. Soc. Lecture Note Series. 148.
\item{[SK]} Sato, Kimura: A classification of irreducible prehomogeneous
 vector spaces and their relative invariants. Nagoya Math. J. 65 (1977), 1-155.
\item{[Sc]} Schofield: Semi-invariants of quivers. J. London Math. Soc. 43 (1991), 383-395.
\item{[S]} Seven: Recognizing Cluster Algebras of Finite type.
  arXiv: math.CO/0406545
\item{[TW]} Tachikawa, Wakamatsu: Tilting functors and stable equivalences
 for selfinjective algebras. J. Algebra (1987).
\item{[Z1]} Zelevinsky: From Littlewood-Richardson coefficients to cluster algebras
   in three lectures. In: Fomin (ed): Symmetric Functions 2001. 253-273.
\item{[Z2]} Zelevinsky: Cluster algebras: Notes for the 2004 IMCC.
     arXiv: math.RT/0407414.
\item{[Zh1]} Zhu, Bin: Equivalences between cluster categories.
     arXiv: math.RT/0511382
\item{[Zh2]} Zhu, Bin: Applications of BGP-reflection functors: isomorphisms of
        cluster algebras. arXiv: math.RT/0511384
	\bigskip\bigskip
\noindent
Claus Michael Ringel, Fakult\"at f\"ur Mathematik, Universit\"at Bielefeld\par
\noindent
POBox 100 131, D-33501 Bielefeld, Germany \par
\noindent
email: {\tt ringel\@math.uni-bielefeld.de}

\bye